\documentclass{compositio}
\usepackage[francais]{babel}
\usepackage[T1]{fontenc}
\usepackage[utf8]{inputenc}
\usepackage{url}
\usepackage{amsmath}
\usepackage{amsfonts}
\usepackage{amssymb}
\usepackage{enumerate}
\newtheorem{thm}{Th\'eor\`eme}[section]
\newtheorem{prop}[thm]{Proposition}
\newtheorem{lemme}[thm]{Lemme}
\newtheorem{cor}[thm]{Corollaire}
\theoremstyle{definition}
\newtheorem{question}[thm]{Question}
\newtheorem{hyp}[thm]{Hypoth\`ese}

\newtheorem{defi}[thm]{D\'efinition}
\theoremstyle{remark}
\newtheorem{nota}[thm]{Notation}
\newtheorem{notas}[thm]{Notations}

\newtheorem{rem}[thm]{Remarque}

\makeatletter
\@addtoreset{equation}{section}
\renewcommand{\theequation}{\arabic{section}.\arabic{equation}}
\makeatother
\newcommand{\termin}[1]{{\em #1}}
\newcommand{\A}{\mathbf{A}}
\newcommand{\F}{\mathbf{F}}
\newcommand{\G}{\mathbf{G}}
\newcommand{\LL}{\mathbf{L}}
\newcommand{\N}{\mathbf{N}}
\newcommand{\PP}{\mathbf{P}}
\newcommand{\Q}{\mathbf{Q}}
\newcommand{\R}{\mathbf{R}}
\newcommand{\Z}{\mathbf{Z}}
\def\bT{\boldsymbol{T}}
\newcommand{\bd}{\boldsymbol{d}}
\newcommand{\dd} {\boldsymbol{d}}
\newcommand{\be}{\boldsymbol{e}}
\newcommand{\ff} {\boldsymbol{f}}
\newcommand{\bn}{\boldsymbol{n}}
\newcommand{\mm} {\boldsymbol{m}}
\newcommand{\nn} {\boldsymbol{n}}
\newcommand{\xx}{\boldsymbol{x}}
\newcommand{\yy}{\boldsymbol{y}}
\newcommand{\cA}{{\cal A}}
\newcommand{\cC}{{\cal C}}
\newcommand{\cD}{{\cal D}}
\newcommand{\cM}{{\cal M}}
\newcommand{\cT}{{\cal T}}
\DeclareMathAlphabet{\eulercal}{U}{eus}{m}{n}
\newcommand{\ecA}{{\eulercal A}}
\newcommand{\ecC}{{\eulercal C}}
\newcommand{\ecF}{{\eulercal F}}
\newcommand{\ecH}{{\eulercal H}}
\newcommand{\ecL}{{\eulercal L}}
\newcommand{\ecM}{{\eulercal M}}
\newcommand{\ecO}{{\eulercal O}}
\newcommand{\ecT}{{\eulercal T}}
\newcommand{\ecX}{{\eulercal X}}
\newcommand{\eq} {\Leftrightarrow}
\newcommand{\longeq} {\Longleftrightarrow}
\newcommand{\longto} {\longrightarrow}
\newcommand{\longlto}{\longleftarrow}
\newcommand{\isom} {\overset{\sim}{\to}}
\newcommand{\longisom} {\overset{\sim}{\longto}}
\DeclareMathOperator{\Pic}{Pic}
\DeclareMathOperator{\Div}{Div}
\DeclareMathOperator{\rg}{rg}
\DeclareMathOperator{\Hom}{Hom}
\DeclareMathOperator{\HOM}{\textbf{Hom}}
\DeclareMathOperator{\Min}{Min}
\DeclareMathOperator{\Supp}{Supp}
\DeclareMathOperator{\Spec}{Spec}
\DeclareMathOperator{\Res}{Res}
\DeclareMathOperator{\Tr}{Tr}
\DeclareMathOperator{\Frob}{Frob}
\DeclareMathOperator{\NS}{NS}
\let\leq\leqslant
\let\geq\geqslant
\newcommand{\cf}{{\it cf.}\ }
\newcommand{\ie}{{\it i.e.}\ }
\newcommand{\cone}{\cC}
\newcommand{\sumu}[1]{\underset{#1}{\sum}}
\newcommand{\produ}[1]{\underset{#1}{\prod}}
\newcommand{\cupu}[1]{\underset{#1}{\cup}}
\newcommand{\capu}[1]{\underset{#1}{\cap}}
\newcommand{\coprodu}[1]{\underset{#1}{\coprod}}
\newcommand{\sqcupu}[1]{\underset{#1}{\sqcup}}
\newcommand{\oplusu}[1]{\underset{#1}{\oplus}}
\newcommand{\otimesu}[1]{\underset{#1}{\otimes}}
\newcommand{\symb}[1]{\left[#1\right]}
\newcommand{\vide} {\varnothing}
\newcommand{\eqdef} {\overset{\text{{\tiny{d\'ef}}}}{=}}
\newcommand{\grsym} {\mathfrak{S}} 
\newcommand{\ind} {\mathbf{1}}
\newcommand{\acc}[2]{\left\langle #1\, ,\,#2 \right\rangle} 
\newcommand{\chif}{\chi_{{}_{\text{\textnormal{\scriptsize form}}}}}
\newcommand{\abs}[1]{\left| #1 \right|} 
\newcommand{\place}[2]{\raisetag{#1}\tag*{$#2$}}
\newcommand{\card}[1]{\left|#1\right|} 
\newcommand{\bigast} {\boldsymbol{(\ast)}}
\newcommand{\cc}{\ecC}
\newcommand{\mk}{\cM_k}
\newcommand{\mkq} {\cM_{k,\Q}}
\newcommand{\mkc} {\cM_k^{\chi}}
\newcommand{\mkcq} {\cM_{k,\Q}^{\chi}}
\newcommand{\mloc}{\cM_{k,\text{loc}}}
\newcommand{\mlocc} {\cM^{\chi}_{k,\text{loc}}}
\newcommand{\mlocq }{\cM_{k,\text{loc},\Q}}
\newcommand{\mloccq} {\cM^{\chi}_{k,\text{loc},\Q}}
\newcommand{\mhat} {\widehat{\cM_k}}
\newcommand{\mhatq}{\widehat{\cM_{k,\Q}}}
\newcommand{\mhatc} {\widehat{\cM^{\chi}_k}}
\newcommand{\mhatcq} {\widehat{\cM^{\chi}_{k,\Q}}}
\newcommand{\muxm} {\mu^{\text{mot}}_{X}}
\newcommand{\mubsm} {\mu^{\bs,\text{mot}}_{X}}
\newcommand{\mubsmp} {\mu^{\bs,\text{mot}}_{\PP^1}}
\newcommand{\mubsc} {\mu^{\bs,\chi}_{X}}
\newcommand{\mubsp} {\mu^{\bs,\chi}_{\PP^1}}
\newcommand{\musc} {\mu^{\chi}_{\Sigma}}
\newcommand{\mubmx} {\mu^{B,\text{mot}}_{X}}
\newcommand{\mubcx} {\mu^{B,\chi}_{X}}
\newcommand{\zetamh} {Z^{\text{mot}}_{h}}
\newcommand{\zetamthz} {Z^{\text{mot}}_{\PP^1,U,h_0}}
\newcommand{\zetacthz} {Z^{\chi}_{\PP^1,U,h_0}}
\newcommand{\zetam}[1]{Z^{\text{mot}}_{#1}}
\newcommand{\zetamc}[1]{Z^{\chi}_{#1}}
\newcommand{\zetaxm} {\zetam{X}}
\newcommand{\zetacm} {\zetam{\cc}}
\newcommand{\zetapm} {\zetam{\PP^1}}
\newcommand{\zetaxmc} {\zetamc{X}}
\newcommand{\zetacmc} {\zetamc{\cc}}
\newcommand{\fil}{\ecF}
\newcommand{\filc} {\fil_{\chi}}
\newcommand{\xs} {X_{\Sigma}}
\newcommand{\bs} {B_{\Sigma}}
\newcommand{\ts} {\cT_{\Sigma}}
\newcommand{\roa} {\rho_{\alpha}}
\newcommand{\da} {d_{\alpha}}
\newcommand{\ea} {e_{\alpha}}
\newcommand{\na} {n_{\alpha}}
\newcommand{\ha} {h_{\alpha}}
\newcommand{\Ea} {E_{\alpha}}
\newcommand{\cDa} {\cD_{\alpha}}
\newcommand{\ecLa} {\ecL_{\alpha}}
\newcommand{\ecMa} {\ecM_{\alpha}}
\newcommand{\hirz}[1]{\ecH_{#1}}
\newcommand{\Gammaf}{\Gamma_{\!\!\ff}}
\newcommand{\cara}[1]{\ecX^{{}^{\,\ast}}\left(#1\right)}
\newcommand{\cocara}[1]{\ecX_{\ast}\left(#1\right)}
\newcommand{\ceff} {C_{\text{eff}}}
\newcommand{\xeff} {X^{0,+}}
\newcommand{\kpff} {K_0(\text{PFF}_k)}
\newcommand{\wt} {\widetilde}
\begin{document}
\title[Produit eul\'erien motivique]{Produit eul\'erien motivique et courbes rationnelles sur les vari\'et\'es
toriques \\
Motivic Euler product and rational curves on toric varieties}
\author{David Bourqui}
\email{david.bourqui@univ-rennes1.fr}
\address{IRMAR \\ Universit\'e de Rennes 1 \\ Campus de Beaulieu \\
35042 Rennes cedex \\ France}
\begin{abstract}
We study the asymptotical behaviour of the moduli space 
of morphisms of given anticanonical degree from
a rational curve to a split toric variety, when the
degree goes to infinity. We obtain in this case a geometric analogue of Manin's conjecture
about rational points of bounded height on varieties defined over a global field.
The study is led through a generating series 
whose coefficients lie in a Grothendieck ring of motives, the motivic height zeta function.
In order to establish convergence properties of this function, we use a notion of eulerian motivic product. 
It relies on a construction of Denef and Loeser which associates a virtual motive to 
a first order logic ring formula.
\\~\\
\indent Nous \'etudions le comportement asymptotique de l'espace des modules des morphismes
de degr\'e anticanonique donn\'e d'une courbe rationelle vers une vari\'et\'e torique d\'eploy\'ee,
lorsque ce degr\'e tend vers l'infini. 
Nous obtenons dans ce cas un analogue g\'eom\'etrique de la conjecture 
de Manin sur le nombre de points de hauteur born\'ee des vari\'et\'es d\'efinies sur un corps global.
L'\'etude se fait via une s\'erie g\'en\'eratrice \`a coefficients dans un anneau de Grothendieck de motifs,
la fonction z\^eta des hauteurs motivique. Afin d'\'etablir des propri\'et\'es de convergence de cette
fonction, nous utilisons une notion de produit eul\'erien motivique, laquelle repose
sur la construction de Denef et Loeser permettant d'associer un motif
virtuel \`a une formule  logique du premier ordre
dans le langage des anneaux.

\end{abstract}
\maketitle

\section{Introduction}\label{sec:intro}
Soit $k$ un corps, $\cc$ une courbe projective, lisse et g\'eom\'etriquement int\`egre d\'efinie sur $k$, 
et $V$ une vari\'et\'e projective et lisse d\'efinie sur $k$. 
On fixe un faisceau $\ecL$ sur $V$ dont la classe dans le groupe de N\'eron-Severi est situ\'ee
\`a l'int\'erieur du c\^one effectif. 

Si le corps $k$ est fini, 
un probl\`eme naturel est d'\'etudier le comportement asymptotique du nombre de morphismes de $\cc$ vers $V$
de $\ecL$-degr\'e donn\'e quand ce degr\'e tend vers l'infini. 
Ce probl\`eme est l'analogue g\'eom\'etrique 
du probl\`eme arithm\'etique du comptage asymptotique du nombre de points de hauteur
born\'ee sur une vari\'et\'e d\'efinie sur un corps de nombres. 
Concernant ces deux probl\`emes, une s\'erie de questions a \'et\'e soulev\'ee par Manin et ses collaborateurs 
vers la fin des ann\'ees 1980, 
lesquelles  
ont depuis \'et\'e \'etudi\'ees pour de larges classes de vari\'et\'es, notamment dans le cas arithm\'etique. 
Le lecteur pourra se reporter \`a \cite{Pey:bki} et \cite{Pey:bordeaux} pour plus de pr\'ecisions et un \'etat des lieux
sur la question en 2001, ainsi qu'à \cite{Bro:manin_dim_2} pour une
description de progrès plus récents dans le cas des surfaces.

Si le corps $k$ est quelconque, on peut plus g\'en\'eralement s'int\'eresser au comportement asymptotique
de la vari\'et\'e param\'etrant 
les morphismes de $\cc$ vers $V$ de $\ecL$-degr\'e donn\'e quand ce degr\'e tend vers l'infini. 
On peut par exemple essayer d'estimer le comportement asymptotique de la dimension 
et du nombre de composantes
g\'eom\'etriques irr\'eductibles de ces espaces de modules.
Une autre fa\c con de concevoir le problème est d'\'etudier une s\'erie g\'en\'eratrice associ\'ee qui est 
\`a coefficients dans l'anneau de Grothendieck des vari\'et\'es (ou des motifs) sur $k$, et qui, lorsque
le corps de base est fini, se sp\'ecialise sur la fonction z\^eta des hauteurs classiques. 
Nous renvoyons
\`a la sous-section \ref{subsec:manin_mot} pour une formulation plus pr\'ecise des questions qu'il semble
l\'egitime de se poser dans ce cas de figure. 
Signalons que nombre de ces questions sont d\^ues \`a Peyre.

Dans ce texte, nous étudions ces questions pour les variétés toriques déployées.
Les principaux résultats obtenus  sont rassemblés dans l'énoncé suivant. 
\begin{thm}\label{thm:theoprinc}
Soit $k$ un corps et $V$ une variété torique déployée sur $k$,
supposée projective et  lisse.
Soit $U$ son orbite ouverte. 
Pour tout entier $d\geq 1$, on note $U_{0,d}$ la variété
quasi-projective paramétrant les $k$-morphisme $\PP^1_k\to V$ dont
l'image rencontre $U$ et de degré anticanonique $d$.
\begin{enumerate}
\item\label{item:1:theoprinc}
Soit $m\geq 1$ un entier. On suppose que $V$ est la $m$-ème surface de Hirzebruch.
Alors la série
$$
(1+\LL\,T)\,(1+\LL\,T+\LL^2\,T^2+\dots+\LL^{\,m+1}T^{\,m+1})\,(1-\LL\,T)^{2}\,\Big(\sum_{d\geq 1} \symb{U_{0,d}}\,T^d\Big)
$$
est un polynôme à coefficients dans l'anneau de Grothendieck des
$k$-variétés, 
dont la valeur en $\LL^{-1}$ est 
$
\LL^{\,2}\,(1-\LL^{-2})^2
$.
\item\label{item:2:theoprinc}
On suppose le corps $k$ de caract\'eristique z\'ero.
La série 
$$
(1-\LL\,T)^{\rg(\Pic(V))}\Big(\sum_{d\geq 1} \chi(\symb{U_{\ecL_0,d}})\,T^d\Big)
$$
(à coefficients dans l'anneau des motifs virtuels)
converge en $T=\LL^{-1}$ vers
\begin{equation}\label{eq:tammot}
\alpha^{\ast}(V)
\,
\LL^{\,\dim(V)}\,
\left(
\frac{1}{1-\LL^{-1}}
\right)^{\,\rg(\Pic(V))}
\exp
\left(
\sum_{n\geq 1}
\Psi^{\chi}_n(\PP^1)
\log
\left(
(1-\LL^{\,-n})^{\rg(\Pic(V))}
\frac
{\Phi^{\chi}_n(V)}
{\LL^{\,-n\,\dim(V)}}
\right)
\right)
\end{equation}
\end{enumerate}
\end{thm}
Précisons les notation utilisées (\cf la sous-section
\ref{subsec:groth}). 
On désigne par $\symb{X}$  la classe d'une $k$-variété $X$
dans l'anneau de Grothendieck des variétés et, si $k$ est de
caractéristique zéro, par $\chi(\symb{X})$ son image dans l'anneau de Grothendieck des motifs de Chow.
Le symbole $\LL$ désigne indifféremment $\symb{\A^1}$ ou $\chi(\symb{\A^1})$.
La convergence s'entend au sens de la topologie définie par la filtration
dimensionnelle, employée initialement dans la théorie de l'intégration
motivique.
Les familles de motifs virtuels $(\Psi^{\chi}_n(\PP^1))_{n\geq 1}$ et
$(\Phi^{\chi}_n(V))_{n\geq 1}$ sont définies à la sous-section
\ref{subsec:seconde} ; on peut les voir 
comme des incarnations motiviques des
notions de nombre de points fermés de degré $n$ et de nombre de
points rationnel à valeurs dans une extension de degré $n$ d'une
variété sur un corps fini. De cette façon \eqref{eq:tammot} peut
s'interpréter comme un analogue motivique de la constante de Peyre
intervenant (au moins conjecturalement) dans l'expression
asymptotique du nombre de points de hauteur bornée sur les variétés de Fano.
L'invariant $\alpha^{\ast}(V)$ apparaissant dans \eqref{eq:tammot} est
défini à la sous-section \ref{subsec:manin_mot}.

Signalons que Peyre a d\'emontr\'e un r\'esultat similaire au théorème
\ref{thm:theoprinc} lorsque $V$
est une vari\'et\'e de drapeaux, la courbe $\ecC$ étant de genre quelconque (\cf  \cite{Pey:ober}).

Les ingr\'edients de la d\'emonstration du théorème \ref{thm:theoprinc}
sont des versions motiviques de ceux que nous avons utilis\'es dans \cite{Bou:vtetor} pour calculer
la fonction z\^eta des hauteurs d'une vari\'et\'e torique d\'eploy\'ee d\'efinie sur un corps global
de caract\'eristique non nulle. Ce sont :
\begin{enumerate}
\item
le lemme \ref{lm:param_mot}, qui explicite
la vari\'et\'e param\'etrant les morphismes d'une courbe rationnelle vers une vari\'et\'e torique d\'eploy\'ee
\`a l'aide de la description de Cox du foncteur des points d'une telle vari\'et\'e (\cf  \cite{Cox:funct}).
Ce lemme est une version g\'eom\'etrique du lemme 2 de \cite{Bou:vtetor},
lui-m\^eme inspir\'e de la m\'ethode utilis\'ee par Salberger 
sur les corps de nombres dans \cite{Sal:tammes}. Dans tout ceci, l'utilisation
du torseur universel au-dessus d'une variété torique 
déployée joue un rôle essentiel.
\item
une formule d'inversion de M\"obius motivique, 
version motivique de la formule d'inversion utilis\'ee dans \cite{Bou:vtetor},
elle-m\^eme adapt\'ee des formules d'inversion utilis\'ees par Peyre et Salberger
dans le cadre de la version arithmétique des conjecture de Manin.
\item
une notion de \og produit eul\'erien motivique\fg\
qui nous permet de d\'emontrer des propri\'et\'es 
de convergence de la s\'erie g\'en\'eratrice associ\'ee \`a la formule d'inversion en question,
et de donner une interpr\'etation du terme principal de la fonction z\^eta similaire
\`a l'interpr\'etation en termes de nombre de Tamagawa dans le cas classique.
Nous faisons ici usage de la construction de Denef et Loeser permettant d'associer
canoniquement un motif virtuel \`a une formule logique du premier ordre.
\end{enumerate}

Nous d\'ecrivons \`a pr\'esent l'organisation de l'article.

Dans la section \ref{sec:prod_eul_mot}, apr\`es quelques rappels, 
nous présentons la notion de produit eulerien motivique et démontrons notamment que la fonction z\^eta de Hasse-Weil
motivique s'\'ecrit sous forme d'un produit eul\'erien motivique.

Dans la section \ref{sec:fonc_mob_mot}, nous introduisons des fonctions d'inversions de M\"obius motiviques
et montrons que les s\'eries g\'en\'eratrices associ\'ees s'\'ecrivent sous forme d'un produit eul\'erien motivique,
ce qui permet d'en d\'egager des propri\'et\'es de convergence.

Dans la section \ref{sec:fonc_haut_mot}, nous définissons la fonction
zêta des hauteurs motivique et précisons quelques questions 
permettant d'esquisser une version motivique des conjectures de Manin.

Enfin, dans la section \ref{sec:calcul}, 
nous d\'ecrivons la vari\'et\'e des morphismes de degr\'e donn\'e de
$\PP^1$ vers une vari\'et\'e torique déployée. 
Utilisant une fonction de M\"obius ad\'equate et les r\'esultats de la section \ref{sec:fonc_mob_mot}, nous en d\'eduisons
la démonstration du théorème \ref{thm:theoprinc}. Ceci montre que
certaines des questions de la section \ref{sec:fonc_haut_mot} ont une
réponse positive dans le cas d'une variété torique déployée.

{\bf Remerciements}

Je remercie Emmanuel Peyre et Antoine Chambert-Loir pour d'utiles discussions.
Je remercie Fran\c cois Loeser de m'avoir indiqu\'e la r\'ef\'erence \cite{GusLueMel:power}.

\section{Fonction z\^eta de Hasse-Weil et produit eul\'erien motivique}\label{sec:prod_eul_mot}

\subsection{Quelques rappels et d\'efinitions}

\subsubsection{Anneaux de Grothendieck de variétés et de motifs}\label{subsec:groth}
Soit $k$ un corps.
On note $\mk$  l'anneau de
Grothendieck de la catégorie des variétés définies sur $k$ (\cf
\cite[\S 13.1.1]{And:mot}). 
Si $X$ est une telle variété, on note
$\symb{X}$ sa classe dans $\mk$.
On note $\LL=\symb{\A^1_k}$ la classe de la droite affine et
$
\mloc=\mk\left[\LL^{-1}\right].
$
Si le corps $k$ est fini de cardinal $q$, l'application qui \`a 
un $k$-schéma $X$ de type fini associe le nombre de points $k$-rationnels
de $X$
induit un morphisme d'anneau
$
\#_k\,:\,\mloc\to\Z\left[q^{\,-1}\right].
$
On munit $\mloc$ de la filtration dimensionnelle introduite par
Kontsevitch dans le cadre de la théorie de l'intégration motivique : 
pour $m\in \Z$, $\fil^{\,m} \mloc$ désigne le sous-groupe de $\mk$ engendr\'e par 
les \'el\'ements de la forme $\LL^{-i}[V]$, 
o\`u $V$ est une $k$-variété et $i$ et $V$ v\'erifient 
$
i-\dim(V)\geq m.
$
On d\'efinit le complété associé
$$
\mhat\eqdef\underset{\longlto}{\lim} \,\,\mloc / \fil^{\,m} \mloc.
$$
On note $\mkq=\mk\otimes \Q$,
$\mlocq=\mloc\otimes \Q$,
$\fil^{\,m}\mlocq=\fil^{\,m}\mloc \otimes \Q$
et 
$
\mhatq=\underset{\longlto}{\lim} \,\,\mlocq / \fil^{\,m} \mlocq.
$ 

Soit $K_0(\text{CHMot}_k)$ l'anneau de Grothendieck
de la catégorie motifs de Chow à coefficients rationnels définis sur
$k$ (\cf  \cite[Chapitre
4 et 13.2.1]{And:mot}).
Si $M$ est  un motif, on note
$\symb{M}$ sa classe dans $K_0(\text{CHMot}_k)$.
Si $k$ est de caract\'eristique z\'ero, 
il existe un unique morphisme
$
\chi\,:\,\mk\to K_0(\text{Mot}_k)
$
tel que la classe $[X]$ d'une vari\'et\'e $X$ projective et lisse sur $k$ 
s'envoie sur la classe du motif de Chow de $X$ (\cf  \cite[Theorem 4]{GilSou:descent} 
ainsi que \cite{GuiNav:crit_ext} et \cite{Bit:univ_eul_car}). 
Nous d\'esignerons par $\mkc$ l'image de $\mk$ par ce morphisme. 
On notera $\LL$ en lieu et place de $\chi(\LL)$.
On note $\mlocc=\mkc[\LL^{-1}]$, 
$\filc^{\,\bullet}$ la filtration image de $\fil^{\,\bullet}$ par $\chi$
et $\mhatc=\underset{\longlto}{\lim} \,\,\mlocc / \filc^{\,m} \mlocc$.
On d\'efinit de mani\`ere analogue 
$\mkcq$,
$\mloccq$,
$\fil^{\,m}\mloccq$
et 
$
\mhatcq
$ 
.

Soit $X$  une $k$-variété quasi-projective. Pour tout $n\geq 1$, on
note $X^{<n>}$ la puissance symétrique $n$-ème de $X$.
Suivant Kapranov (\cf \cite{Kap:ell}), on d\'efinit
$$
\zetaxm(T)
\eqdef
\sum_{n\geq 0}
\symb{X^{<n>}}\,
T^{\,n}\quad \in \mk[[T]].
$$
Si $k$ est fini, 
 $\#_k(\zetaxm)$ 
est la fonction zêta de Hasse-Weil classique de $X$.
Pour un corps de base quelconque, $\zetaxm(T)$ est baptisée fonction zêta de Hasse-Weil motivique.
Par exemple, si $X=\PP^1$, on a pour tout $n\geq 0$
$
\left(\PP^1\right)^{<n>}\isom \PP^n
$
d'o\`u
$
\zetapm(T)=\frac{1}{(1-T)(1-\LL\,T)}.
$
En genre supérieur, on a le résultat suivant dû à Kapranov
(\cf  \cite[Theorem 1.1.9]{Kap:ell} et \cite[Theorem 3.7]{LaLu:rationality_criteria})
\begin{thm}
Soit $\cc$ une $k$-courbe projective, lisse, 
g\'eom\'etriquement int\`egre, de genre $g$, et telle 
que $\Pic^1(\cc)(k)$ soit non vide.  
Il existe alors un polyn\^ome $P_{\cc}$ \`a coefficients dans $\mk$
de degr\'e $2\,g$ tel que  
\begin{equation}\label{eq:zetacm}
(1-T)(1-\LL\,T)\zetacm(T)=P_{\cc}(T).
\end{equation}
\end{thm}

\subsubsection{Motif virtuel associé à une formule}
Concernant les rappels qui suivent, on renvoie \`a
\cite{DeLo:def_sets_motives},
\cite{DeLo:grot_pff} et \cite{Nic:rel:motive} pour plus de d\'etails.
Dans ce texte, on appelle \termin{formule \`a coefficients dans $k$} (voire \termin{formule} si le corps $k$ est clairement 
indiqu\'e par le contexte) une formule du premier ordre dans le
langage des anneaux \`a coefficients dans $k$.
Pour toute formule $\varphi$ \`a coefficients dans $k$ en $n$ variables libres et 
toute extension $K$ de $k$ on notera $\varphi(K)$ le sous-ensemble
de $K^n$ constitu\'e des \'el\'ements de $K^n$ satisfaisant $\varphi$. 
Si $X$ est une variété quasi-affine définie sur $k$, on appellera \termin{formule sur $X$}
toute formule \`a coefficients dans $k$ en $n$ variables libres de la forme $\varphi \wedge \varphi_X$
où $\varphi$ est une formule en $n$ variables libres et $\varphi_X$ une formule d\'efinissant les équations
d'un plongement de $X$ dans l'espace affine $\A^n$.

Un \termin{corps pseudo-fini} est un corps parfait,
pseudo-algébriquement clos et admettant dans 
une clôture algébrique fixée
une unique extension de degré $n$ pour tout $n\geq 1$.

Soit $d\geq 1$ et  $\varphi$, $\psi$ des formules à coefficients dans $k$
en les variables libres $(x_1,\dots,x_m)$ et $(y_1,\dots,y_n)$ respectivement.
On dit que \termin{$\varphi$ est un $d$-revêtement de $\psi$}
s'il existe une formule 
$\theta$ en les variables libres $(x_1,\dots,x_m,y_1,\dots,y_n)$ telle
que pour tout corps pseudo-fini $K$  contenant $k$, l'ensemble
$\theta(K)\subset K^n\times K^m$
est le graphe d'une application $d$ pour $1$ de $\varphi(K)$
sur $\psi(K)$. Deux formules sont dites logiquement équivalentes si
l'une est un $1$-revêtement de l'autre.

On note $\kpff$ l'anneau de Grothendieck de la th\'eorie des corps
pseudo-finis sur $k$. Son groupe sous-jacent est engendr\'e par les symboles $\symb{\varphi}$, o\`u $\varphi$ est une formule 
\`a coefficients dans $k$. Ces générateurs satisfont les relations $\symb{\varphi}=\symb{\psi}$ si $\varphi$
et $\psi$ sont logiquement \'equivalentes et 
$
\symb{\varphi\vee\psi}+\symb{\varphi\wedge \psi}
=
\symb{\varphi}+\symb{\psi}
$
si $\varphi$ et $\psi$ ont les m\^emes variables libres.
Le produit est défini par $\symb{\varphi}\,\symb{\psi}\eqdef
\symb{\varphi\vee \psi}$ pour toutes formules $\varphi$ et $\psi$
ayant des ensembles de variables libres disjoints.

\begin{thm}\label{thm:denefloeser}
Soit $k$ un corps de caractéristique zéro.
Il existe un unique morphisme d'anneaux
\begin{equation}
\chif\,:\,\kpff\longto \mkcq
\end{equation}
qui envoie la classe d'une formule qui est une conjonction d'\'equations polyn\^omiales
sur la classe de la vari\'et\'e affine d\'efinie par ces \'equations 
et qui satisfait pour toutes formules $\varphi$ et $\psi$ telles que $\varphi$
est un $d$-revêtement de $\psi$ la relation
\begin{equation}\label{rel:d:rev}
\chif\left(\symb{\varphi}\right)=d\,\chif\left(\symb{\psi}\right).
\end{equation}
\end{thm}
\begin{rem}
Soient $X$ et $Y$ des vari\'et\'es affines normales irr\'eductibles
et $X\to Y$ un revêtement galoisien étale de groupe $G$.
Pour tout sous-groupe cyclique $C$ on note $\varphi_{X,Y,C}$ une formule
sur $Y$ telle que, pour tout corps pseudo-fini $K$ contenant $k$, $\varphi_{X,Y,C}(K)$
s'identifie \`a l'ensemble des \'el\'ements de $Y(K)$ 
qui se rel\`event \`a un \'el\'ement
de $(X/C)(K)$ mais pas \`a un \'el\'ement de $(X/D)(K)$ pour tout
sous-groupe strict $D$ de $C$,
en d'autre termes 
qui admettent $C$ comme groupe de d\'ecomposition dans le revêtement $X\to Y$. 
Une telle formule est appelée \termin{formule galoisienne}.
La relation \eqref{rel:d:rev} entraîne alors la relation 
\begin{equation}\label{eq:rel_form_gal}
\chif\left(\symb{\varphi_{X,Y,C}}\right)
=
\frac{
\card{C}
}
{
\card{N_{G}(C)}
}
\chif\left(\symb{\varphi_{X,X/C,C}}\right).
\end{equation}
Denef et Loeser ont démontré l'existence et l'unicité d'un morphisme $\chif$
vérifiant la relation \eqref{eq:rel_form_gal} pour toute formule
galoisienne (\cf \cite[Theorem 2.1]{DeLo:grot_pff}). 
Le fait qu'un tel morphisme vérifie en outre la condition  \eqref{rel:d:rev} est
énoncé sans preuve dans \cite{Ha:whatis}. Cette propriété est démontrée (et étendue à
un cadre relatif) par Nicaise dans \cite{Nic:rel:motive}
(\cf notamment le lemme 8.5). Pour une démonstration élémentaire du
fait que la relation \eqref{eq:rel_form_gal} entraîne la relation
\eqref{rel:d:rev}, on peut consulter \cite{Bou:DL}.
\end{rem}

\subsection{Produit eul\'erien motivique : premi\`ere approche}

On cherche un analogue motivique de la décomposition de la fonction
zêta de Hasse-Weil classique en produit eulérien.
Soit $k$ un corps. On définit pour toute $k$-variété quasi-projective $X$
une famille $(\Phi_n(X))_{n\geq 1}$
d'éléments de $\mk$ par la relation
\begin{equation}\label{eq:prod_eul_mot_1}
\sum_{n\geq 1}\Phi_n(X)\,T^n
=
T\,\frac{d}{dT} \log \zetaxm(T)
\end{equation}
et une famille $(\Psi_n(X))_{n\geq 1}$  d'éléments de $\mkq$
par les relations
\begin{equation}\label{eq:prod_eul_mot_2}
\forall n\geq 1, 
\quad
\Phi_n(X)=\sum_{d|n}d\,\Psi_d(X).
\end{equation}
\begin{lemme}\label{lm:prop:phin:psin}
Soit $k$ un corps et $X$ une $k$-variété quasi-projective.
\begin{enumerate}
\item\label{item:0:lm:prop:phin:psin}
On suppose $k$ fini. Pout tout $n\geq 1$, $\#_k \Phi_n(X)$
(respectivement $\#_k \Psi_n(X)$) est le nombre de points de $X$ à
valeurs dans une extension de degré $n$ de $k$ (respectivement le
nombre de points fermés de degré $n$ de $X$).
\item\label{item:1:lm:prop:phin:psin}
On a  $\Phi_1(X)=\Psi_1(X)=[X]$.
\item\label{item:2:lm:prop:phin:psin}
Pour tout $n\geq 1$, on a la relation 
\begin{equation}\label{eq:form_phi_n}
\Phi_n(X)=
\sum_{k=1}^n
(-1)^{k+1}
\,
\frac{n}{k}
\,
\sum_
{\substack{
(m_1,\dots,m_k)\in (\N_{>0})^k
\\ ~\\
m_1+\dots+m_k=n
}}
\,
\prod_{i=1}^k
\left[X^{<m_i>}\right].
\end{equation}
\item\label{lm:psi_in_fil}\label{item:3:lm:prop:phin:psin}
Pour tout $n\geq 1$, $\Phi_n(X)$ et $\Psi_n(X)$ appartiennent \`a $\fil^{-n\,\dim(X)}\mlocq$.
\item\label{lm:phi(ad)}\label{item:4:lm:prop:phin:psin}
Pour tout $d\geq 1$ et tout $n\geq 1$, on a
$
\Phi_n(\A^d)=\LL^{n\,d}.
$
\item\label{item:5:lm:prop:phin:psin}
On a la relation
\begin{align}\label{eq:prod_eul_mot_3}
\zetaxm(T)
=
\exp\left(\sum_{n\geq 1}\Psi_n(X)\,\log\left(\frac{1}{1-T^n}\right)\right).
\end{align}
\end{enumerate}
\end{lemme}
\begin{proof}
Le fait que $\#_k \zetaxm$ co\"\i ncide avec la fonction zêta de
Hasse-Weil classique et les propriétés standards d'icelle montrent le point
 \ref{item:0:lm:prop:phin:psin}.
Le point \ref{item:1:lm:prop:phin:psin} découle immédiatement des
définitions et le point \ref{item:2:lm:prop:phin:psin} d'un calcul
élémentaire. Le point \ref{item:3:lm:prop:phin:psin} se déduit du
point \ref{item:2:lm:prop:phin:psin} et des relations \eqref{eq:prod_eul_mot_2}.
Le point \ref{item:4:lm:prop:phin:psin} découle du fait qu'on a pour  tout $n\geq 1$
la relation $\symb{\left(\A^d\right)^{<n>}}=\LL^{n\,d}$ (\cf
\cite[Lemma 4.4]{Got:hilbscheme}). 
Le point \ref{item:5:lm:prop:phin:psin} découle des définitions par un calcul standard.
\end{proof}

La relation \eqref{eq:prod_eul_mot_3} peut être vue comme une décomposition
en \og produit eulérien motivique\fg. On peut généraliser ainsi cette 
notion~: 
si $P$ est un \'el\'ement de $\mkq[[T]]$ v\'erifiant $P(0)=1$
on d\'efinit le produit eul\'erien motivique associ\'e comme \'etant
$$
\Pi^{\text{mot}}_{X,P}(T)
\eqdef
\exp\left[\sum_{n\geq 1}\Psi_n(X)\,\log(P(T^n))\right]\quad\in \mkq[[T]].
$$

\subsection{Produit eul\'erien motivique : seconde approche}\label{subsec:seconde}

L'approche de la section pr\'ec\'edente 
nous semble limit\'ee 
d\`es qu'il s'agit de montrer 
que
d'autres s\'eries que la fonction z\^eta de Hasse-Weil motivique 
(telles que celles \'etudi\'ees \`a la section \ref{sec:fonc_mob_mot})
s'\'ecrivent sous forme d'un produit eul\'erien motivique.
Si $k$ est de caract\'eristique z\'ero, 
on peut utiliser le théorème  \ref{thm:denefloeser} 
pour donner une d\'efinition naturelle
de la familles  $(\Psi_n(X))$
{\em en tant qu'\'el\'ements de $\mkcq$}. Compte tenu de nos
objectifs, cette définition  s'avèrera beaucoup plus
maniable.

\subsubsection{Construction}

\begin{notas}\label{notas:x0n}
Soit $k$ un corps. Soit $n\geq 1$ un entier. Pour toute $k$-vari\'et\'e quasi-projective $X$, 
on note $(X^n)_0$ l'ouvert de $X^n$ constitu\'e des $n$-uplets d'\'el\'ements
deux \`a deux distincts, et $(X^{<n>})_0$ l'ouvert de $X^{<n>}$ image de $(X^n)_0$ 
par le morphisme 
naturel $X^n\to X^{<n>}$. 
Ce dernier morphisme induit un rev\^etement 
galoisien étale $(X^n)_0\to \left(X^{<n>}\right)_0$
de groupe $\grsym_n$.
\end{notas}
Soit $k$ un corps, $X$ une $k$-variété affine et irr\'eductible et  $n\geq 1$. 
On note 
$\pi_n$ le $\grsym_n$-rev\^etement étale $(X^n)_0\to \left(X^{<n>}\right)_0$,
$\sigma_n$ un $n$-cycle de $\grsym_n$ et 
$\psi_n(X)$ une formule sur $X^{<n>}$ telle que, 
pour tout corps pseudo-fini 
$K$ contenant $k$,
$\psi_n(X)(K)$
est l'ensemble des \'el\'ements de $X^{<n>}(K)$ qui sont dans $\left(X^{<n>}\right)_0(K)$ 
et admettent un groupe
de d\'ecomposition dans $(X^n)_0\to \left(X^{<n>}\right)_0$ engendr\'e par $\sigma_n$. 
Ce dernier ensemble
s'identifie naturellement 
\`a l'ensemble des points ferm\'es de 
degr\'e $n$ sur $X_K$.
L'expression de la formule $\psi_n(X)$ d\'epend du choix
du plongement de $X$ dans un espace affine, mais 
l'existence d'un $k$-isomorphisme entre deux tels plongements
montre que $\overline{\psi}_n(X)\eqdef\symb{\psi_n(X)}$ 
ne d\'epend que de la classe d'isomorphisme
de $X$.
\begin{lemme}\label{lm:rel_psi_UF_aff}
Soit 
$U$ un ouvert affine de $X$ et $F=X\setminus U$.
On a, pour tout $n\geq 1$,
$$
\overline{\psi}_n(X)=\overline{\psi}_n(F)+\overline{\psi}_n(U).
$$
\end{lemme}
\begin{proof}[D\'emonstration]
On note $\pi_U$ (respectivement $\pi_F$) le morphisme naturel
$U^{<n>}\to X^{<n>}$ (respectivement $F^{<n>}\to X^{<n>}$).
La formule $\psi_n(X)$ s'\'ecrit alors $\psi^{'}\vee \psi^{''}$,
o\`u $\psi^{'}$ (respectivement $\psi^{''}$)
est une formule dont l'interpr\'etation
dans un corps pseudo-fini 
$K$ contenant $k$
est l'ensemble des \'el\'ements de $X^{<n>}(K)$ qui satisfont
$\psi_n(X)$ et sont dans l'image de $\pi_U$
(respectivement $\pi_F$).
De tels \'el\'ements sont en bijection avec l'ensemble
des \'el\'ements de $U^{<n>}(K)$ satisfaisant $\psi_n(U)$.
Cette bijection est donn\'ee par le $k$-morphisme de vari\'et\'es
affines $\pi_U$, et donc $\psi^{'}$ et $\psi_n(U)$ sont logiquement
\'equivalentes. De m\^eme $\psi^{''}$ et $\psi_n(F)$ sont logiquement
\'equivalentes, d'o\`u le r\'esultat.
\end{proof}
Soit $X$ une $k$-variété quelconque 
et $X=\cupu{i\in I}X_i$ un recouvrement ouvert affine.
On pose
$$
\overline{\psi}_n(X)
=
\sum_{\vide\neq J \subset I}
\overline{\psi}_n
\Big(
\bigcap_{i\in J}X_i\,
\setminus 
\bigcup_{i\notin J}X_i
\Big),
$$
ce qui, d'après le lemme \ref{lm:rel_psi_UF_aff}, ne d\'epend pas du recouvrement choisi. 
De ce même lemme, on déduit aussitôt le résultat suivant.
\begin{lemme}\label{lm:rel_psi_UF}
Soit $X$ une $k$-vari\'et\'e, 
$U$ un ouvert de $X$ et $F=X\setminus U$.
On a pour tout $n\geq 1$
$$
\overline{\psi}_n(X)=\overline{\psi}_n(U)+\overline{\psi}_n(F).
$$
\end{lemme}
\begin{cor}\label{cor:psin_constr}
Soit $n\geq 1$.
L'application qui \`a $X$ associe $\overline{\psi}_n(X)$
s'\'etend en un morphisme de groupes
$$
\overline{\psi}_n\,:\,\mk\longto \kpff.
$$
En particulier, 
si $A$ est un sous-ensemble constructible d'une $k$-vari\'et\'e, 
$\overline{\psi}_n(A)$ est bien d\'efini.
Si $(A_i)_{i\in I}$ est une famille finie de sous-ensembles
constructibles d'une vari\'et\'e $X$ sur $k$, on a
$$
\overline{\psi}_n\left(\cupu{i\in I}A_i\right)
=
\sum_{\vide\,\neq \,J \,\subset \,I}
\overline{\psi}_n
\left(
\capu{i\in J}A_i\,
\setminus 
\cupu{i\notin J}A_i
\right)\\
=
\sum_{\vide\,\neq \,J \,\subset \,I}
(-1)^{\,1+\card{J}}
\,
\overline{\psi}_n
\left(
\capu{i\in J}A_i 
\right).
$$
\end{cor}
\begin{nota}
Soit $k$ un corps de caractéristique zéro.
Pour toute $k$-vari\'et\'e $X$ 
et tout $n\geq 1$,
on pose
$$
\Psi^{\chi}_n(X)=\chi(\overline{\psi}_n(X))
$$
et
$$
\Phi^{\chi}_n(X)
=
\sum_{d|n}d\,\,\Psi^{\chi}_d(X).
$$
\end{nota}
On peut donc voir $\Phi_n^{\chi}(X)$ comme l'image
par $\chi$ de la classe d'une hypoth\'etique formule 
dont l'interpr\'etation dans 
tout corps pseudo-fini $K$ contenant $k$ d\'efinirait l'ensemble des points 
de $X$ \`a valeur dans l'unique extension de degr\'e $n$ de $K$. 
De la d\'efinition de $\Phi_n^{\chi}(X)$ et de \cite[proposition 3.6.1 et \S 3.3] {DeLo:def_sets_motives} 
on  d\'eduit d'ailleurs ais\'ement la proposition suivante.
\begin{prop}
Soir $k$ un corps de type fini sur $\Q$.
Soit $R$ un anneau int\`egre et normal de type fini sur $\Z$, de corps des fractions $k$.
Si $x$ est un point ferm\'e de $\Spec(R)$, on note $\F_x$ le corps r\'esiduel en $x$ et 
$\Frob_x$ le frobenius en $x$. 
Pour toute $k$-variété $X$, il existe un \'el\'ement non nul $f$ de $R$ tel que pour tout points ferm\'e $x$ de $\Spec(R_f)$
on ait
$$
\Tr\Frob_x(\Phi_n^{\chi}(X))=\card{X_{\F_x}(\F_{x,n})}.
$$
\end{prop}

\subsubsection{Propri\'et\'es}
\begin{prop}\label{prop:phiuf}
Soit $k$ un corps de caractéristique zéro.
Soit $X$ une $k$-vari\'et\'e et $U$ un ouvert de $X$.
On a pour tout $n\geq 1$
$$
\Phi_n^{\chi}(X)=\Phi_n^{\chi}(U)+\Phi_n^{\chi}(X\setminus U).
$$
Si $A$ est un sous-ensemble constructible de $X$, $\Phi_n^{\chi}(A)$
est bien d\'efini et si $(A_i)_{i\in I}$ est une famille finie de sous-ensembles
constructibles de $X$ on a une formule similaire \`a celle du corollaire \ref{cor:psin_constr}.
\end{prop}
\begin{proof}[Démonstration]
Ceci découle de la définition de $\Phi_n^{\chi}$ et des
propriétés analogues de $\Psi_n^{\chi}$.
\end{proof}
\begin{prop}\label{prop:phixy}
Soit $k$ un corps de caractéristique zéro.
Soit $n\geq 1$ un entier.
Soit $X$ et $Y$ des vari\'et\'es sur $k$.
On a
$$
\Phi^{\chi}_n(X\times Y)=\Phi^{\chi}_n(X)\,\Phi^{\chi}_n(Y).
$$
Supposons en outre qu'il existe  un morphisme $X\to Y$ 
qui est une fibration localement Zariski triviale de fibre $Z$.
On a a alors
$$
\Phi^{\chi}_n(X)=\Phi^{\chi}_n(Y)\,\Phi^{\chi}_n(Z).
$$
L'application qui \`a $X$ associe $\Phi^{\chi}_n(X)$ 
s'\'etend en un morphisme d'anneaux
$$
\Phi^{\chi}_n(X)\,:\,\mk\longto \mkc.
$$
\end{prop}
Cette proposition d\'ecoule de la d\'efinition de $\Phi_n^{\chi}$
et du lemme \ref{lm:rel} ci-dessous.
\begin{rem}
Evgeny Gorsky m'a signalé que la proposition \ref{prop:phixy} découlait
du fait (démontré par \cite{Hei:func_eq}) que la structure de $\lambda$-anneau définie sur
l'anneau de Grothendieck des motifs par la fonction zêta de Hasse-Weil
motivique était spéciale. La démonstration proposée ici est de nature
plus arithmétique.
\end{rem}
\begin{lemme}\label{lm:rel}
Soit $k$ un corps de caractéristique zéro.
Soit $X$ et $Y$ des vari\'et\'es sur $k$.
Pour tout $n\geq 1$, 
on a
$$
\Psi^{\chi}_n(X\times Y)
=
\sum_{
\substack{
d|n,\,e|n \\
~\\
d\vee e=n
}
}
\frac{d\,e}{n}\,\Psi^{\chi}_d(X)\,\Psi^{\chi}_e(Y).
$$
\end{lemme}
\begin{proof}[D\'emonstration]
Gr\^ace au lemme \ref{lm:rel_psi_UF}, en prenant des recouvrements ouverts affines et en stratifiant
on peut supposer $X$ et $Y$ affines, normales et irr\'eductibles.

Dans toute la d\'emonstration, pour tout entier $n$, 
on identifie $\grsym_n$ au groupe des bijections de $\Z/n\,\Z$,
et on note $\sigma_n$ le $n$-cycle $i\mapsto i+1$.  

Soit $n\geq 1$.
Soit $d$ et $e$ des diviseurs de $n$ tels que $d \vee e=n$.
Nous utilisons les notations \ref{notas:x0n}.
Soit $Z_{d,e}\eqdef(X^d)_0\times (Y^e)_0$. 
L'action naturelles
du groupe $\grsym_d\times \grsym_e$ sur $Z_{d,e}$ 
induit un $\grsym_d\times \grsym_e$-rev\^etement 
étale
\begin{equation*}
Z_{d,e}\longto (X^{<d>})_0\times (Y^{<e>})_0.
\end{equation*}
Les injections diagonales $(X^d)_0 \to (X^d)^{\frac{n}{d}}$
et $(Y^e)_0 \to (Y^e)^{\frac{n}{e}}$ induisent un morphisme
\begin{equation*}
\pi_{d,e}\,:\,Z_{d,e} \longto (X\times Y)_0^{<n>}.
\end{equation*}

On note $C_{d,e}$
le sous-groupe cyclique d'ordre $n$ de $\grsym_d\times\grsym_e$ engendr\'e par $(\sigma_d,\sigma_e)$.
Les formules $\psi_d(X)\wedge\psi_e(Y)$ et $\varphi_{Z_{d,e},(X^{<d>})_0\times (Y^{<e>})_0,C_{d,e}}$
sont alors logiquement \'equivalentes.

Soit $G_{d,e}$ le sous-groupe maximal 
de $\grsym_d\times \grsym_e$ tel que le morphisme
$\pi_{d,e}$ se factorise \`a travers $Z_{d,e}/G_{d,e}$. 
On peut d\'ecrire $G_{d,e}$ de la mani\`ere
suivante.
On note $f_d$ 
(respectivement $f_e$) 
le morphisme naturel $\Z/n\Z\to \Z/d\Z$ 
(respectivement $\Z/n\Z \to \Z/e\Z$). 
Alors un \'el\'ement $(\sigma_1,\sigma_2)$ de $\grsym_d\times \grsym_e$ est dans 
$G_{d,e}$ si et seulement s'il existe un \'el\'ement $\sigma$ de $\grsym_n$ v\'erifiant
$$
\forall i\in \Z/n\Z,\quad \left\{
\begin{array}{rcl}
\sigma_1f_d(i)&=&f_d(\sigma(i))\\
\sigma_2f_e(i)&=&f_e(\sigma(i)).
\end{array}
\right.
$$
En d'autres termes, si on identifie $\Z/n\Z$ \`a un sous-ensemble de 
$\Z/d\Z\times \Z/e\Z$ via $f_d\times f_e$, $G_{d,e}$ est le sous-groupe
de  $\grsym_d\times \grsym_e$ constitu\'e des \'el\'ements qui stabilisent 
$\Z/n\Z$. 
En particulier $G_{d,e}$ s'identifie \`a un sous-groupe de $\grsym_n$,
not\'e $\grsym_n^{d,e}$.
Notons que 
$G_{d,e}$ contient $C_{d,e}$ (l'\'el\'ement de $\grsym_n^{d,e}$ correspondant \`a $(\sigma_d,\sigma_e)$
est $\sigma_d$).

Soit $K$ un corps pseudo-fini contenant $k$, 
$\overline{K}$ une cl\^oture alg\'ebrique
de $K$, 
$\tau$ un g\'en\'erateur topologique du groupe
de Galois absolu de $K$, et, pour tout entier $d$,  $K_d$ l'unique extension de degr\'e $d$
de $K$ dans $\overline{K}$, et $K'_d$ l'ensemble des g\'en\'erateurs de $K_d$, \ie  l'ensemble
des \'el\'ements de $\overline{K}$ dont l'orbite sous $\tau$ est de cardinal $d$.

Soit $z$ un \'el\'ement de $(X\times Y)_0^{<n>}(K)$ satisfaisant $\psi_n(X\times Y)$.
Ceci signifie que
$z$ s'identifie \`a un ensemble du type $\{(\tau^ix,\tau^iy)_{i\in \Z/n\Z}\}$,
o\`u  $(x,y)$ est un \'el\'ement de $(X\times Y)(\overline{K})$
tel que l'\'egalit\'e $(\tau^ix,\tau^iy)=(x,y)$ ait lieu si et seulement si $n$ divise $i$.
Soit $d$ et $e$ tels que $x\in X(K'_d)$ et $y\in Y(K'_e)$. 
On 
a alors n\'ecessairement $d\vee e=n$. 
Le couple $(d,e)$ est en fait l'unique couple
v\'erifiant $d\vee e=n$ et tel que $z$ se rel\`eve \`a un point g\'eom\'etrique
de $Z_{d,e}$.
Montrons que $z$ se rel\`eve en fait \`a un unique \'el\'ement 
de $(Z_{d,e}/G_{d,e})(K)$, et que cet \'el\'ement admet $C_{d,e}$ comme groupe de d\'ecomposition
dans le revêtement $Z_{d,e}\to Z_{d,e}/G_{d,e}$.
L'ensemble des points g\'eom\'etriques de $Z_{d,e}$ qui s'envoient sur $z$ 
est l'ensemble des éléments  
$
\left((x_j)_{j\in \Z/d\Z},(y_k)_{k\in \Z/e\Z}\right)\in (X^d\times Y^e) (\overline{K})
$
qui v\'erifient la propriété : il existe un \'el\'ement $\mu$ de $\grsym_n$ tel qu'on ait
$$
\forall i\in \Z/n\Z,\quad \left\{
\begin{array}{rcl}
x_{f_d(i)}&=&\tau^{\mu(i)}x\\
y_{f_e(i)}&=&\tau^{\mu(i)}y.
\end{array}
\right.
$$
Un tel \'el\'ement $\mu$ est alors n\'ecessairement dans $\grsym_n^{d,e}$.
On voit donc que  $\pi_{d,e}^{-1}(z)$ est une orbite sous $G_{d,e}$
et on v\'erifie par ailleurs facilement qu'elle est stable sous $\tau$. Ceci montre que $z$ 
se rel\`eve \`a un unique \'el\'ement de $(Z_{d,e}/G_{d,e})(K)$.

Montrons que cet \'el\'ement admet $C_{d,e}$ comme groupe de d\'ecomposition
dans le revêtement $Z_{d,e}\to Z_{d,e}/G_{d,e}$.
Notons
$$
(\xx,\yy)=\left((\tau^j x)_{j\in \Z/d\Z}\,,\,(\tau^k y)_{k\in \Z/e\Z}\right)
\in \pi_{d,e}^{-1}(z).
$$
Il suffit de montrer que la $C_{d,e}$-orbite de $(\xx,\yy)$ est $\tau$-stable 
et que pour tout sous-groupe strict $C'$ de $C_{d,e}$ la $C'$-orbite
de $(\xx,\yy)$
n'est pas $\tau$-stable. 
Ceci est imm\'ediat compte tenu du fait que $(\sigma_d,\sigma_e)(\xx,\yy)=\tau.(\xx,\yy)$
et que $(\sigma_d,\sigma_e)$ engendre $C_{d,e}$.

Montrons \`a pr\'esent qu'un \'el\'ement de $(Z_{d,e}/G_{d,e})(K)$ admettant $C_{d,e}$
pour groupe de d\'ecomposition dans le revêtement $Z_{d,e}\to Z_{d,e}/G_{d,e}$ 
s'envoie par $\pi_{d,e}$ sur un \'el\'ement
de $(X\times Y)_0^{<n>}(K)$ satisfaisant $\psi_n(X\times Y)$.
Un tel \'el\'ement se rel\`eve \`a un point g\'eom\'etrique $\left((x_j)\times(y_k)\right)$
de $Z_{d,e}$ v\'erifiant : il existe un entier $l$ premier \`a $n$ tel que 
pour tout $(k,j)$ on a $(x_{j+l},y_{k+l})=(\tau\,x_j,\tau \,y_l)$.
Pour un entier $m$ premier \`a $n$ convenable, cet \'el\'ement s'\'ecrit
donc $\left((\tau^{i\,m}\,x)(\tau^{k\,m}\,y)\right)$ avec 
$x\in X(K'_d)$ et $y\in Y(K'_e)$,
et son image dans $(X\times Y)_0^{<n>}$ est un point $K$-rationnel satisfaisant 
$\psi_n(X\times Y)$.

On note alors $\theta_{d,e}$ une formule 
sur $(X\times Y)_0^{<n>}$
telle que, 
pour corps pseudo-fini $K$ contenant $k$, 
$\theta_{d,e}(K)$
s'identifie \`a l'ensemble des \'el\'ements de $(X\times Y)_0^{<n>}(K)$ qui se rel\`event \`a un \'el\'ement
de $(Z_{d,e}/G_{d,e})(K)$ admettant $C_{d,e}$ comme groupe de d\'ecomposition. 
Ce qui pr\'ec\`ede montre que, d'une part, pour tout $(d,e)$ vérifiant $d\vee e=n$, 
les formules $\theta_{d,e}$ et $\varphi_{Z_{d,e},Z_{d,e}/G_{d,e},C_{d,e}}$ sont logiquement \'equivalentes
et d'autre part que les formules $(\theta_{d,e})_{d\vee e=n}$ forment une partition de $\psi_n(X\times Y)$.

D'apr\`es \eqref{eq:rel_form_gal} 
et ce qui pr\'ec\`ede, on a
\begin{align*}
\chif([\theta_{d,e}])
&=
\chif([\varphi_{Z_{d,e},Z_{d,e}/G_{d,e},C_{d,e}}])\\
&=
\frac{
\card{N_{\grsym_d\times\grsym_e}(C_{d,e})}
}
{
\card{N_{G_{d,e}}(C_{d,e})}
}
\chif([\varphi_{Z_{d,e},X^{<d>}_0\times Y^{<e>}_0,C_{d,e}}])\\
&=
\frac{
\card{N_{\grsym_d\times\grsym_e}(C_{d,e})}
}
{
\card{N_{G_{d,e}}(C_{d,e})}
}
\,
\Psi_d^{\chi}(X)\,\Psi_e^{\chi}(Y).
\end{align*}
Il suffit donc pour terminer la d\'emonstration de 
montrer la relation
$$
\frac{
\card{N_{\grsym_d\times\grsym_e}(C_{d,e})}
}
{
\card{N_{G_{d,e}}(C_{d,e})}
}
=
\frac{de}{n}.
$$
Un \'el\'ement $(\sigma_1,\sigma_2)$ 
de $\grsym_d\times\grsym_e$
est dans $N_{\grsym_d\times\grsym_e}(C_{d,e})$
si et seulement s'il existe un \'el\'ement $l$ de $(\Z/n\Z)^{\ast}$
tel que 
$$
(\sigma_1,\sigma_2)(\sigma_d,\sigma_e)=(\sigma_d,\sigma_e)^l(\sigma_1,\sigma_2)
$$
\ie  si et seulement si on a
$$
\left\{
\begin{array}{rcl}
\forall j\in \Z/d\Z,\quad \sigma_1(j+1)&=&\sigma_1(j)+l\\
\forall k\in \Z/e\Z,\quad \sigma_2(k+1)&=&\sigma_2(k)+l.
\end{array}
\right.
$$
Un tel \'el\'ement est ainsi enti\`erement d\'etermin\'e par la donn\'ee de $l\in(\Z/n\Z)^{\ast}$
et du couple $(\sigma_1(0),\sigma_2(0))$. Il sera dans $G_{d,e}$ si et seulement
si on a en outre $(\sigma_1(0),\sigma_2(0))\in \Z/n\Z$. 
Ainsi on a 
$$
\card{N_{\grsym_d\times\grsym_e}(C_{d,e})}
=
d\,e\,\card{(\Z/n\Z)^{\ast}}
$$
et
$$
\card{N_{G_{d,e}}(C_{d,e})}
=
n\,\card{(\Z/n\Z)^{\ast}}.
$$
On en déduit le r\'esultat annonc\'e.
\end{proof}
\subsection{Décomposition de la fonction zêta de Hasse-Weil motivique
en produit eulerien motivique}
Soit $k$ un corps de caractéristique zéro et $X$ une $k$-variété quasi-projective.
\`A partir de la d\'efinition de $\Phi_n^{\chi}$, un calcul standard
montre la relation
$$
\exp
\left(
\sum_{n\geq 1}\frac{\Phi^{\chi}_n(X)}{n}\,T^n
\right)
=
\exp\left(
\sum_{n\geq 1}\Psi^{\chi}_n(X)\,\log\left(\frac{1}{1-T^n}\right)
\right).
$$
Nous allons montrer que cette derni\`ere expression 
est \'egale \`a $\zetaxmc(T)$.

\begin{notas}\label{notas:binom}
Soit $A$ une $\Q$-algèbre.
Pour tout \'el\'ement $x$ de $A$ et tout $n\geq 1$, 
on pose
$$
\binom{x}{n}
=
\frac{
\underset{i=0}{\overset{n-1}\prod}
(x-i)
}
{n!}
.
$$
Soit $r\geq 1$ et 
$\ff=(f_1,\dots,f_r)\in (\N_{>0})^r$
vérifiant $f_1\leq \dots \leq f_r$.
On d\'efinit une partition 
$\{1,\dots,r\}=\coprodu{\gamma\in \Gammaf} I_{\gamma}$
par la condition 
$$
\forall i,j\in \{1,\dots,r\}, 
\quad
f_i=f_j\longeq \exists \gamma,\quad i,j\in I_{\gamma}.
$$
Pour $\gamma\in \Gammaf$, 
on pose $f_{\gamma}=f_i$ o\`u $i$ est un \'el\'ement de $\gamma$,
et $n_{\gamma}=\card{I_{\gamma}}$.

Enfin si $(x_n)$ une suite d'éléments de $A$ on pose 
$$
\left(x_{\ff}\right)
=
\prod_{\gamma\in \Gammaf} 
\binom{x_{f_{\gamma}}}
{n_{\gamma}}.
$$
\end{notas}
Le lemme suivant découle d'un calcul élémentaire.
\begin{lemme}\label{lm:form_expl_prod_eul}
Soit $A$ une $\Q$-algèbre, $E$ un ensemble fini non vide et  $P=1+\sum_{\nn\in \N^E\setminus \{0\}}a_{\nn}\,\bT^{\,\nn}$
un \'el\'ement de $A[[(T_e)_{e\in E}]]$.
On a alors pour toute suite $(x_n)$ d'élements de $A$ la relation
\begin{multline}
\exp\left(\sum_{n\geq 1}\,a_n\,\log(P(T_e^n)_{e\in E})\right)
\\
=1+\sum_{\mm\in \N^E\setminus \{0\}}
\left(
\sum_{r\geq 1}
\quad
\sum_{
\substack{
\ff\in (\N_{>0})^r
\\~\\
f_1\leq \dots \leq f_r
}
}
(x_{\ff})                                                                                          
\quad
\sum_{
\substack{
(\nn_1,\dots,\nn_r)\in (\N^E\setminus \{0\})^r
\\~\\
\sum \nn_i\,f_i=\mm
}
}
\quad
\prod_{i=1}^r a_{\nn_i}
\right)
\,
\bT^{\,\mm}.
\end{multline}
\end{lemme}

\begin{lemme}\label{lm:sansnom}
Soit $k$ un corps. Soit $n\geq 1$ et $\varphi$ une formule  \`a coefficients dans $k$
en les variables libres $(x_1,\dots,x_n)$.
On pose $\psi_1=\varphi$. 
Pour $m\geq 2$, soit $\psi_m$ la formule d'anneau
en les variables libres 
$(x_{i,j})_{\substack{i=1,\dots,n\\j=1,\dots,m}}$ donn\'ee par 
$$
\left(\bigwedge_{j=1}^m \varphi(x_{1,j},\dots,x_{n_j})
\right)
\bigwedge
\left(
\bigwedge_{
\substack{
j,k\in \{1,\dots,m\}
\\j\neq k
} 
}
(x_{1,j},\dots,x_{n_j})
\neq
(x_{1,k},\dots,x_{n_k})
\right).
$$
On a alors pour tout $m\geq 1$ la relation
$$
\left[
\psi_m
\right]
=
\prod_{j=0}^{m-1}
\left(\symb{\varphi}-j\right).
$$
\end{lemme}
\begin{proof}[D\'emonstration]
Soit $m\geq 2$. 
Les formules en les $n\,m$ variables libres $(\xx_1,\dots,\xx_m)$
$$
\psi_{m-1}(\xx_1,\dots,\xx_{m-1})\wedge \varphi(\xx_m)
$$
et
$$
\psi_m
\vee 
\bigvee_{j=1}^{m-1}
\left(
\psi_{m-1}(\xx_1,\dots,\xx_{m-1})
\wedge 
(\xx_m=\xx_j)
\right)
$$
sont logiquement \'equivalentes.
Pour $j=1,\dots m-1$ la formule 
$$
\psi_{m-1}(\xx_1,\dots,\xx_{m-1})
\wedge 
(\xx_m=\xx_j)
$$
est logiquement \'equivalente \`a $\psi_{m-1}$.
On a alors
$$
\left[
\psi_m
\right]
+
(m-1)
\left[
\psi_{m-1}
\right]
=
\left[
\psi_{m-1}
\right]
\left[
\varphi
\right]
$$
soit
$$
\left[
\psi_m
\right]
=
\left[
\psi_{m-1}
\right]
\left(
\left[
\varphi
\right]
-m+1)
\right)
$$
d'o\`u le r\'esultat en raisonnant par r\'ecurrence sur $m$.
\end{proof}
\begin{prop}\label{prop:prod_eul_mot}
Soit $k$ un corps de caractéristique zéro et $X$ une $k$-variété quasi-projective.
Dans l'anneau $\mkcq[[T]]$, on a l'\'egalit\'e 
$$
\zetaxmc(T)
=\exp
\left[
\sum_{n\geq 1}\Psi^{\chi}_n(X)\,\log\left(\frac{1}{1-T^n}\right)
\right].
$$ 
\end{prop}
\begin{proof}
Posons
$
\Pi_X(T)
=\exp
\Big[
\sum_{n\geq 1}\Psi^{\chi}_n(X)\,\log\left(\frac{1}{1-T^n}\right)
\Big].
$
Soit $F$ un ferm\'e de $X$ et $U=X\setminus F$.
D'apr\`es le lemme \ref{lm:rel_psi_UF}, 
on a
$
\Pi_X(T)=\Pi_U(T)\,\Pi_F(T).
$
Par ailleurs,
on a 
$
\zetaxmc(T)=\zetamc{U}(T)\,\zetamc{F}(T)
$ (\cf \cite[\S 13.3.1]{And:mot}).
Ainsi, en prenant des recouvrements ouverts affines et en stratifiant, 
on est ramen\'e \`a d\'emontrer la proposition dans le cas o\`u $X$ est affine, normale, irr\'eductible. 
Pour $r\geq 1$ et $\ff\in (\N_{>0})^r$,
on note 
$$
\cA_{\ff,m}
\eqdef
\left\{
(n_1,\dots,n_r)\in (\N_{>0})^r
,\quad
\sum_{i=1}^r n_i\,f_i=m
\right\}.
$$
D'apr\`es le lemme \ref{lm:form_expl_prod_eul},
il s'agit donc de
d\'emontrer pour tout $m\geq 1$ la relation 
\begin{equation}\label{eq:rel_xm}
\chi\left(\left[X^{<m>}\right]\right)
=
\sum_{r>0}
\quad
\sum_{
\substack{
\ff=(f_1,\dots,f_r)\in \N_{>0}^r
\\~\\
f_1\leq \dots \leq f_r
}
}
\left(\Psi^{\chi}_{\ff}(X)\right)
\,
\card{\cA_{\ff,m}}.
\end{equation}
Pour $n\geq 1$, on d\'esigne 
par $X^{(0)}_n$ l'ensemble des points ferm\'es
de $X$ de degr\'e $n$.
La formule \eqref{eq:rel_xm} est le pendant motivique de l'\'egalité
\begin{equation}\label{eq:rel_xm_kfini}
\card{X^{<m>}(k)}
=
\sum_{r>0}
\quad
\sum_{
\substack{
\ff=(f_1,\dots,f_r)\in \N_{>0}^r
\\~\\
f_1\leq \dots \leq f_r
}
}
\left(\card{X^{(0)}_{\ff}}\right)
\,
\card{\cA_{\ff,m}},
\end{equation}
qui est valable si $k$ est un corps fini.
Dans ce cadre, la relation \eqref{eq:rel_xm_kfini} est une cons\'equence
de la décomposition de la fonction zêta de Hasse-Weil en produit
eulérien, 
mais elle peut aussi se retrouver via un argument combinatoire direct.
La preuve de la relation  \eqref{eq:rel_xm} qui suit
est une adaptation  motivique d'un tel argument combinatoire.

Soit $m\geq 1$, $r\geq 1$ et
$\ff\in (\N_{>0})^r$ 
tel que  $f_1\leq \dots \leq f_r$.
On utilise les notations \ref{notas:binom}.
On a une action naturelle de 
$\grsym_{\Gammaf }
\eqdef
\underset{\gamma\in \Gammaf }{\prod}
\grsym_{h_{\gamma}}$
sur 
$\cA_{\ff,m}$, ainsi que 
sur
$\underset{i=1}{\overset{r}{\prod}} \left(X^{<f_i>}\right)_0$.
Soit $Z_{\ff}$ l'ouvert $\grsym_{\Gammaf }$-stable de 
$\underset{i=1}{\overset{r}{\prod}} \left(X^{<f_i>}\right)_0$
donn\'e par 
$$
\prod_{\gamma\in \Gammaf } \left(\prod_{i\in I_{\gamma}} \left(X^{<f_i>}\right)_0\right)_0.
$$ 
Soit $\varphi_{\ff}$ une formule sur $Z_{\ff}$ ayant la propriété
suivante : 
pour tout corps pseudo-fini $K$ contenant $k$,
$\varphi_{\ff}(K)$ s'identifie \`a l'ensemble des \'el\'ements 
$(y_1,\dots,y_r)$ de $Z_{\ff}(K)$
tels que, pour tout $i$, $y_i$ satisfait la formule $\psi_{f_i}(X)$. 
D'apr\`es le lemme \ref{lm:sansnom}, 
on a donc
\begin{equation}\label{eq:psifimi1}
\symb{\varphi_{\ff}}
=
\prod_{\gamma\in \Gammaf }
\,
\prod_{j=0}^{h_{\gamma}-1}
\,
\left(\symb{\psi_{f_{\gamma}}(X)}-j\right).
\end{equation}
Soit $\nn\in \cA_{\ff,m}$.
On note $\grsym_{\nn}$ le stabilisateur de $\nn$
sous l'action de $\grsym_{\Gammaf }$, et 
$$
\pi_{\ff,\nn}\,:\,
Z_{\ff}
\longto 
X^{<m>}
$$
le $k$-morphisme 
qui envoie le $r$-uplet de z\'ero-cycles $(C_1,\dots,C_r)$
sur $\sum_i n_i\,C_i$. Ce morphisme se factorise 
\`a travers $Z_{\ff}/\grsym_{\nn}$.
Soit  $\psi_{\ff,\nn}$ une formule sur $X^{<m>}$ 
telle que, pour tout corps pseudo-fini $K$ contenant $k$,
$\psi_{\ff,\nn}(K)$ est l'ensemble des \'el\'ements de 
$X^{<m>}(K)$ qui sont l'image par $\pi_{\ff,\nn}$ d'un \'el\'ement 
$(y_1,\dots,y_r)$ de $Z_{\ff}(K)$ satisfaisant $\varphi_{\ff}$.
Un \'el\'ement de $\psi_{\ff,\nn}(K)$
est donc un z\'ero-cycle $K$-rationnel s'\'ecrivant $\sum_{i=1}^r n_i\,P_i$ o\`u, pour tout $i$,
$P_i$ est un point fermé de degré $f_i$ de $X_K$,  et $P_i\neq P_j$ si $f_i=f_j$.
L'ensemble des pr\'eimages de cet \'el\'ement par $\pi_{\ff,\nn}$ forme
 une $\grsym_{\nn}$-orbite. Ainsi le morphisme $Z_{\ff}/\grsym_{\nn}\to X^{<m>}$ 
induit une bijection entre $\psi_{(\ff,\nn)}(K)$
et les \'el\'ements de $(Z_{\ff}/\grsym_{(\nn)})(K)$ qui se rel\`event \`a un \'el\'ement 
de $Z_{\ff}(K)$ satisfaisant $\varphi_{\ff}$.
Donc le morphisme $\pi_{\ff,\nn}$ fait de $\varphi_{\ff}$
un $\card{\grsym_{\nn}}$-revêtement de $\psi_{(\ff,\nn)}$.
D'après le théorème \ref{thm:denefloeser}, on a alors
\begin{equation*}
\chif\left(
\symb{\psi_{(\ff,\nn)}}
\right)
=
\frac{1}{\card{\grsym_{\nn}}}
\chif\left(
\left[\varphi_{(f_i)}\right]
\right).
\end{equation*}
En  notant $\ecA^0_{\ff,m}$ un syst\`eme de repr\'esentants
de $\ecA_{\ff,m}$ modulo l'action de $\grsym_{\Gammaf }$,
on en d\'eduit 
$$
\sum_{\nn\in \ecA^0_{\ff,m}}
\chif\left(
\symb{\psi_{\ff,\nn}}
\right)
=
\Big(
\sum_{\nn\in \ecA^0_{\ff,m}}
\frac{1}{\card{\grsym_{\nn}}}
\Big)
\,
\chif\left(
\symb{\varphi_{\ff}}
\right)
=
\frac{\card{\ecA_{\ff,m}}}
{\card {\grsym_{\Gammaf}}}
\chif\left(
\symb{\varphi_{\ff}}
\right).
$$
D'apr\`es \eqref{eq:psifimi1}, on a donc
\begin{equation}\label{eq:psifimi2}
\sum_{\bn\in \ecA^0_{\ff,m}}
\chif\left(
\left[\psi_{(\ff,\nn)}\right]
\right)
=
\left(\psi^{\chi}_{\ff}(X)\right)\,\card{\ecA_{\ff,m}}.
\end{equation}

L'interprétation de $\psi_{\ff,\nn}(K)$ en termes de zéro-cycles
utilisée ci-dessus montre par ailleurs que tout \'el\'ement de $X^{<m>}(K)$ satisfait $\psi_{\ff,\nn}$ pour
un unique $\ff$ et un $\nn\in \cA_{\ff,m}$ unique modulo l'action de
$\grsym_{\Gammaf}$.
Ainsi les formules
$$
\left(
\psi_{\ff,\nn}
\right)_{
\substack{
r>0, 
\\ \\
\ff\in \N_{>0}^r,
\\ \\
f_1\leq \dots \leq f_r,
\\ \\
\nn\in \ecA^0_{\ff,m}.
}
}
$$
forment une partition de $X^{<m>}$.
Ceci conclut la démonstration de la relation
\eqref{eq:rel_xm}.
\end{proof}
\begin{cor}\label{cor:psichi_chipsi}
Soit $k$ un corps de caractéristique zéro et $X$ une $k$-variété quasi-projective.
Pour tout $n\geq 1$ on a $\Phi_n^{\chi}(X)=\chi\left(\Phi_n(X)\right)$
et $\Psi_n^{\chi}(X)=\chi\left(\Psi_n(X)\right)$.
\end{cor}
\begin{proof}[D\'emonstration]
Ceci d\'ecoule de la proposition \ref{prop:prod_eul_mot} et des d\'efinitions
de $\Phi_n(X)$ et $\Psi_n(X)$.
\end{proof}
\begin{cor}\label{cor:psichi_in_fil}
Soit $k$ un corps de caractéristique zéro et $X$ une $k$-variété quasi-projective.
Pour tout $n\geq 1$, $\Psi^{\chi}_n(X)$ est un \'el\'ement
de $\fil^{-n\,\dim(X)}\mkcq$.
\end{cor}
\begin{proof}
Ceci d\'ecoule du corollaire \ref{cor:psichi_chipsi} 
et du point  \ref{lm:psi_in_fil} du lemme \ref{lm:prop:phin:psin}.
\end{proof}
\begin{cor}\label{cor:phichiad}
Pour tout $d\geq 1$ et tout $n\geq 1$, on a $\Phi_n^{\chi}(\A^{\,d})=\LL^{\,n\,d}$.
\end{cor}
\begin{proof}[D\'emonstration]
Ceci d\'ecoule du corollaire \ref{cor:psichi_chipsi} et du point
\ref{lm:phi(ad)}
du lemme \ref{lm:prop:phin:psin}.
\end{proof}
\begin{rem}
Le corollaire \ref{cor:psichi_in_fil} peut s'obtenir de mani\`ere plus
directe à partir de la définition de $\Psi_n^{\chi}(X)$ comme formule
galoisienne et de la relation \eqref{eq:rel_form_gal}.
\end{rem}
\begin{cor}\label{cor:convergence_chi}
Soit $k$ un corps de caractéristique zéro et $X$ une $k$-variété quasi-projective.
Soit $n_0\geq 1$ un entier et $P=1+\sum_{n\geq n_0}a_n\,T^n$ un
élément de $\mkcq[[T]]$.  
\'Ecrivons
$$
\exp\left[
\sum_{n\geq 1}
\Psi_n^{\chi}(X)\,\log(P(T^n))
\right]
=
\sum_{n\geq 0}
\alpha_n\,T^n.
$$
Soit
$
\left(V_{n}\right)_{n\geq 0}
$ 
une suite d'\'el\'ements de $\mhatcq$. 
On suppose qu'il existe une application
$\varphi\,:\,\N\to \Z$ 
v\'erifiant :
\begin{enumerate}
\item
$
\varphi(n)-\frac{n\,\dim(X)}{n_0}
\underset{n\to +\infty}{\longto} +\infty\quad;
$
\item
pour tout $n\in \N$, on a
$
V_{n}
\in 
\fil^{\,\varphi(n)}\,\mhatcq.
$ 
\end{enumerate}
Alors la s\'erie
$
\sumu{n\in \N}{\alpha}_n\,V_{n}
$
converge dans $\mhatcq$.
\end{cor}
\begin{proof}[D\'emonstration]
D'après le lemme \ref{lm:form_expl_prod_eul}, on a pour tout $n\geq 0$ la relation
$$
\alpha_n
=
\sum_{r>0}
\quad
\sum_{
\substack{
\ff=(f_1,\dots,f_r)\in \N_{>0}^r
\\~\\
f_1\leq \dots \leq f_r
}
}
\left(\Psi^{\chi}_{\ff}(X)\right)
\quad
\sum_{
\substack{
(n_1,\dots,n_r)\in \N_{>0}^r
\\~\\
\sum n_i\,f_i=n
}
}
\quad
\prod_{i=1}^r a_{n_i}.
$$
Les indices $\ff$ intervenant dans la somme v\'erifient tous
$
n_0\abs{\ff}\leq n
$.
D'après le corollaire 
\ref{cor:psichi_in_fil} 
on a $\left(\Psi^{\chi}_{\ff}(X)\right)\in \fil^{-\abs{\ff}\,\dim(X)}\mkcq$.
Ainsi,  
on a 
$
\alpha_n\in 
\fil^{\,-\frac{n\,\dim(X)}{n_0}}\,\mkcq.
$
On en d\'eduit le r\'esultat.
\end{proof}
\begin{rem}
On obtient un résultat analogue dans $\mkq$ en remplaçant 
$\exp\left[
\sum_{n\geq 1}
\Psi_n^{\chi}(X)\,\log(P(T^n))
\right]$
par
$\exp\left[
\sum_{n\geq 1}
\Psi_n(X)\,\log(P(T^n))
\right]$.
On utilise alors le point \ref{lm:psi_in_fil}  du lemme  \ref{lm:prop:phin:psin}.
\end{rem}
\begin{rem}
Soi $k$ un corps et $X$ une $k$-variété quasi-projective.
Soit $(A_n)_{n\geq 1}$ une famille de $k$-variétés quasi-projectives
et
$
P(T)=1+\sum_{n\geq 1}\symb{A_n}\,T^n$.
Les auteurs de \cite{GusLueMel:power} définissent alors la \og puissance $[X]$-\`eme de $P(T)$ \fg\ 
par la formule
\begin{equation}\label{eq:expr1-pipmot}
P(T)^{\symb{X}}=
1+
\sum_{k=1}^{\infty}
\Big\{
\sum_{\rho\geq 1}
\sum_{
\substack{
(k_1,\dots,k_\rho)\in \N^{\,\rho}
\\~\\
\sum j\,k_j=k
\\~\\
k_{\rho}\neq 0
}
}
\Big[
\Big(
\left(X^{\sum_j k_j}\right)_0 
\times 
\prod_{j=1}^\rho A_j^{k_j}
\Big)
/
\prod_{j=1}^{\rho} \grsym_{k_j}
\Big]
\Big\}
T^k
\end{equation}
où $\prod_j \grsym_{k_j}$ agit de mani\`ere diagonale sur chacun des facteurs
$\left(X^{\sum_j k_j}\right)_0$ et $\prod_{j=1}^\rho A_j^{k_j}$.

Par ailleurs, d'après le lemme \ref{lm:form_expl_prod_eul},
on a
\begin{equation}\label{eq:expr2-pipmot}
\Pi^{\text{mot}}_{X,P}(T)
=
1+
\sum_{k=1}^{\infty}
\left\{
\sum_{r\geq 1}
\quad
\sum_{
\substack{
\ff=(f_1,\dots,f_r)\in \N_{>0}^r
\\~\\
f_1\leq \dots \leq f_r
}
}
\left(\Psi_{\ff}(X)\right)
\,\,
\sum_{
\substack{
(n_1,\dots,n_r)\in \N_{>0}^r
\\~\\
\sum n_i\,f_i=k
}
}
\quad
\prod_{i=1}^r \symb{A_{n_i}}
\right\}
\,
T^{k}.
\end{equation}
Si $k$ est de caractéristique zéro,
la \og d\'ecomposition
arithm\'etique \fg\ 
(donn\'ee par la relation \eqref{eq:rel_form_gal}) 
des motifs des quotients de vari\'et\'es apparaissant dans l'expression 
\eqref{eq:expr1-pipmot} permet alors de montrer 
que $\chi(\Pi^{\text{mot}}_{X,P}(T))$ co\"\i ncide avec $\chi(P(T)^{\symb{X}})$.
En ce sens la notion de produit eulerien motivique correspond à une
décomposition arithmétique de la puissance formelle définie dans \cite{GusLueMel:power}.

Explicitons ce qui se passe au niveau de la partie de
la décomposition arithmétique correspondant \`a l'image des points
rationnels dans le quotient (dans ce qui suit tout se passe au niveau
de l'anneau $\mkc$, on omet d'écrire $\chi$ pour alléger l'écriture) :
il s'agit, pour $k\geq 1$ donn\'e, de comparer d'une part l'expression
$$
\sum_{\rho\geq 1}
\quad
\sum_{
\substack{
(k_1,\dots,k_\rho)\in \N^{\,\rho}
\\~\\
\sum j\,k_j=k
\\~\\
k_{\rho}\neq 0
}
}
\left[
\left(X^{\sum_j k_j}\right)_0 
\times 
\prod_{j=1}^\rho A_j^{k_j}
\right]
/
\prod_{j=1}^{\rho} k_j!
$$
qui s'\'ecrit encore (\cf  lemme \ref{lm:sansnom})
\begin{equation}\label{eq:expr1-pipmot-geom}
\sum_{\rho\geq 1}
\quad
\sum_{
\substack{
(k_1,\dots,k_\rho)\in \N^{\,\rho}
\\~\\
\sum j\,k_j=k
\\~\\
k_{\rho}\neq 0
}
}
\frac{
[X]\,([X]-1)\,\dots([X]-\sum_jk_j+1)
\, 
\prod_{j=1}^\rho \symb{A_j}^{k_j}
}
{\card{\prod_j \grsym_{k_j}}}
\end{equation}
et les termes du coefficient de $T^k$ dans \eqref{eq:expr2-pipmot} 
correspondant au cas o\`u
tous les $f_i$ sont \'egaux \`a $1$, soit
$$
\sum_{r\geq 1}
\quad
\left(\Psi_{\underbrace{(1,\dots,1)}_{r\text{ r\'ep\'etitions}}}(X)\right)
\,\,
\sum_{
\substack{
(n_1,\dots,n_r)\in \N_{>0}^r
\\~\\
\sum n_i=k
}
}
\quad
\prod_{i=1}^r \symb{A_{n_i}}
$$
ce qui s'\'ecrit aussi
\begin{equation}\label{eq:expr2-pipmot-geom}
\sum_{r\geq 1}
\quad
\frac
{[X]\,([X]-1)\,\dots([X]-r+1)}
{r!}
\,\,
\sum_{
\substack{
(n_1,\dots,n_r)\in \N_{>0}^r
\\~\\
\sum n_i=k
}
}
\quad
\prod_{i=1}^r \symb{A_{n_i}}
\end{equation}
L'\'egalit\'e de \eqref{eq:expr1-pipmot-geom} et \eqref{eq:expr2-pipmot-geom}
d\'ecoule alors d'un argument combinatoire \'el\'ementaire 
(si $(k_1,\dots,k_{\rho})\in \N^{\rho}$ avec $k_{\rho}\neq 0$, soit $r=\sum k_j$ ; il existe
$\frac{r!}{\prod k_j!}$ $r$-uplets distincts $(n_1,\dots,n_r)$ de $\N_{>0}^r$
v\'erifiant $k_j=\card{\{i,\,n_i=j\}}$).
\end{rem}

\section{Une formule d'inversion de M\"obius motivique}\label{sec:fonc_mob_mot}

Dans cette section, nous introduisons un analogue motivique de fonctions d'inversion de M\"obius
utilis\'ees pour traiter des probl\`emes de comptage de points de hauteur born\'ee sur les vari\'et\'es 
toriques 
(\cf  \cite{Pey:duke}, \cite{Sal:tammes}, \cite{dlB:compter_torique}, \cite{Bou:vtetor}).

Soit $E$ un ensemble fini non vide. 
On munit $\{0,1\}^E$ de l'ordre partiel usuel.
Soit  $B$ un sous-ensemble  de $\{0,1\}^E$ v\'erifiant la propriété
suivante : 
si $\nn\in B$ et $\nn'\geq \nn$ alors $\nn'\in B$. 
On note $B^{\min}$
l'ensemble des \'el\'ements minimaux de $B$, et
 $A=\{0,1\}^E\setminus B$. 
On d\'efinit une fonction $\mu^0_B\,:\,\{0,1\}^E\to \Z$ par la relation
\begin{equation}\label{eq:mu0b}
\forall \nn\in \{0,1\}^E,\quad 
\ind_A(\nn)
=
\sum_{0\leq \nn' \leq \nn}\mu^0_B(\nn').
\end{equation}
Pour $\nn\in \{0,1\}^E$
on pose
$$
\ell_B(\nn)
=
\card{\{(\nn')\in B^{\min},\quad \nn'\leq \nn\}}.
$$
On vérfie qu'on a alors
$$
\forall \nn\in \{0,1\}^E,
\quad
\mu^0_B(\nn)
=
\left\{
\begin{array}{lcl}
1 & \text{si} & \nn=0\\
0 & \text{si} & \nn\in A\setminus \{0\}\\
(-1)^{\,\ell_B(\nn)} 
& \text{si} & \nn\in B\setminus \{0\}.
\end{array}
\right.
$$
On d\'efinit un \'el\'ement $P_B$ de $\Z[T_e]_{e\in E}$ par 
$$
P_B(T_e)
=
\sum_{\nn\in \{0,1\}^E}\mu_B^0(\nn)\,\prod_{e\in E}\,T_e^{n_e}
$$
et un \'el\'ement $Q_B$ de $\Z[[T_e]]_{e\in E}$ par
$$
Q_B(T_e)
=
\frac{P_B(T_e)}{\produ{e\in E}(1-T_e)}
=
P_B(T_e)\,
\sum_{\dd\in \N^E}
\,\,
\prod_{e\in E}T_e^{d_e}.
$$
Pour $\dd=(d_e)\in \N^E$ on d\'efinit 
$\wt{\dd}=(\wt{d}_e)\in \{0,1\}^E$
par $\wt{d}_e=1$ si et seulement si $d_e\geq 1$.
\begin{lemme}\label{lm:rel_qb}
On a 
$$
Q_B(T_e)
=
\sum_{
\dd\in N^E
}
\ind_{A}\left(\wt{\dd}\right)
\,
\prod_{e\in E}
T_e^{d_e}.
$$
\end{lemme}
\begin{proof}[D\'emonstration]
Ceci d\'ecoule imm\'ediatement de la 
d\'efinition de $Q_B$ et de la relation \eqref{eq:mu0b}.
\end{proof}

Soit $k$ un corps et $X$ une $k$-vari\'et\'e quasi-projective. 
Pour toute extension $K$ de $k$,
on note $\xeff(K)$ le mono\"\i de des z\'ero-cycles effectifs $K$-rationnels 
sur $X$, et, pour $d\in \N$, $\xeff_d(K)$ le sous-ensemble de $\xeff(K)$
constitu\'e des \'el\'ements de degr\'e $d$. Ainsi $\xeff_d(K)$ s'identifie \`a $X^{<d>}(K)$.

Soit $F$ un sous-ensemble de $E$. 
Pour $\dd\in \N^E$, 
la vari\'et\'e $\produ{e\in E}X^{\,<d_e>}$ 
des $E$-uples de z\'ero-cycles effectifs 
de $X$ de degr\'e $\dd$ 
contient un ouvert non vide $X_{\dd,F}$ d\'efini par la condition 
$\underset{e\in F}{\cap}\Supp(C_e)=\vide$.
On note alors $X^B_{\dd}$ l'ouvert de 
$\produ{e\in E}X^{\,<d_e>}$ 
d\'efini par 
$$
X^B_{\dd}
=
\bigcap_{\nn\in B}
X_{\dd,\{\nn=1\}}
=
\bigcap_{\nn\in B^{\min}}
X_{\dd,\{\nn=1\}}
.
$$
Pour toute extension $K$ de $k$, on a donc  l'égalité
$$
X^B_{\dd}(K)
=
\left\{
(C_e)\in \prod_{e\in E}\xeff_{d_e}(K),
\quad
\forall \nn\in B,
\quad
\bigcap_{e\in E,\,\,n_e=1}\,\Supp\left(C_e\right)=\vide
\right\}.
$$

\subsection{Le cas d'un corps fini}

Soit $k$ un corps fini et $X$ une $k$-vari\'et\'e quasi-projective. 
Soit 
$$
\ecA_X^B
=
\left\{
(C_e)
\in 
\left(
\xeff(k)
\right)^E,
\quad
\forall n\in B,\quad
\bigcap_{e\in E,\,\,n_e=1}\,\Supp\left(C_e\right)=\vide
\right\}.
$$
Il existe alors une unique fonction 
$
\wt{\mu}^B_X\,:\,\left(\xeff(k)\right)^E\to \Z
$
v\'erifiant la condition\footnote{
On munit $\left(\xeff(k) \right)^E$ de l'ordre partiel usuel.
}
\begin{equation}\label{eq:def_muxb}
\forall (C_e)\in \left(\xeff(k)\right)^E,
\quad
\ind_{\ecA_X^B}(C_e)
=
\sum_{0\leq (C'_e)\leq (C_e)}\,
\wt{\mu}^B_X((C'_e)).
\end{equation}
\begin{prop}\label{prop:rel_muxb}
On a les d\'ecompositions en produit eul\'erien
\begin{equation}\label{eq:zwtmu}
\sum_{(C_e)\in \left(\xeff(k)\right)^E}
\wt{\mu}^B_X((C_e))
\,
\prod_{e\in E}\,T_e^{\deg(C_e)}
=
\prod_{x\in X^{(0)}}P_B\left((T_e^{\deg(x)})\right)
\end{equation}
et
\begin{equation}\label{eq:zxbt}
\sum_{(C_e)\in \left(\xeff(k)\right)^E}
\ind_{\ecA_X^B}(C_e)
\,
\prod_{e\in E}\,T_e^{\deg(C_e)}
=
\prod_{x\in X^{(0)}}Q_B\left((T_e^{\deg(x)})\right).
\end{equation}
\end{prop}
\begin{proof}[D\'emonstration]
D'après \eqref{eq:def_muxb}, le membre de gauche
de  \eqref{eq:zxbt} est égal au membre de gauche de \eqref{eq:zwtmu}
multiplié par $\prod_{e\in E}Z_X(T_e)$.
Ainsi l'une des deux relations \eqref{eq:zwtmu} ou \eqref{eq:zxbt}
entraîne aussitôt l'autre.
La démonstration de ces relations est classique, 
et se base sur le fait que
$
\wt{\mu}^B
$
est une fonction multiplicative.
La proposition 1 (p. 180) de \cite{Bou:vtetor} 
(elle-m\^eme inspir\'ee de la proposition analogue dans le cas des corps de nombres
que l'on trouve dans \cite{Sal:tammes} ou \cite{Pey:ecoledete})
traite le cas particulier 
d'une fonction d'inversion de M\"obius associ\'ee \`a un \'eventail 
(\cf  sous-section \ref{subsec:mob_ev}) et la preuve est la m\^eme dans le cas g\'en\'eral.
\end{proof}
D\'efinissons \`a pr\'esent une fonction 
$\mu^B_X\,:\,\N^E\to \Z$
en posant
$$
\mu^B_X((d_e))
=
\sum_{(C_e),\,\deg(C_e)=d_e}
\wt{\mu}_X^B(C_e).
$$
On a alors
\begin{equation}\label{eq:rel_mub_cl}
\forall \dd\in \N^E,\quad
\card{X^B_{(d_e)}(k)}
=
\sum_{0\leq \dd'\leq \dd}
\mu^B_X(\dd')\,\prod_{e\in E}\,\card{X^{<d_e-d'_e>}(k)}.
\end{equation}
Pour $n\geq 1$, on d\'esigne 
par $X^{(0)}_n$ l'ensemble des points ferm\'es
de $X$ de degr\'e $n$.
De la proposition \ref{prop:rel_muxb}, 
on d\'eduit les relations
\begin{equation}\label{eq:zeta_mub_cl}
\sum_{\dd\in \N^E}\,\mu_X^B(\dd)\prod_{e\in E}T_e^{d_e}
=
\prod_{n\geq 1}
P_{B}(T_e^n)^{\card{X^{(0)}_n}}
\end{equation}
et
\begin{equation}\label{eq:zeta_xb_cl}
\sum_{\dd\in \N^E}\,
\card{X^B_{(d_e)}(k)}
\,
\prod_{e\in E}T_e^{d_e}
=
\prod_{n\geq 1}
Q_{B}(T_e^n)^{\card{X^{(0)}_n}}.
\end{equation}
Ce sont ces relations dont on veut obtenir une version motivique.

\subsection{Un analogue motivique}\label{subsec:ad_mot}

On considère dans cette sous-section un corps $k$ de caract\'eristique z\'ero.
On va démontrer des relations  dans l'anneau de motifs virtuels  $\mkc$.
Le problème de la démonstraton de relations analogues dans l'anneau $\mk$
est discut\'e \`a la sous-section \ref{subsec:ad_mot_gen}.
Soit $X$ une $k$-variété quasi-projective.
On mime la relation \eqref{eq:rel_mub_cl}
et on d\'efinit une \termin{fonction de M\"obius motivique}
$
\mubcx\,:\,\N^{E}\rightarrow \mkc  
$
par la relation
\begin{equation}\label{eq:rel_mub_mot}
\forall\,\dd\in \N^{E},
\quad
\chi\left(
\left[X^B_{\dd}\right]
\right)
=
\sum_{0\leq \dd^\prime\leq \dd}
\mubcx(\dd^\prime)\,
\prod_{e\in E}
\chi\left(\left[X^{<d_e-d'_e>}\right]\right).
\end{equation}
Si on pose
$$
Z^B_{X}((T_e))
=
\sum_{\dd\in \N^E}\,
\chi\left(\left[X^B_{\dd}\right]\right)
\,
\prod_{e\in E}T_e^{d_e}
$$
et
$$
Z_{\mubcx}((T_e))
=
\sum_{\dd\in \N^E}
\,
\mubcx(\dd)
\,
\prod_{e\in E}T_e^{d_e},
$$
on a donc la relation
\begin{equation}\label{eq:rel_zb_zmu_mot}
Z^B_X((T_e))
=
Z_{\mubcx}((T_e))
\,
\prod_{e\in E}
\zetaxmc(T_e).
\end{equation}
\begin{thm}\label{thm:prod_eul_mot_zmubcx}
Soit $k$ un corps de caractéristique zéro et $X$ une $k$-variété quasi-projective.
Soit $B$ un sous-ensemble de $\{0,1\}^E$ v\'erifiant la propriété
suivante : si $\nn\in B$ et $\nn'\geq \nn$ alors $\nn'\in B$.
On a les d\'ecompositions en produit eul\'erien motivique
\begin{equation}\label{eq:zmubcx}
Z_{\mubcx}((T_e))
=
\exp
\left[
\sum_{n\geq 1}
\Psi^{\chi}_n(X)\,\log\left(P_{B}(T_e^n)\right)
\right]
\end{equation}
et
\begin{equation}\label{eq:zbx}
Z^B_X((T_e))
=
\exp
\left[
\sum_{n\geq 1}
\Psi^{\chi}_n(X)\,\log\left(Q_{B}(T_e^n)\right)
\right]
.
\end{equation}
\end{thm}
\begin{proof}[D\'emonstration] 
Compte tenu de \eqref{eq:rel_zb_zmu_mot} et de la 
proposition \ref{prop:prod_eul_mot}, 
l'une des deux relations \eqref{eq:zmubcx}
et \eqref{eq:zbx} entraîne aussitôt l'autre.

Montrons la  relation \eqref{eq:zbx}.
Soit
$
\Pi^B_X(\bT)
=
\exp
\Big[
\sum_{n\geq 1}
\Psi^{\chi}_n(X)\,\log\left(Q_{B}(T_e^n)\right)
\Big].
$
Soit $F$ un ferm\'e de $X$ et $U=X\setminus F$.
On a une stratification de $X^B_{\dd}$
en sous-vari\'et\'es localement ferm\'ees
\begin{align*}
X^B_{\dd}&
=\coprod_{0\leq \dd' \leq \dd}\,\left(\prod_{e\in E} U^{<d'_e>}\times F^{<d_e-d'e>}\right)
\bigcap X^B_{\dd}
=\coprod_{0\leq \dd' \leq \dd}U^B_{\dd'}\times F^B_{\dd'}.
\end{align*}
On en déduit la relation
$
Z^B_{X}(T)=Z^B_{U}(T)\,Z^B_{F}(T).
$
Par ailleurs le lemme \ref{lm:rel_psi_UF}
entra\^\i ne la relation
$
\Pi^B_{X}(T)=\Pi^B_{U}(T)\,\Pi^B_{F}(T).
$
Ainsi, en prenant des recouvrements ouverts affines et en stratifiant, 
on est ramen\'e \`a démontrer la  relation \eqref{eq:zbx}
dans le cas o\`u $X$ est affine, normale, irr\'eductible.
Pour $\dd\in \N^E$, $r\geq 1$ et $\ff\in \N_{>0}^r$,
on note
$$
\ecA^A_{\ff,\dd}
\eqdef
\left\{
(\nn_1,\dots,\nn_r)
\in 
(\N^E\setminus \{0\})^r,
\quad
\sum \nn_i\,f_i=\dd\,
\text{ et }
\,
\forall i=1\,\dots,r,\quad \wt{\nn_i}\in A
\right\}.
$$
D'apr\`es les lemmes 
\ref{lm:form_expl_prod_eul} 
et 
\ref{lm:rel_qb}, 
démontrer \eqref{eq:zbx}  
revient \`a établir
pour tout $\dd\in \N^E$
la relation
\begin{equation}\label{eq:rel_xbd}
\chi
\left(
\left[
X^B_{\dd}
\right]
\right)
=
\sum_{r\geq 1}
\quad
\sum_{
\substack{
\ff=(f_1,\dots,f_r)\in \N_{>0}^r
\\
f_1\leq \dots \leq f_r
}
}
\left(\Psi_{\ff}(X)\right)                                                                                          
\,
\card{\ecA^A_{\ff,\dd}}
\end{equation}
Fixons $\dd\neq 0$. Soit  $r\geq 1$ et $\ff\in (\N_{>0})^r$ 
tel que $f_1\leq \dots \leq f_r$. 
On utilise les notations \ref{notas:binom} et on reprend 
les notations et la d\'emarche de la preuve de la proposition
\ref{prop:prod_eul_mot}. 
Le groupe $\grsym_{\Gammaf}$ agit 
sur $(\N^E\setminus \{0\})^{r}$.
L'ensemble $\ecA^A_{\ff,\dd}$ est stable sous cette action.
Soit $(\nn_i)\in \ecA^A_{\ff,\dd}$.
On note $\grsym_{(\nn_i)}$ le stabilisateur de 
$(\nn_i)$ sous l'action de $\grsym_{\Gammaf}$,
et
$$
\pi_{\ff,(\nn_i)}\,:\,
Z_{\ff}
\longto 
\prod_{e\in E} X^{<d_e>}
$$
le $k$-morpisme qui envoie le $r$-uplet de z\'ero-cycles $(C_1,\dots,C_r)$
sur $\left(\sum_i n_{i,e}\,C_i\right)_{e\in E}$. 
Ce morphisme se factorise 
\`a travers $Z_{\ff}/\grsym_{(\nn_i)}$.

Soit $\psi_{\ff,(\nn_i)}$ une formule 
vérifiant la propriété suivante :
 pour tout corps pseudo-fini $K$ contenant $k$,
$\psi_{\ff,(\nn_i)}(K)$ est l'ensemble des \'el\'ements de 
$\produ{e\in E} X^{<d_e>}(K)$ qui sont l'image par $\pi_{\ff,(\nn_i)}$ d'un \'el\'ement 
$(y_1,\dots,y_r)$ de $Z_{\ff}(K)$ satisfaisant $\varphi_{\ff}$.
La condition $\wt{\nn_i}\in A$ dans la d\'efinition de $\ecA^A_{\ff,\dd}$ 
entra\^\i ne que de tels \'el\'ements sont en 
particulier dans $X^B_{\dd}(K)$.

Un \'el\'ement de $\produ{e\in E}X^{<d_e>}(K)$ satisfaisant $\psi_{\ff,(\nn_i)}$
est donc un $E$-uplet de z\'ero-cycles $K$-rationnels s'\'ecrivant $(\sum_{i=1}^r n_{i,e}\,C_i)$ o\`u, pour tout $i$,
$C_i$ est un point fermé de $X_K$ de degré $f_i$  et $C_i\neq C_j$ si $f_i=f_j$. 
Les pr\'eimages d'un tel \'el\'ement dans $Z_{\ff}$ forment une orbite sous $\grsym_{(\nn_i)}$.
Ainsi le morphisme $Z_{\ff}/\grsym_{(\nn_i)}\to \produ{e\in E} X^{<d_e>}$ induit
une bijection entre  $\psi_{\ff,(\nn_i)}(K)$
et les \'el\'ements de $(Z_{(\ff}/\grsym_{(\nn_i)})(K)$ qui se rel\`event \`a un \'el\'ement 
de $Z_{\ff}(K)$ satisfaisant $\varphi_{\ff}$.
Donc le morphisme $\pi_{\ff,(\nn_i)}$
 fait de $\varphi_{\ff}$
un $\card{\grsym_{(\nn_i)}}$-revêtement de $\psi_{\ff,(\nn_i)}$.
D'après le théorème \ref{thm:denefloeser}, on a alors
\begin{equation*}
\chif\left(
\symb{\psi_{\ff,(\nn_i)}}
\right)
=
\frac{1}{\card{\grsym_{(\nn_i)}}}
\chif\left(
\symb{\varphi_{\ff}}
\right).
\end{equation*}
L'interprétation de $\psi_{\ff,(\nn)_i}(K)$ en termes de zéro-cycles
utilisée ci-dessus montre par ailleurs que tout \'el\'ement de $X^B_{\dd}(K)$ satisfait $\psi_{\ff,(\nn_i)}$ pour
un unique $\ff$ et un $(\nn_i)\in \cA_{\ff,\dd}$ unique modulo l'action de
$\grsym_{\Gammaf}$.
On peut alors conclure comme dans la preuve de la proposition \ref{prop:prod_eul_mot}.
\end{proof}
\begin{cor}\label{cor:muxb}
Soit $k$ un corps de caractéristique zéro et $X$ une $k$-variété quasi-projective.
Soit $B$ un sous-ensemble de $\{0,1\}^E$ v\'erifiant la propriété
suivante : si $\nn\in B$ et $\nn'\geq \nn$ alors $\nn'\in B$.
Soit $\nu_B$ la valuation de $P_B-1$
et 
$
\left(V_{\dd}\right)_{\dd\in \N^E}
$ 
une famille d'\'el\'ements de $\mhatcq$.
On suppose qu'il existe une application $\varphi\,:\,\N\to \Z$ 
v\'erifiant :
\begin{enumerate}
\item
$
\varphi(n)-\frac{n\,\dim(X)}{\nu_B}
\underset{n\to +\infty}{\longto} +\infty\quad ;
$
\item
pour tout $\dd\in \N^E$, 
$
V_{\dd}
\in 
\fil^{\,\varphi(\sum d_e)}\,\mhatcq.
$
\end{enumerate}
Alors la s\'erie
$
\sumu{\dd\in \N^{E}}\mubcx(\dd)\,V_{\dd}
$
converge dans $\mhatcq$.
\end{cor}
\begin{proof}
Ceci se déduit du théorème \ref{thm:prod_eul_mot_zmubcx} grâce au
lemme
\ref{lm:form_expl_prod_eul}
et une adaptation aisée de la preuve du corollaire \ref{cor:convergence_chi}.
\end{proof}
\subsection{Questions dans l'anneau de Grothendieck des variétés}\label{subsec:ad_mot_gen}
Dans cette sous-section, $k$ est un corps quelconque. 
On reprend les notations de l'introduction de la  section \ref{sec:fonc_mob_mot}.
Soit $X$ une $k$-variété quasi-projective.
On d\'efinit une fonction
$
\mubmx\,:\,\N^{E}\rightarrow \mk
$
par la relation
\begin{equation}\label{eq:rel_mub_mot_2}
\forall\,\dd\in \N^{E},\quad
\symb{X^B_{\dd}}
=
\sum_{0\leq \dd^{\,\prime}\leq \dd}
\mubmx(\dd^{\,\prime})\,
\prod_{e\in E}
\left[X^{<d_e-d'_e>}\right].
\end{equation}
Si $k$ est de caract\'eristique z\'ero, on a donc pour tout $\dd$
la relation $\chi(\mubmx(\dd))=\mubcx(\dd)$.
Si $k$ est fini on a pour tout $\dd$ la relation $\#_k\,\mubmx(\dd)=\mu^B_X(\dd)$.
Au vu du th\'eor\`eme \ref{thm:prod_eul_mot_zmubcx},
on peut se poser la question suivante.
\begin{question}\label{qu:prod_eul_mot_zmubmx}
La relation
$$
\sum_{\dd\in \N^E}
\,
\mubmx(\dd)
\,
\prod_{e\in E}T_e^{d_e}
=
\exp
\left[
\sum_{n\geq 1}
\Psi_n(X)\,\log\left(P_{B}(T_e^n)\right)
\right]
$$
est-elle vérifiée ?
\end{question}
Une r\'eponse positive \`a la question \ref{qu:prod_eul_mot_zmubmx}
fournirait une autre d\'emonstration du th\'eor\`eme \ref{thm:prod_eul_mot_zmubcx}.
Elle entra\^\i nerait \'egalement une r\'eponse positive \`a la question suivante :
\begin{question}\label{ques:muxb}
Soit $\nu_B$ la valuation de $P_B$.
Soit
$
\left(V_{\dd}\right)_{\dd\in \N^E}
$ 
une famille d'\'el\'ements de $\mhatq$. 
On suppose qu'il existe une application
$\varphi\,:\,\N\to \Z$
v\'erifiant :
\begin{enumerate}
\item
$
\varphi(n)-\frac{n\,\dim(X)}{\nu_B}
\underset{n\to +\infty}{\longto} +\infty\quad;
$
\item
pour tout $(d_e)\in \N^{E}$,
$
V_{\dd}\in 
\fil^{\,\,\varphi(\sum d_e)}\,\mhatq.
$
\end{enumerate}
Est-il vrai que la s\'erie
$$
\sum_{(\dd)\in \N^{E}}\mubmx(\dd)\,V_{\dd}
$$
converge dans $\mhatq$ ?
\end{question}
On montre ci-dessous que les r\'eponses aux questions
\ref{qu:prod_eul_mot_zmubmx} et \ref{ques:muxb} sont positives pour un cas particulier d'ensemble $B$.
Dans le cadre de l'application à l'étude des fonctions zêta des
hauteurs motivique,
ceci permet  d'obtenir des résultats dans l'anneau $\mk$ dans le cas des espaces projectifs et des surfaces de Hirzebruch.
Pour traiter le cas d'une vari\'et\'e torique g\'en\'erale, faute de pouvoir
montrer que la r\'eponse est positive, nous devrons nous contenter de
résultats dans l'anneau $\mkc$ 
pour pouvoir utiliser le th\'eor\`eme \ref{thm:prod_eul_mot_zmubcx} 
et le corollaire \ref{cor:muxb}.

\subsection{Calcul de $\mubmx$ dans un cas particulier}
Soit $k$ un corps et $X$ une $k$-variété quasi-projective.
On d\'efinit une fonction $\muxm\,:\,\N\to \mk$ par la relation
$$
\left(\sum_{d\geq 0}\muxm(d)T^{d}\right)\zetaxm(T)=1.
$$
On a donc pour tout $d\geq 0$ la relation
\begin{equation}\label{eq:inv_mu_x}
\sum_{0\leq r \leq d}\muxm(r)\left[X^{<d-r>}\right]=0.
\end{equation}
Soit $E$ un ensemble fini non vide.
On d\'efinit une fonction $\mu_X^{E}\,:\,\N^{E}\to \mk$ en posant 
\begin{equation}\label{eq:muxi}
\mu_X^{E}(\dd)
=
\left\{\begin{array}{l}\muxm(d)\quad\text{si}\quad{d_e=d}\text{ pour tout }e\in E,
\\
0\quad\text{sinon.}\end{array}\right.
\end{equation}
Une d\'efinition \'equivalente de $\mu_X^{E}$ est d'imposer la relation 
$$
\left(
\sum_{\dd\in \N^E}\mu_X^{E}(\dd)\,\bT^{\,\dd}
\right)\,
\zetaxm\left(\prod_{e\in E}T_e\right)
=1.
$$
\begin{prop}\label{prop:cas_part_B}
Soit $B$ le sous-ensemble de $\{0,1\}^E$ réduit à l'élément constant
égal à $1$.
Alors $\mubmx=\mu^E_{X}$.
\end{prop}
\begin{proof}[D\'emonstration]
Soit 
$
W_{\dd,\delta}
$
la sous-vari\'et\'e 
de $\produ{e\in E}X^{<d_e>}$ des $E$-uples 
de z\'ero-cycles effectifs dont l'intersection est de degr\'e $\delta$.
En particulier $W_{\dd\, ,\,0}=X_{(\dd,E}$. 
On a un isomorphisme 
$$
W_{\dd,\delta}
\isom X^{\,<\delta>}\times X_{(d_e-\delta),E}
$$
et une fibration en sous-vari\'et\'es localement ferm\'ees
$$
\prod_{e\in E}X^{<d_e>}
=
\coprod_{0\leq d \leq \Min(d_e) }W_{(d_e),\delta}.
$$
On a donc
$$
\prod_{e\in E}
\left[X^{\,<d_e>}\right]
=
\sum_{\delta=1}^{\Min(d_e)} 
\left[X^{<\delta>}\right]\,\left[X_{(d_e-\delta),E}\right].
$$
Pour chaque $E$-uple $\dd^{\,\prime}$ 
vérifiant $0\leq \dd^{\,\prime}\leq \dd$.
écrivons la relation ci-dessus et multiplions la  par $\mu^E_{X}\left(\dd^{\,\prime}\right)$.
En sommant toutes les relations obtenues, 
et compte tenu de  \eqref{eq:inv_mu_x},
on obtient la relation
$$
\left[X_{\dd,E}\right]
=
\sum_{0\leq \dd^{\,'}\leq \dd} 
\mu^E_{X}\left(\dd^{\,\prime}\right)
\,
\prod_{e\in E}
\left[X^{<d_e-d'_e>}\right].
$$
Or, l'hypothèse sur $B$ entraîne   pour tout $\dd$
l'égalité 
$
X^B_{\dd}
=
X_{\dd,E}$. On en déduit le résultat.
\end{proof}
\begin{cor}\label{cor:cas_part_B}
Soit $B$ un sous-ensemble de $\{0,1\}^E$ 
tel que si $\nn\in B$ et $\nn'\geq \nn$ alors $\nn'\in B$.
On suppose en outre que $B$ v\'erifie l'hypoth\`ese suivante : 
il existe une partition 
$
E=E_{\beta}\,\coprod\,\sqcupu{\gamma\in \Gamma}E_{\gamma}
$
de $E$ telle qu'on ait
$$
B^{\min}
=
\left\{
\nn\in \N^E,
\quad
\exists \gamma\in \Gamma,
\quad
\left(n_e=1 \eq e\in E_{\gamma}\right)
\right\}.
$$
Alors on a 
$$
\forall \dd\in \N^E,
\quad
\mubmx(\dd)
=
\prod_{e\in E_{\beta}}
[X^{<d_e>}]
\prod_{\gamma\in \Gamma}\,
\mu^{E_{\gamma}}_X
((d_e)_{e\in E_{\gamma}}).
$$
\end{cor}

\subsection{Fonction de M\"obius motivique associ\'ee \`a un \'eventail}\label{subsec:mob_ev}

Soit $N$ un $\Z$-module libre de rang fini dont on notera $r$ le rang.
On rappelle bri\`evement la notion d'\'eventail de $N$.
Une partie $\sigma$ de $N\otimes \R$ est un \termin{c\^one poly\'edral rationnel} 
de  $N\otimes \R$
si elle s'écrit
$
\sigma=\sum_{i\in I}\R_{\geq 0}\,m_i
$
o\`u $I$ est un ensemble fini et les $(m_i)$ sont dans $N$. 
Un c\^one polyédral rationnel $\sigma$ est dit \termin{strictement convexe} si $\sigma\cap -\sigma=\{0\}$.
Un \termin{\'eventail} de $N$ est un ensemble fini $\Sigma$ 
de c\^ones polyédraux rationnels strictement convexes de $N\otimes \R$, 
v\'erifiant les conditions suivantes : 
\begin{itemize}
\item
toute face d'un c\^one de $\Sigma$ est un c\^one de $\Sigma$, 
\item
l'intersection de deux c\^ones de $\Sigma$ est une face de chacun des deux c\^ones. 
\end{itemize}
Un \'eventail $\Sigma$ est dit \termin{r\'egulier} si tout c\^one de $\Sigma$ 
est engendr\'e par une partie d'une $\Z$-base de $N$, 
et \termin{complet} si les c\^ones de $\Sigma$ recouvrent $N\otimes \R$.

Soit $\Sigma$ un \'eventail non r\'eduit \`a $\{0\}$. 
On note $\Sigma(1)$ l'ensemble des rayons de $\Sigma$, \ie 
l'ensemble des c\^ones de dimension 1 de $\Sigma$.
Pour $\sigma\in \Sigma$, on note $\sigma(1)$ l'ensemble des
\'el\'ements de $\Sigma(1)$ qui sont des faces de $\sigma$. 
Pour $\alpha\in\Sigma(1)$ on abr\'egera les notations $\alpha\in \sigma(1)$
et $\alpha\notin \sigma(1)$ en 
$\alpha\in \sigma$ et $\alpha\notin \sigma$.

Soit 
$\bs$ le sous-ensemble de $\{0,1\}^{\Sigma(1)}$ d\'efini par
$$
\bs
=
\left\{
\nn\in \{0,1\}^{\Sigma(1)},\quad
\forall \sigma\in \Sigma, \quad
\exists \,\alpha\notin\sigma,
\quad
\na=1
\right\}.
$$
Il est clair que si $\nn\in \bs$ et $\nn'\geq \nn$ alors $\nn'\in \bs$.
\begin{lemme}\label{val:pbs}
La valuation de $P_{\bs}-1$ est sup\'erieure ou \'egale \`a 2.
\end{lemme}
\begin{proof}[D\'emonstration]
Soit $\alpha_0\in \Sigma(1)$ et $(\na)\in \{0,1\}^{\Sigma(1)}$
v\'erifiant $n_{\alpha_0}=1$ et $n_{\alpha}=0$ pour $\alpha\neq\alpha_0$.
Il s'agit de montrer que $\na\notin \bs$. Mais ceci d\'ecoule aussit\^ot 
de la d\'efinition de $\bs$ et du fait que 
$\alpha$ est une face de $\alpha$.
\end{proof}
On d\'eduit du lemme \ref{val:pbs} et du corollaire \ref{cor:muxb} 
le critère de convergence suivant.
\begin{cor}\label{cor:musc}
Soit $k$ un corps de caractéristique zéro et $X$ une $k$-variété quasi-projective.
Soit $\Sigma$ un éventail.
Soit
$
\left(V_{\dd}\right)_{\dd\in \N^{\Sigma(1)}}
$ 
une famille d'\'el\'ements de $\mhatcq$ v\'erifiant
$$
\forall \,\,\dd\in \N^{\Sigma(1)},\quad V_{\dd}
\in 
\fil^{\,\,\dim(X)\,\abs{\dd}}\,\,\mhatcq.
$$
Alors la s\'erie
$
\sumu{\dd\in \N^{{\Sigma(1)}}}\mubsc(\dd)\,V_{\dd}
$
converge dans $\mhatcq$.
\end{cor}
\paragraph{Le cas des espaces projectifs}
Soit $n\geq 1$ et $\Sigma$ l'\'eventail de $\Z^n\otimes \R$
dont les rayons sont engendr\'es par les \'el\'ements de la base canonique
$(e_i)_{1\leq i \leq n}$ de $\Z^n$ et l'\'el\'ement $\sum_i e_i$. 
La vari\'et\'e torique associ\'ee est l'espace projectif de dimension $n$.
On a $\bs=\{(1,\dots,1)\}$.
Soit $k$ un corps et $X$ une $k$-variété quasi-projective. D'apr\`es la proposition
\ref{prop:cas_part_B} on a $\mubsm=\mu^{\Sigma(1)}_X$.
En particulier, les r\'eponses aux questions \ref{qu:prod_eul_mot_zmubmx}
et  \ref{ques:muxb} sont positives dans ce cas.

\paragraph{Le cas des surfaces de Hirzebruch}\label{par:musm_hirz}
Soit $m\geq 0$ un entier et $\Sigma$ l'\'eventail de $\Z^2\otimes \R$ 
dont les rayons sont engendr\'es par $\rho_1=(1,0)$, $\rho_2=(-1,m)$, $\rho_3=(0,1)$, $\rho_4=-\rho_3$. 
La vari\'et\'e torique associ\'ee est la $m$-\`eme surface de Hirzebruch. 
On a $\bs^{\min}=\{(1,0,1,0),(0,1,0,1)\}$.
Soit $k$ un corps et $X$ une $k$-variété quasi-projective.
D'apr\`es le corollaire \ref{cor:cas_part_B}, 
on a 
$$
\mubsm(d_1,d_2,d_3,d_4)
=
\mu^{\{1,3\}}_X(d_1,d_3)
\,
\mu^{\{2,4\}}_X(d_2,d_4).
$$
L\`a encore, les r\'eponses aux questions \ref{qu:prod_eul_mot_zmubmx} 
et  \ref{ques:muxb}
sont positives.

\section{Fonction z\^eta des hauteurs motivique}\label{sec:fonc_haut_mot} 

Soit $k$ un corps et $\cc$ une $k$-courbe projective, lisse
et g\'eom\'etriquement int\`egre. 
Soit $K$ le corps des fonctions de $\cc$. 
Consid\'erons une vari\'et\'e projective $V$ d\'efinie\footnote{
On pourrait en fait consid\'erer des vari\'et\'es d\'efinies sur $K$,
\ie  des familles non constantes, 
mais nous nous limiterons ici au cas
particulier o\`u le corps de d\'efinition est le corps des constantes.
}
sur $k$. 
Les \'el\'ements de $V(K)$ s'identifient donc aux $k$-morphismes $x\,:\,\cc\rightarrow V$. 
On fixe un fibr\'e en droites $\ecL$ sur $V$. 

\subsection{Le cas classique}

On suppose le corps  $k$ fini. La formule
$
h_{\ecL}(x)=\deg\left(x^{\ast}\ecL\right)
$
d\'efinit alors une hauteur d'Arakelov (logarithmique) sur $V(K)$, 
relative au faisceau $\ecL$. Si on suppose que la classe de $\ecL$ est \`a l'int\'erieur du c\^one effectif,
il existe 
un ouvert non vide $U_0$ de $V$ tel que pour tout ouvert $U$ de $U_0$ et tout entier $d\geq 0$, 
l'ensemble 
$
\left\{x\in U(K),\,h_{\ecL}(x)=d\right\}
$ 
est fini (\cf  \cite[Corollaire 2.7.3 et Remarque 2.7.4]{Pey:ecoledete}).
La fonction z\^eta des hauteurs associ\'ee \`a un tel ouvert $U$  est d\'efinie par
$$
Z_{\cc,U,h_{\ecL}}(T)=\sum_{x\in U(K)}\,T^{\,h_{\ecL}(x)}.
$$
Si $V$ n'est pas de type g\'en\'eral, 
le comportement analytique attendu de cette s\'erie est 
d\'ecrit par la version g\'eom\'etrique des conjectures de Manin et al.

\subsection{Le cas motivique}\label{subsec:cas_motiv}

Le corps $k$ est supposé quelconque.
On peut toujours d\'efinir une fonction hauteur $h_{\ecL}\,:\,V(K)\to \Z$ 
par la formule $h_{\ecL}(x)=\deg\left(x^{\ast}\ecL\right)$.
Si $L$ est une extension de $k$, on note $K_L$ le corps des fonctions 
de la courbe $\cc\times_k L$. Si $s$ est un point d'un sch\'ema $S$, on notera $\kappa_s$ le corps r\'esiduel en $s$.
Si $\varphi$ est un morphisme de source un $S$-sch\'ema, 
on notera $\varphi_s$ le morphisme déduit de $\varphi$ par le
changement de base $\Spec(\kappa_s)\to S$.

On note $\HOM^{\ecL,d}_{k}(\cc,V)$ 
le foncteur qui \`a un $k$-sch\'ema $S$
associe
$$
\left\{
\varphi\in \Hom_k(\cc_S,V),
\quad
\forall \,s\in S,\quad 
\deg\left(\varphi_s^{\ast}(\ecL)\right)=d
\right\}.
$$
Pour tout ouvert $U$ de $V$, on note $\HOM^{U,\ecL,d}_{k}(\cc,V)$
le foncteur qui à un $k$-sch\'ema $S$
associe
$$
\left\{
\varphi\in \HOM^{\ecL,d}_{k}(\cc,V)(S),
\quad
\forall \,s\in S,\quad 
\varphi_s\in U(K_{\kappa_s})
\right\}.
$$

\begin{lemme}\label{lm:repr}
On suppose que la classe de $\ecL$ est \`a l'int\'erieur du c\^one
effectif de $V$.
Alors il existe un ouvert $U_0$ non vide de $V$ tel que pour tout ouvert
$U$ de $U_0$ et pour tout $d\geq 1$, le foncteur $\HOM^{U,\ecL,d}_{k}(\cc,V)$
est repr\'esentable par un $k$-sch\'ema quasi-projectif.
\end{lemme}
\begin{proof}[D\'emonstration]
Supposons $\ecL$ ample. D'apr\`es \cite[4.c]{Gro:Hilb},
$\HOM^{\ecL,d}_{k}(\cc,V)$ est repr\'esentable par un $k$-sch\'ema quasi-projectif.
Or, pour tout ouvert $U$ de $V$,  $\HOM^{U,\ecL,d}_{k}(\cc,V)$ est un sous-foncteur ouvert
de $\HOM^{\ecL,d}_{k}(\cc,V)$. On en déduit le r\'esultat quand $\ecL$
est ample.

Dans le cas g\'en\'eral, $V$ \'etant projective, on peut fixer un fibr\'e en droite tr\`es ample
$\ecL_0$ sur $V$. Comme la classe de $\ecL$ est \`a l'int\'erieur du c\^one effectif, 
il existe un entier $N\geq 1$ et un fibr\'e en droites effectif $\ecL_1$ tel que
$\ecL^{\otimes N}=\ecL_0\otimes \ecL_1$.
Soit $U_0$ le compl\'ementaire des points-base de $\ecL_1$ et $U$ un ouvert contenu dans $U_0$.
Pour toute extension $L$ de $k$ et tout \'el\'ement $\varphi$ de
$U(K_L)$ (\ie  tout morphisme
$\varphi\,:\,\cc_L\to V$ dont l'image rencontre $U$), on a alors $\deg(\varphi^{\ast}(\ecL_1))\geq 0$.
Ainsi si un tel $\varphi$ v\'erifie $\deg(\varphi^{\ast}(\ecL))=d$ on a $\deg(\varphi^{\ast}(\ecL_0))\leq \frac{d}{N}.$
Donc le foncteur $\HOM^{U,\ecL,d}_{k}(\cc,V)$ s'identifie \`a un sous-foncteur ouvert 
du foncteur $\coprodu{d'\leq \frac{d}{N}} \HOM^{\ecL_0,d'}_{k}(\cc,V)$. 
Ce dernier foncteur \'etant repr\'esentable par un $k$-sch\'ema quasi-projectif, on a le r\'esultat.
\end{proof}
\begin{rem}
Hormis l'utilisation du th\'eor\`eme de repr\'esentabilit\'e de Grothendieck, la preuve
de ce lemme est formellement la m\^eme que celle de son analogue \og classique\fg.
Notons aussi que si $U$ est un ouvert de $V$ dont le groupe de Picard est trivial,
$U$ est n\'ecessairement contenu dans l'ouvert $U_0$ de la d\'emonstration, en 
particulier $\HOM^{U,\ecL,d}_{k}(\cc,V)$ est repr\'esentable par un sch\'ema quasi-projectif.
\end{rem}

On suppose désormais que la classe de $\ecL$ est \`a l'int\'erieur du c\^one effectif.
Pour tout ouvert $U$ assez petit de $V$ et tout entier $d\geq 0$, 
on note $U_{\ecL,d}$ le sch\'ema quasi-projectif repr\'esentant le foncteur $\HOM^{U,\ecL,d}_{k}(\cc,V)$.
On souhaite notamment \'etudier le \og comportement asympotique\fg\ de
$U_{\ecL,d}$ quand $d$ tend vers $+\infty$.
Pour cela, on consid\`erera notamment 
l'\'el\'ement de $\mk[[T]]$ d\'efini par 
$$
Z^{\text{mot}}_{\cc,U,h_{\ecL}}(T)
\eqdef
\sum_{d\geq 0}\,\left[U_{\ecL,d}\right]\,T^{\,d}.
$$
Si $k$ est fini, la s\'erie 
$Z_{\cc,U,h_{\ecL}}(T)$
est l'image de $Z^{\text{mot}}_{\cc,U,h_{\ecL}}$ par le morphisme $\#_k$.
Si $k$ est de caract\'eristique z\'ero, 
l'image  dans $\mkc[[T]]$ 
de $Z^{\text{mot}}_{\cc,U,h_{\ecL}}(T)$
par le morphisme $\chi$ 
sera
not\'ee
$
Z^{\chi}_{\cc,U,h_{\ecL}}(T)
$
.

Il est naturel de se demander s'il n'existe pas pour 
les séries $Z^{\text{mot}}_{\cc,U,h_{\ecL}}$ ou $Z^{\chi}_{\cc,U,h_{\ecL}}$ 
des analogues motiviques des r\'esultats obtenus ou conjectur\'es pour 
la série $Z_{\cc,U,h_{\ecL}}$. On pr\'ecise cette question dans un cas
particulier \`a la section suivante.

\subsection{Un analogue motivique de la conjecture de Manin}\label{subsec:manin_mot}

On suppose d\'esormais que la vari\'et\'e $V$ v\'erifie 
les hypoth\`eses suivantes :
\begin{hyp}\label{hyp:fano}
\begin{enumerate}
\item
La classe du faisceau anticanonique de $V$ est 
situ\'ee \`a l'int\'erieur du c\^one effectif.
\item
L'ensemble $V(K)$ est Zariski dense. 
En d'autres termes,
pour tout ouvert non vide $U$ de $V$, 
il existe un $k$-morphisme
$\cc \to V$ dont $U$ rencontre l'image.
\end{enumerate}
\end{hyp}
Dans tout ce qui suit,
on note $\ecL_0$ le faisceau anticanonique de $V$
et on consid\`ere la hauteur $h_{\ecL_0}$ associ\'ee \`a $\ecL_0$, not\'ee $h_0$.
On note aussi $U_{0,d}=U_{\ecL_0,d}$. 

Supposons le corps $k$ fini de cardinal $q$. Une partie
de la version g\'eom\'etrique de la conjecture de Manin peut alors 
se traduire par les deux questions suivantes.
\begin{question}\label{qu:conj_man_clas_1}
Existe-t-il un ouvert non vide $U$
de $V$
tel que la s\'erie 
$Z_{\cc,U,h_{0}}(q^{-s})$
converge absolument pour $\Re(s)>1$ ?
\end{question}
\begin{question}\label{qu:conj_man_clas_2}
Existe-t-il un ouvert non vide $U$ 
satisfaisant les exigences de la question \ref{qu:conj_man_clas_1}
et un $\varepsilon>0$ tel que la fonction holomorphe sur $\Re(s)>1$ 
d\'efinie par $Z_{\cc,U,h_{0}}(q^{-s})$
se prolonge en une fonction m\'eromorphe sur $\Re(s)>1-\varepsilon$
avec un p\^ole d'ordre $\rg(\NS(V))$ en $s=1$ ?
\end{question}

Supposons \`a pr\'esent $k$ quelconque.
La question qui suit est un analogue motivique na\"\i f de la question \ref{qu:conj_man_clas_1}.
\begin{question}\label{qu:conj_manin_mot_1_naif}
Existe-t-il un ouvert non vide $U$ de $V$ tel que, pour tout entier $\kappa\geq 2$, 
la s\'erie $Z^{\text{mot}}_{\cc,U,h_{0}}(\LL^{-\kappa})$ converge dans $\mhat$ ? 
\end{question}
Cette question peut se reformuler ainsi :
existe-t-il un ouvert non vide $U$ tel que pour tout entier $\kappa\geq 2$ on a 
$$
\lim_{d\to \infty} \dim(U_{0,d})-\kappa\,d = -\infty\quad ?
$$
On voit ainsi que pour obtenir un analogue plus fid\`ele de la
question \ref{qu:conj_man_clas_1}, cette dernière condition
devrait être exigée pour tout réel $\kappa>1$.
\begin{question}\label{qu:conj_manin_mot_1}
Existe-t-il un ouvert non vide $U$ de $V$ tel que
$$
\overline{\lim_{d\to \infty}} \,\, \,\,\frac{\dim(U_{0,d})}{d}\leq 1 \quad ?
$$
\end{question}

Pour toute vari\'et\'e $X$, on note $\rho(X)$ le nombre de composantes g\'eom\'etriques irr\'eductibles de dimension maximale de $X$.
Au vu des estimations de Lang-Weil et 
du comportement asympotique de $\card{U_{0,d}(k)}$ donn\'ee par des th\'eor\`emes taub\'eriens
standards lorsque la r\'eponse \`a la question \ref{qu:conj_man_clas_2} est positive, 
un analogue naturel de la question \ref{qu:conj_man_clas_2} est la
question suivante.
\begin{question}\label{qu:conj_manin_mot_irr}
Existe-t-il un ouvert non vide $U$ de $V$ v\'erifiant
$$
\overline{\lim_{d\to \infty}} \,\, \,\,\frac{\dim(U_{0,d})}{d}=1 
$$
et
$$
\overline{\lim_{d\to \infty}} \,\, \,\,\frac{\log(\rho(U_{0,d}))}{\log(d)}=\rg(\NS(V)) \quad ?
$$
\end{question}
Un autre analogue possible est la question suivante.
\begin{question}\label{qu:conj_manin_mot_2}
Existe-t-il un ouvert non vide $U$ de $V$ v\'erifiant
$$
\overline{\lim_{d\to \infty}} \,\, \,\,\frac{\dim(U_{0,d})}{d}=1 
$$
et tel que la s\'erie d\'efinie par le produit 
$(1-\LL\,T)^{\rg(\NS(V))}\,Z^{\text{mot}}_{\cc,U,h_{0}}(T)$
converge pour $T=\LL^{-1}$?
\end{question}
On note $\overline{k}$ une cl\^oture alg\'ebrique de $k$, et $k^{\text{s\'ep}}$
la cl\^oture s\'eparable de $k$ dans $\overline{k}$.
D\'esormais, en plus des hypoth\`eses \ref{hyp:fano}, on suppose 
que la vari\'et\'e $V$ v\'erifie les hypoth\`eses suivantes (\cf  les hypoth\`eses
du paragraphe 2.1 de \cite{Pey:var_drap} o\`u l'on notera que la vari\'et\'e $V$
est suppos\'ee d\'efinie seulement sur le corps des fonctions de $\cc$):

\begin{hyp}\label{hyp:ceff}
\begin{enumerate}
\item
Les groupes de cohomologie $H^1(V,\ecO_V)$ et $H^2(V,\ecO_V)$ 
sont nuls.
\item
Le groupe $\Pic(V\times_k k^{\text{s\'ep}})$ est libre de rang fini 
et co\"\i ncide avec
$\Pic(V\times_k \overline{k})$.
\item\label{item:galabstriv}
L'action du groupe de Galois absolu sur $\Pic(V\times_k k^{\text{s\'ep}})$ est triviale.
\item
Le c\^one effectif de $V$ (not\'e $\ceff(V)$) est polyédral rationnel.
\item
Si $l$ est un nombre premier distinct de la caract\'eristique de $k$, la partie $l$-primaire
de $Br(V\times_k \overline{k})$ est finie.
\end{enumerate}
\end{hyp}
Le fait d'avoir suppos\'e la vari\'et\'e $V$ d\'efinie sur le corps des constantes de $\cc$
et l'hypothèse \ref{item:galabstriv} ci-dessus signifie que nous nous pla\c cons dans le cas arithm\'etiquement le plus simple.
Notons que les vari\'et\'es toriques d\'eploy\'ees v\'erifient les hypoth\`eses \ref{hyp:fano} et \ref{hyp:ceff}.

On pose 
$$
Z_{\ceff(V),\ecL_0}(T)
\eqdef
\sum_{y\in \ceff(V)^{\vee}\cap \Pic(V)^{\vee}}\,T^{\,\acc{y}{\ecL_0}}.
$$
En \'ecrivant $\ceff(V)^{\vee}$ comme le support d'un \'eventail r\'egulier,
on voit que $Z_{\ceff(V),\ecL_0}(T)$ est une fraction rationnelle \`a coefficients 
rationnels en $T$,
dont $1$ est un p\^ole d'ordre $\rg(\Pic(V))$.
On pose
\begin{equation*}
\alpha^{\ast}(V)
\eqdef
\lim_{
z\to 1
}
(z-1)^{\,\rg(\Pic(V))}\,Z_{\ceff(V),\ecL_0}(z).
\end{equation*}
C'est un nombre rationnel non nul.

Supposons le corps $k$ fini de cardinal $q$. 
Une version g\'eom\'etrique de la conjecture de Manin 
raffin\'ee par Peyre
peut alors s'\'enoncer ainsi :
\begin{question}\label{qu:conj_man_clas_3}
Si $U$ est un ouvert satisfaisant les exigences
des questions \ref{qu:conj_man_clas_1}
et
\ref{qu:conj_man_clas_2},
a-t-on
\begin{multline*}
\lim_{s\rightarrow 1}(s-1)^{\rg(\Pic(V))}
Z_{\cc,U,h_{0}}(q^{-s})
=\\
\alpha^{\ast}(V)
\,
q^{\,(1-g_{\cc})\,\dim(V)}
\,
\left(\Res_{s=1}Z_{\cc}(q^{-s})\right)^{\rg(\Pic(V))}
\times\,
\prod_{v\in \cc^{(0)}}
(1-\card{k_v}^{-1})^{\rg(\Pic(V))}
\,\,
\frac
{\card{V(k_v)}}
{\card{k_v}^{\dim(V)}}\quad ?
\end{multline*}
\end{question}

La convergence du produit eul\'erien apparaissant dans l'expression 
conjecturale du terme principal en $s=1$ est démontr\'ee par Peyre
en utilisant les conjectures de Weil d\'emontr\'ees par Deligne. Ainsi, 
sous les hypoth\`eses \'enonc\'ees ci-dessus,
Peyre montre que l'analogue de l'hypoth\`ese de Riemann 
entra\^\i ne l'estimation asymptotique 
\begin{equation}\label{eq:asymp_vkv}
\card{V(k_v)}
=
\card{k_v}^{\dim V}+
\rg(\Pic(V)\,\card{k_v}^{\dim V-1}
+
\underset{\deg(v)\to \infty}{\ecO}
\left(\card{k_v}^{\dim V-\frac{3}{2}}\right).
\end{equation}
La convergence en découle.

La r\'eponse \`a la question \ref{qu:conj_man_clas_2} est positive
dans le cas des vari\'et\'es toriques (\cite{Bou:vtetor}) et 
dans le cas des vari\'et\'es de drapeaux (\cite{Pey:var_drap}).

Au vu de \eqref{eq:asymp_vkv}, on peut se poser la question suivante :

\begin{question}\label{qu:asymp_vkv_mot}
A-t-on 
$$
\forall d\geq 1,\quad
\Phi_d(V)-\LL^{d\,\dim V}-
\rg(\Pic(V)\,\LL^{d\,(\dim V-1)}
\in  
\fil^{\,d(\frac{3}{2}-\dim(V))}\mk
\quad ?
$$
\end{question}

Nous donnons \`a pr\'esent un analogue motivique 
de la question \ref{qu:conj_man_clas_2}. 

\begin{question}\label{qu:conj_manin_mot_3}
On suppose que les r\'eponses aux questions 
\ref{qu:asymp_vkv_mot}
et \ref{qu:conj_manin_mot_2}
sont positives.
Soit $U$ un ouvert satisfaisant les exigences de la question
\ref{qu:conj_manin_mot_2}. 
La s\'erie
$$
(1-\LL\,T)^{\rg(\Pic(V))}\,Z^{\text{mot}}_{\cc,U,h_{0}}
$$
converge-t-elle dans $\mhatq$ en $T=\LL^{-1}$ 
vers 
\begin{multline*}
\alpha^{\ast}(V)
\,
\LL^{\,(1-g_{\cc})\,\dim(V)}
\,
\left(
\left[
(1-\LL\,T)\,\zetacm(T)
\right]
(\LL^{-1})
\right)^{\,\rg(\Pic(V))}\\
\\
\times\,
\exp
\left(
\sum_{n\geq 1}
\Psi_n(\cc)
\log
\left(
(1-\LL^{\,-n})^{\rg(\Pic(V))}
\frac{\Phi_n(V)}{\LL^{\,-n\,\dim(V)}}
\right)
\right) \quad ?
\end{multline*}
\end{question}
La convergence du produit eul\'erien motivique est assur\'ee par la r\'eponse 
positive \`a la question \ref{qu:asymp_vkv_mot}.

Si le corps $k$ est fini, 
une r\'eponse positive \`a la question \ref{qu:conj_manin_mot_3}  
ne permet pas a priori d'obtenir par sp\'ecialisation via $\#_k$ des r\'esultats sur 
la fonction z\^eta des hauteurs usuelle. 
Le terme principal motivique propos\'e ici se sp\'ecialise formellement sur le terme principal 
de la fonction z\^eta des hauteurs usuelle, 
mais la fonction $\#_k$ n'est pas d\'efinie sur $\mhat$.

Une telle sp\'ecialisation sera cependant fructueuse dans le cas o\`u la
s\'erie $Z^{\text{mot}}_{\cc,U,h_{0}}(T)$ est une fraction rationnelle en $T$.
Cette situation se produit pour les espaces projectifs et, 
comme on le montre ci-dessous, 
pour les surfaces de Hirzebruch si $\cc=\PP^1$. 
Cependant, pour une vari\'et\'e torique d\'eploy\'ee quelconque,
nous suspectons qu'en g\'en\'eral $Z_{\cc,U,h_{0}}(T)$ 
(et donc $Z^{\text{mot}}_{\cc,U,h_{0}}(T)$) 
n'est pas une fraction rationnelle en $T$.

Supposons \`a pr\'esent $k$ de caract\'eristique z\'ero.
Une variante \'evidente des questions 
\ref{qu:asymp_vkv_mot}
 et
\ref{qu:conj_manin_mot_3} est donn\'ee par les \'enonc\'es suivants.

\begin{question}\label{qu:asymp_vkv_mot_chi}
A-t-on
$$
\forall d\geq 1,\quad 
\Phi^{\chi}_d(V)-\LL^{d\,\dim V}-
\rg(\Pic(V)\,\LL^{d\,(\dim V-1)}
\in  
\fil^{\,d(\frac{3}{2}-\dim(V))}\mkcq
\quad ?
$$
\end{question}
\begin{question}\label{qu:conj_manin_mot_chi_3}
On suppose que les r\'eponses aux questions 
\ref{qu:asymp_vkv_mot_chi}
et \ref{qu:conj_manin_mot_1}
sont positives.
Soit $U$ un ouvert satisfaisant les exigences de la question
\ref{qu:conj_manin_mot_1}. 
La s\'erie
$$
(1-\LL\,T)^{\rg(\Pic(V))}\,Z^{\text{mot},\chi}_{\cc,U,h_{0}}(T)
$$
converge-t-elle dans $\mhatcq$ en $T=\LL^{-1}$ 
vers
\begin{multline*}
\alpha^{\ast}(V)
\,
\LL^{\,(1-g_{\cc})\dim(V)}\,
\left(
\left[
(1-\LL\,T)\,\zetacmc(T)
\right]
\left(
\LL^{-1}
\right)
\right)^{\,\rg(\Pic(V))}\\
\\
\times
\,
\exp
\left(
\sum_{n\geq 1}
\Psi^{\chi}_n(\cc)
\log
\left(
(1-\LL^{\,-n})^{\rg(\Pic(V))}
\frac
{\Phi^{\chi}_n(V)}
{\LL^{\,-n\,\dim(X)}}
\right)
\right)\quad ?
\end{multline*}
\end{question}

Dans la section suivante,\textbf{\textit{
en supposant que la courbe $\cc$ est rationelle,
}}
nous démontrons les r\'esultats suivants :
la r\'eponse \`a la question \ref{qu:conj_manin_mot_1} est 
positive si $V$ est une vari\'et\'e torique d\'eploy\'ee, 
en prenant pour $U$ l'orbite ouverte (corollaire \ref{cor:rep_pos_ques_1}) ;
la r\'eponse \`a la question \ref{qu:conj_manin_mot_3} est positive 
si $V$ est une surface de Hirzebruch, en prenant pour $U$ l'orbite
ouverte  (on montre en fait un résultat plus fort, à savoir  
le point \ref{item:1:theoprinc} du théorème
\ref{thm:theoprinc}, dont l'énoncé est repris à la sous-section \ref{subsec:hirz}).
Enfin, si on suppose $k$ de caract\'eristique z\'ero,
nous montrons que la r\'eponse \`a la question \ref{qu:conj_manin_mot_chi_3} est 
positive si $V$ est une vari\'et\'e torique d\'eploy\'ee, en prenant pour $U$
l'orbite ouverte (à savoir le point \ref{item:2:theoprinc} du théorème
\ref{thm:theoprinc}).

Les preuves du corollaire \ref{cor:rep_pos_ques_1} et du th\'eor\`eme \ref{thm:theoprinc}
 s'appuient sur le lemme de param\'etrisation \ref{lm:param_mot}.
La preuve du point \ref{item:2:theoprinc} du théorème
\ref{thm:theoprinc}  utilise en outre le th\'eor\`eme \ref{thm:prod_eul_mot_zmubcx}
et le corollaire \ref{cor:muxb}. 
Une r\'eponse positive \`a la question \ref{qu:prod_eul_mot_zmubmx} 
permettrait de montrer, par les mêmes méthodes, que la r\'eponse \`a la question \ref{qu:conj_manin_mot_3} est 
positive pour toute vari\'et\'e torique d\'eploy\'ee. 

\begin{rem}
Pour une réponse partielle à la question 
\ref{qu:asymp_vkv_mot_chi} lorque la variété $V$ n'est plus  nécessairement
torique, on pourra consulter \cite{Bou:tammot}.
\end{rem}

\section{Démonstration des r\'esultats annonc\'es}\label{sec:calcul}

\subsection{Quelques rappels sur les vari\'et\'es toriques}\label{subsec:mot_vtetor}

Nous nous contentons de citer les r\'esultats qui nous seront utiles, ce qui nous permet de fixer quelques notations.
Nous renvoyons le lecteur aux r\'ef\'erences classiques sur le sujet 
(par exemple \cite{Oda:conv}, \cite{Ful:toric}, \cite{Ew}) pour plus de d\'etails.  
Soit $r\geq 1$ un entier et $U\isom\G_{m}^r$ un tore d\'eploy\'e de dimension $r$ d\'efini sur $k$.
Soit $\cara{U}$ son groupe des caract\`eres et $\cocara{U}$ son groupe des cocaract\`eres.
On note $\acc{.}{.}$ l'accouplement naturel entre ces deux $\Z$-modules.
Soit $k$ un corps.
\`A tout \'eventail $\Sigma$ de $\cocara{U}$ est associ\'e une
$k$-vari\'et\'e normale irr\'eductible $\xs$, 
munie d'une action de $U$ et poss\'edant une orbite ouverte isomorphe \`a $U$, 
en d'autres termes une vari\'et\'e torique d\'eploy\'ee
d\'efinie sur $k$. On suppose désormais l'\'eventail $\Sigma$ projectif et r\'egulier. 
La vari\'et\'e $\xs$  est alors projective et lisse.

Pour $\alpha\in\Sigma(1)$ nous notons $\roa\in \cocara{U}$ le g\'en\'erateur de $\alpha$, 
$\mathfrak{D}_{\alpha}$ le diviseur $U$-invariant associ\'e 
et $\cD_{\alpha}$ sa classe dans le groupe de Picard de $\xs$. 
Le diviseur $\sum_{\alpha\in\Sigma(1)}\mathfrak{D}_{\alpha}$ est alors
un diviseur anticanonique.
Le morphisme qui à un \'el\'ement $m$ de $\cara{U}$ associe $(\acc{m}{\roa})_{\alpha\in\Sigma(1)}$
induit une suite exacte 
\begin{equation}\label{eq:suite_exacte_pic_tor}
0\longto \cara{U} \longto 
\oplusu{\alpha\in \Sigma(1)} 
\Z\,\mathfrak{D}_{\alpha}
\longto \Pic(\xs)\longto 0.
\end{equation}

\subsection{Param\'etrisation des morphismes}\label{subsec:param_mot}

On conserve les notations de la sous-section \ref{subsec:mot_vtetor}. 
Soit
$$
\ts=
\A_k^{\Sigma(1)}
\setminus
\underset{\sigma\in \Sigma}{\bigcap}
\left\{
\prod_{\alpha\notin \sigma}x_{\alpha}=0\right\}.
$$
La suite exacte \eqref{eq:suite_exacte_pic_tor}
induit une suite exacte de tores
\begin{equation}\label{eq:exsqtores}
0\longto T_{\Pic(\xs)}\longto \G_{m}^{\Sigma(1)}\overset{\pi}{\longto} U\longto 0.
\end{equation}
L'action diagonale de  $\G_m^{\Sigma(1)}$
sur $\ts$ induit par restriction une action de 
$T_{\Pic(\xs)}$ sur $\ts$. 
D'apr\`es \cite{Cox:hom_coo_ring}, le morphisme  
$\pi$
s'étend en un morphisme équivariant $\ts\longto \xs$ 
qui est un quotient g\'eom\'etrique de $\ts$ par 
$T_{\Pic(\xs)}$. Un tel quotient fournit ainsi un syst\`eme de coordonn\'ees $T_{\Pic(\xs)}$-homog\`enes sur $\xs$.
Ceci g\'en\'eralise les coordonnées homogènes classiques de l'espace projectif. 
Signalons que le morphisme $\ts\to \xs$ fait de $\ts$ un torseur universel au-dessus de $\xs$
(\cf  \cite[proposition 8.5]{Sal:tammes} et \cite[Appendix]{Mad:var_tor}). 
Ces coordonnées homogènes permettent \`a Cox 
de donner dans \cite{Cox:funct} une description du foncteur des points
de $\xs$, que nous rappelons ci-dessous.
Elle g\'en\'eralise la description bien connue du foncteur des points 
de l'espace projectif.
Il en découle pour tout $k$-schéma $S$ une description simple des
$k$-morphismes de $\PP^1_S$ dans $\xs$, 
que nous allons utiliser
 pour expliciter
le sch\'ema quasi-projectif qui repr\'esente le foncteur $\HOM^{U,\ecL_0,d}_{k}(\PP^1,\xs)$
d\'efini dans la sous-section \ref{subsec:cas_motiv}.

\begin{defi}
Soit $S$ un $k$-sch\'ema.
Une \termin{$\Sigma$-collection sur $S$} est la donn\'ee pour tout 
$\alpha\in \Sigma(1)$ d'un fibr\'e en droites $\ecLa$
sur $S$ et d'une section globale $v_{\alpha}$ de $\ecLa$
ainsi que d'une famille $(c_m)_{m\in \cara{U}}$ d'isomorphismes
$$
c_m\,:\,\otimesu{\alpha} \,\ecLa^{\acc{m}{\roa}}\longisom \ecO_{S},
$$
ces donn\'ees \'etant astreintes \`a v\'erifier les conditions suivantes :
\begin{enumerate}
\item\label{item:1:defi:collec}
pour tous $m,m'$ dans $\cara{U}$, on a
$
c_m\otimes c_{m'}=c_{m+m'}\quad;
$
\item\label{item:2:defi:collec}
pour tout $\alpha\in \Sigma(1)$, la section 
$v_{\alpha}$ induit un morphisme
$\ecO_S\to \ecLa$
et par dualit\'e un morphisme
$\ecLa^{-1}\to \ecO_S$ ;
le morphisme induit
$$
\bigoplus_{\sigma\in \Sigma} 
\underset{\alpha\notin\sigma}{\otimes}
\ecLa^{-1}\longto \ecO_S
$$
est surjectif.
\end{enumerate}
Un isomorphisme entre deux $\Sigma$-collections 
$\left((\ecLa,v_{\alpha}),(c_m)\right)$ et
$\left((\ecL'_{\alpha},s'_{\alpha}),(c'_m)\right)$ 
est une famille d'isomorphismes $\ecLa\isom \ecLa'$
envoyant $s_{\alpha}$ sur $s'_{\alpha}$ et $c_m$ sur $c'_m$.
\end{defi}
Cox démontre alors le th\'eor\`eme suivant (\cite[Theorem 1.1]{Cox:funct}).
\begin{thm}\label{thm:cox}
Le foncteur qui \`a un $k$-sch\'ema $S$ associe l'ensemble des classes
d'isomorphisme de $\Sigma$-collections
sur $S$  est repr\'esent\'e par la vari\'et\'e torique $\xs$.
\end{thm}
Tr\`es grossi\`erement, l'id\'ee de la démonstration est la suivante : 
\`a la $\Sigma$-collection $\left((\ecLa,v_{\alpha}),(c_m)\right)$
on fait correspondre le morphisme qui \`a  $s\in S$ associe le 
\og point de coordonn\'ees
homog\`enes $(v_{\alpha}(s))$ \fg.
La condition \ref{item:1:defi:collec} et la suite exacte \eqref{eq:suite_exacte_pic_tor} montrent
que le $\Sigma(1)$-uple
$(v_{\alpha}(s))$ est bien d\'efini modulo l'action de $T_{\Pic(\xs)}$.
La condition \ref{item:2:defi:collec} assure que  $(v_{\alpha}(s))$ est dans $\ts$.
On obtient ainsi un morphisme $S\to \xs$. R\'eciproquement, \`a un morphisme $\pi\,:\,S\to \xs$, 
on associe la $\Sigma$-collection $(\pi^{\ast}\ecO(\mathfrak{D}_{\alpha}),\pi^{\ast}v_{\alpha},\pi^{\ast}c_m)$
o\`u $v_{\alpha}$ est la section canonique de $\ecO(\mathfrak{D}_{\alpha})$ et les trivialisations $c_m$
sont donn\'ees par la suite exacte \eqref{eq:suite_exacte_pic_tor}.

\begin{nota}
On note $\N^{\Sigma(1)}_{\bigast}$
le sous-mono\"\i de de $\N^{\Sigma(1)}$ 
constitu\'e des \'el\'ements $\dd$  v\'erifiant
\begin{equation}\label{eq:bigast}
\forall \,m\in \cara{U},
\quad
\sum_{\alpha\in \Sigma(1)}\acc{m}{\roa} \,\da=0.
\end{equation}
En d'autres termes, 
si on identifie 
$\Pic(\xs)^{\vee}$ 
à un sous-groupe de $\Z^{\Sigma(1)}$
via le dual de la suite exacte
\eqref{eq:suite_exacte_pic_tor},
$\N^{\Sigma(1)}_{\bigast}$
est l'intersection de 
 $\N^{\Sigma(1)}$ et de $\Pic(\xs)^{\vee}$. 
\end{nota}
\begin{nota}
Pour  $\bd\in \Z^{\Sigma(1)}$, on note $\HOM^{U,\bd}_{k}(\PP^1,\xs)$
le foncteur qui \`a un $k$-sch\'ema $S$ associe
l'ensemble des \'el\'ements 
$\varphi$ de $\Hom_k(\PP^1_S,X_{\Sigma})$
tels que, pour tout $s\in S$,
$\varphi_s\in U(\kappa_s(t))$ (\ie  l'image
de $\varphi_s$ rencontre $U$)
et pour tout $\alpha\in \Sigma(1)$,
$
\deg\left(\varphi_s^{\ast}(\mathfrak{D}_{\alpha})\right)=\da.
$
\end{nota}
Si $L$ est une extension de $k$, on a une suite exacte 
$$
U(L(t))\longto\Hom(\cara{U},\Div(\PP^1_L))\overset{\deg}{\longto}\Hom(\cara{U},\Z)\rightarrow 0.
$$
Si $\varphi$ est un \'el\'ement de $U(L(t))$, son image dans $\Hom(\cara{U},\Div(\PP^1_L))$ est 
$$
m\longmapsto \sum_{\alpha\in \Sigma(1)}\acc{m}{\roa}\,\varphi^{\ast}({\mathfrak D}_{\alpha}).
$$
La suite exacte ci-desus montre alors que $\left(\deg(\varphi^{\ast}({\mathfrak D}_{\alpha})\right)$ 
est un \'el\'ement de $\N^{\Sigma(1)}_{\bigast}$. 
Ainsi si $\bd$ n'appartient pas à $\N^{\Sigma(1)}_{\bigast}$,
le foncteur $\HOM^{U,\bd}_{k}(\PP^1,\xs)$ est vide .

Le but de ce qui suit utiliser le th\'eor\`eme \ref{thm:cox} 
est d'expliciter une vari\'et\'e repr\'esentant
$\HOM^{U,\bd}_{k}(\PP^1,\xs)$ lorsque 
$\bd$ appartient à $\N^{\Sigma(1)}_{\bigast}$.

Si $S$ est un $k$-sch\'ema, on note $p_{1,S}$ et $p_{2,S}$ les projections 
de $\PP^1_S$
vers $\PP^1_k$ et $S$ respectivement. 
On a alors un isomorphisme 
$$
p_{1,S}^{\ast}+p_{2,S}^{\ast}\,:\,
\Z\oplus \Pic(S)\isom \Pic(\PP^1_S)
.$$

\begin{defi}
Soit $S$ un $k$-sch\'ema et $\bd\in \N^{\Sigma(1)}_{\bigast}$.
Une \termin{$(\PP^1,\Sigma,\bd)$-collection non d\'eg\'en\'er\'ee sur $S$} est la donn\'ee 
d'une $\Sigma$-collection $((\ecMa,u_{\alpha}),(c_m))$ sur $\PP^1_S$
telle que pour tout $\alpha$, $u_{\alpha}$ est non nulle et telle que la projection
$\Pic(\PP^1_S)\to \Z$ envoie la classe de $\ecMa$ sur $\da$.
Un isomorphisme entre deux $(\PP^1,\Sigma,\bd)$-collections non d\'eg\'en\'er\'ees sur $S$ 
est un isomorphisme entre ces deux objets en tant que $\Sigma$-collections sur $\PP^1_S$.
\end{defi}
Le foncteur qui à $S$ associe l'ensemble des classes d'isomorphisme de
$(\PP^1,\Sigma,\bd)$-collections non d\'eg\'en\'er\'ees sur $S$ 
s'identifie à un sous-foncteur de $\HOM^{\bd}_{k}(\PP^1,\xs)$.
Un examen des arguments de la d\'emonstration du th\'eor\`eme \ref{thm:cox}
permet de montrer le lemme suivant 
(le point important est que l'image r\'eciproque de $U$ dans $\ts$ est l'ouvert $\prod x_{\alpha}\neq 0$).
\begin{lemme}\label{lm:p1_s_da}
\`A toute $(\PP^1,\Sigma,\bd)$-collection 
non d\'eg\'en\'er\'ee sur un $k$-sch\'ema $S$
on peut associer de mani\`ere fonctorielle en $S$ un $k$-morphisme de $\PP^1_S$ vers $\xs$.

Ceci induit un isomorphisme entre le foncteur qui \`a un $k$-sch\'ema $S$ 
associe l'ensemble des classes d'isomorphisme
$(\PP^1,\Sigma,\bd)$-collections non d\'eg\'en\'er\'ees sur $S$ 
et le foncteur $\HOM^{U,\bd}_{k}(\PP^1,\xs)$.
\end{lemme}
\begin{defi}
Soit $S$ un $k$-sch\'ema et $\bd$ un \'el\'ement de $\N^{\Sigma(1)}_{\bigast}$.
Une \termin{$(\Sigma,\bd)$-collection du second type sur $S$} 
est la donn\'ee pour tout 
$\alpha\in \Sigma(1)$ d'un fibr\'e en droites $\ecLa$ sur $S$
et d'un $(\da+1)$-uple $(s_{\alpha,i})_{i=0,\dots,\da}$
de sections globales de $\ecLa$ qui engendrent $\ecLa$.

Un isomorphisme entre deux $(\Sigma,\bd)$-collections du second type
$\left(\ecLa,(s_{\alpha,i})\right)$ et
$\left(\ecL'_{\alpha},(s'_{\alpha,i})\right)$ 
est une famille d'isomorphismes $\ecLa\isom \ecL'_{\alpha}$
envoyant $s_{\alpha,i}$ sur $s'_{\alpha,i}$.
\end{defi}
\begin{defi}
Soit $\left(\ecLa,(s_{\alpha,i})\right)$
une $(\Sigma,\bd)$-collection du second type 
sur $S$ .
Soit $s\in S$. 
Pour tout $\alpha\in \Sigma(1)$, fixons un isomorphisme 
$
H^0(\kappa_s,s^{\ast}\ecLa)\isom \kappa_s.
$
On identifie alors 
l'image de $(s_{\alpha,i})$ dans $H^0(\kappa_s,s^{\ast}\ecLa)^{\da+1}$
\`a un polyn\^ome $P_{\alpha,s}$ homog\`ene en deux variables de degr\'e $\da$.

La $(\Sigma,\bd)$-collection  $\left(\ecLa,(s_{\alpha,i})\right)$
est dite \termin{non d\'eg\'en\'er\'ee} si elle v\'erifie la condition
suivante : pour tout $s\in S$, les polyn\^omes 
$$
\left(\prod_{\alpha\notin \sigma}P_{\alpha,s}\right)_{\sigma\in \Sigma}
$$
n'ont pas de z\'ero commun non trivial dans une cl\^oture alg\'ebrique de $\kappa_s$.

Un isomorphisme entre deux $(\Sigma,\bd)$-collections du second type
non d\'eg\'en\'er\'ees est un isomorphisme entre ces deux objets en tant que $(\Sigma,\bd)$-collections du second type.
\end{defi}
\begin{defi}
Soit $M$ un sous-module de $\Z^{\Sigma(1)}$ tel que le quotient de $\Z^{\Sigma(1)}$
par $M$ soit sans torsion. Une $(\Sigma,\bd)$-collection du second type  $M$-trivialisée
sur un $k$-sch\'ema $S$
est un couple 
$$
\left((\ecLa,(s_{\alpha,i})),(c_m)_{m\in M}\right)
$$
où $(\ecLa,(s_{\alpha,i}))$ est une $(\Sigma,\bd)$-collection du
second type sur $S$ et $(c_m)_{m\in M}$ est une famille d'isomorphismes
$$
c_m\,:\,\otimesu{\alpha} \,\ecLa^{m_{\alpha}}\longisom \ecO_{S},
$$
telle que, pour tout $m$ et $m'$ dans $M$, on a
$
c_m\otimes c_{m'}=c_{m+m'}
$.

Un isomorphisme entre deux $(\Sigma,\bd)$-collections du second type
$M$-trivialisées
$\left((\ecLa,(s_{\alpha,i}),(c_m)\right)$ et
$\left((\ecL'_{\alpha},(s'_{\alpha,i})),(c'_m)\right)$ 
est une famille d'isomorphismes $\ecLa\isom \ecL'_{\alpha}$
envoyant $s_{\alpha,i}$ sur $s'_{\alpha,i}$ et $c_m$ sur $c'_m$.
\end{defi}

Soit $\bd\in \N^{\Sigma(1)}_{\bigast}$.
L'action diagonale du tore  $\G_{m}^{\Sigma(1)}$ sur 
$\produ{\alpha}\left(\A_k^{\,\da+1}\setminus \{0\}\right)$
induit un $\G_{m}^{\Sigma(1)}$-torseur
$$
\produ{\alpha\in \Sigma(1)}\left(\A_k^{\,\da+1}\setminus \{0\}\right)
\longto 
\prod_{\alpha\in \Sigma(1)} \P_k^{\da}.
$$
On a d\'efini au d\'ebut de la section \ref{sec:fonc_mob_mot} 
un ouvert $\left(\PP^1\right)^{B_{\Sigma}}_{\bd}$  
de $\prod \P_k^{\da}$
(\cf  la sous-section \ref{subsec:mob_ev} pour la d\'efinition de l'ensemble $B_{\Sigma}$
associ\'e \`a l'\'eventail $\Sigma$). 
On note $\widetilde{\left(\PP^1\right)^{B_{\Sigma}}_{\bd}}$ 
l'image r\'eciproque de $\left(\PP^1\right)^{B_{\Sigma}}_{\bd}$
dans $\produ{\alpha\in \Sigma(1)}\left(\A_k^{\,\da+1}\setminus \{0\}\right)$.

\begin{lemme}
Soit $(P_{\alpha})$ une famille de polyn\^omes homog\`enes en deux variables 
\`a coefficients dans un corps $L$.
Les conditions suivantes sont \'equivalentes :
\begin{enumerate}
\item
pour tout $n_{\alpha}\in B_{\Sigma}$, les polyn\^omes
$
\left(P_{\alpha}\right)_{\alpha,\,n_{\alpha}=1}
$
n'ont pas de z\'ero commun non trivial dans une cl\^oture alg\'ebrique de $L$.
\item
les polyn\^omes 
$
\left(\prod_{\alpha\notin \sigma}P_{\alpha}\right)_{\sigma\in \Sigma}
$
n'ont pas de z\'ero commun non trivial dans une cl\^oture alg\'ebrique de $L$.
\end{enumerate}
\end{lemme}
\begin{proof}[D\'emonstration]
Supposons qu'il existe $(n_{\alpha})\in B_{\Sigma}$ et un z\'ero commun aux 
polynômes $\left(P_{\alpha}\right)_{\alpha,\,n_{\alpha}=1}$.
Par d\'efinition de $B_{\Sigma}$, pour tout c\^one $\sigma$ un tel
z\'ero est alors un z\'ero de 
$\prod_{\alpha\notin \sigma}P_{\alpha}$. 
R\'eciproquement, si les polyn\^omes 
$
\left(\prod_{\alpha\notin \sigma}P_{\alpha}\right)_{\sigma\in \Sigma}
$
ont un z\'ero commun $z$, on d\'efinit $n_{\alpha}\in \{0,1\}^{\Sigma(1)}$
par $n_{\alpha}=1$ si et seulement si $z$ est un z\'ero de $P_{\alpha}$. 
Alors $(n_{\alpha})\in \bs$ et $z$ est 
un z\'ero commun aux polynômes $\left(P_{\alpha}\right)_{\alpha,\,n_{\alpha}=1}$.
\end{proof}
Il est bien connu que le foncteur qui à un $k$-sch\'ema $S$ associe 
l'ensemble des classes d'isomorphisme de $(\Sigma,\bd)$-collections du second type sur $S$
 est repr\'esent\'e par $\prod \P_k^{\da}$.
De la d\'efinition de $\left(\PP^1\right)^{B_{\Sigma}}_{\bd}$  
et du lemme précédent on déduit le corollaire suivant.
\begin{cor}
L'ouvert $\left(\PP^1\right)^{B_{\Sigma}}_{\bd}$  
s'identifie au sous-foncteur ouvert de $\HOM(S,\prod \PP^{\da})$
qui \`a un $k$-sch\'ema $S$ associe l'ensemble des 
classes d'isomorphisme de
$(\Sigma,\bd)$-collections du second type sur $S$
non d\'eg\'en\'er\'ees.
\end{cor}
\begin{lemme}\label{lm:rep_quot}
Soit $M$ un sous-module de $\Z^{\Sigma(1)}$ tel que le quotient de $\Z^{\Sigma(1)}$
par $M$, not\'e $P$, soit sans torsion. Soit $T_P$ le sous-tore de $\G_m^{\Sigma(1)}$
associ\'e \`a $P$, $T_M$ le tore quotient de $\G_m^{\Sigma(1)}$
associ\'e \`a $M$ et $\pi\,:\,\G_m^{\Sigma(1)}\to T_M$ le morphisme quotient.

Le foncteur qui a un $k$-sch\'ema $S$ associe 
l'ensemble des classes d'isomorphisme 
$(\Sigma,\bd)$-collections du second type sur $S$
$M$-trivialisées (respectivement $M$-trivialisées non-d\'eg\'en\'er\'ees) 
est repr\'esent\'e 
par $\produ{\alpha\in \Sigma(1)}\left(\A_k^{\,\da+1}\setminus \{0\}\right)/T_P$
(respectivement par $\widetilde{\left(\PP^1\right)^{B_{\Sigma}}_{\bd}}/T_P$).
\end{lemme}
\begin{rem}
Le résultat est classique si $M=\Z^{\Sigma(1)}$.
\end{rem}
\begin{proof}[D\'emonstration]
Si $\ecL$ est un fibr\'e en droites sur une vari\'et\'e $X$, 
on note $\widetilde{\ecL}$ le $\G_m$-torseur au dessus de $X$ obtenu en retirant la section nulle \`a l'espace total du fibr\'e.
On identifie $\A_k^{d+1}\setminus\{0\}\to\P_k^d$ au $\G_m$-torseur $\widetilde{\ecO_{\P_k^{d}}(1)}$.
On note $X_P$ la vari\'et\'e
$\prod_{\alpha\in \Sigma(1)}\left(\A_k^{\,\da+1}\setminus \{0\}\right)/T_P$
et $\pi_M$ le $T_M$-torseur
$
X_P \to \prod_{\alpha\in \Sigma(1)} \PP^{\da}_k.
$

On d\'efinit 
sur $X_P$
une $(\Sigma,\bd)$-collection du second type $M$-triviale universelle.
On pose 
$\ecLa=\pi_M^{\ast} \ecO_{\P_k^{\da}}(1)$. Pour $i=0,\dots,\da$, 
on pose $s_{\alpha,i}=\pi_M^{\ast}v_{\alpha,i}$, o\`u $v_{\alpha,i}$ est la base canonique
de $H^0(\PP^{\da}_k,\ecO_{\P_k^{\da}}(1))$. 
Il reste \`a d\'efinir les trivialisations $c_m$. 
Soit $m\,:\,T_M\to \G_m$ un \'el\'ement de $M$. 
Le produit contract\'e 
$
\ecT_m=
X_P
\times^{T_M,m} \G_m
$
est un $\G_m$-torseur au-dessus de $\prod \P_k^{\da}$,
canoniquement isomorphe au produit contract\'e
$$
\left(\produ{\alpha\in \Sigma(1)}\A_k^{\,\da+1}\setminus \{0\}\right)
\times^{\G_m^{\Sigma(1)},m\circ \pi} \G_m.
$$
Ainsi $\ecT_m$ s'identifie canoniquement au 
$\G_m$-torseur $\widetilde{\otimes \ecO_{\P_k^{\da}}(m_{\alpha})}$. 
Donc le tir\'e en arri\`ere de $\ecT_m$ sur $X_P$ 
s'identifie canoniquement \`a $\widetilde{\otimes \ecLa^{m_{\alpha}}}$.
Mais par ailleurs ce tiré en arrière
 est canoniquement isomorphe au $\G_m$-torseur trivial, 
d'o\`u la trivialisation $c_m$. 
Par construction on a $c_m\otimes c_{m'}=c_{m+m'}$.

Soit $S$ un $k$-sch\'ema. \`A tout morphisme $S\to X_P$ 
on associe la $(\Sigma,\bd)$-collection du second type $M$-trivialisée obtenue en tirant en arri\`ere
la $(\Sigma,\bd)$-collection du second type $M$-trivialisée universelle.
Le fait que ceci d\'efinisse une bijection entre l'ensemble des points de $X_P$ \`a valeurs dans $S$ et l'ensemble
des classes d'isomorphisme 
$(\Sigma,\bd)$-collections du second type $M$-trivialisées sur $S$ se d\'emontre alors
de la même fa\c con que Cox démontre le th\'eor\`eme principal de \cite{Cox:funct}. Alternativement, on 
peut exploiter la structure de vari\'et\'e torique sur $X_P$ pour appliquer directement le r\'esultat de Cox.

Comme la bijection décrite ci-dessus induit une bijection entre les points de $\widetilde{\left(\PP^1\right)^{B_{\Sigma}}_{\bd}}/T_P$
\`a valeurs dans $S$ et l'ensemble des classes d'isomorphisme de $(\Sigma,\bd)$-collections du
second type $M$-trivialisées non d\'eg\'en\'er\'ees sur $S$, 
on obtient le résultat pour $\widetilde{\left(\PP^1\right)^{B_{\Sigma}}_{\bd}}/T_P$.
\end{proof}

\begin{prop}\label{prop:repda}
Soit $\bd\in \N^{\Sigma(1)}_{\bigast}$.
La vari\'et\'e $\widetilde{\left(\PP^1\right)^{B_{\Sigma}}_{\bd}}/T_{\Pic(\xs)}$ 
repr\'esente 
le foncteur qui \`a un $k$-sch\'ema $S$ associe l'ensemble des
classes d'isomorphisme de $(\Sigma,\bd)$-collections du second type $\cara{U}$-trivialisées non d\'eg\'en\'er\'ees.
Cette vari\'et\'e repr\'esente \'egalement le foncteur
$\HOM^{U,\bd}_{k}(\PP^1,\xs)$.
\end{prop}
\begin{proof}[D\'emonstration]
Compte tenu de la suite exacte 
\eqref{eq:suite_exacte_pic_tor}, la premi\`ere assertion de la proposition d\'ecoule 
du lemme \ref{lm:rep_quot} appliqué à $M=\cara{U}$.

Pour montrer la deuxi\`eme assertion, il suffit de construire une bijection fonctorielle en $S$
entre les classes d'isomorphisme de $(\Sigma,\bd)$-collections du second type 
$\cara{U}$-trivialisées non-d\'eg\'en\'er\'ees sur $S$
et les classes d'isomorphismes de $(\PP^1,\Sigma,\bd)$-collections non d\'eg\'en\'er\'ees sur $S$.

Soit $\left((\ecLa,(s_{\alpha,i})),(c_m)\right)$ 
une $(\Sigma,\bd)$-collection du second type sur $S$ supposée
$\cara{U}$-trivialisée et non-d\'eg\'en\'er\'ee. 
On lui associe la $(\PP^1,\Sigma,\bd)$-collections non d\'eg\'en\'er\'ee 
$((\ecMa,u_{\alpha}),(c'_m))$ suivante :
on pose 
$\ecMa=p_{1,S}^{\ast}\ecO_{\PP^1_k}\bd\otimes p_{2,S}^{\ast}\ecLa$.
On peut identifier $H^0(\PP^1_S,\ecMa)$
\`a 
$$
H^0(\PP^1,\ecO_{\PP^1}(\da)) \otimes H^0(S,\ecLa)\isom   H^0(S,\ecLa)^{\da+1}.
$$ 
On pose alors
$
u_{\alpha}\eqdef(s_{\alpha,0},\dots,s_{\alpha,\da}).
$
Soit $m\in \cara{U}$. Comme $\bd\in \N^{\Sigma(1)}_{\bigast}$,
on a un isomorphisme
$$
c''_m\,:\,\otimes \ecO_{\PP^1_k}(\da)^{\acc{m}{\roa}} \isom \ecO_{\PP^1_k}.
$$
On pose $c'_m=p_{1,S}^{\ast}c''_m\otimes p_{2,S}^{\ast}c_m$.

R\'eciproquement, soit $((\ecMa,u_{\alpha}),c'_m)$ une classe
d'isomorphisme de $(\PP^1,\Sigma,\bd)$-collections non d\'eg\'en\'er\'ees
sur $S$.
On peut supposer qu'on a 
$$
\ecMa=p_{1,S}^{\ast}\ecO_{\PP^1_k}(\da)\otimes p_{2,S}^{\ast}\ecLa.
$$ 
On associe 
à la classe d'isomorphisme ci-dessus
la $(\Sigma,\bd)$ 
collection $\cara{U}$-trivialisée 
$\left((\ecLa,u_{\alpha}),(p_{2,S}^{\ast}c'_m)\right)$. 
On v\'erifie ais\'ement que ceci fournit la bijection cherch\'ee.
\end{proof}

\begin{nota}
Si $d\in \N$, 
on note $n_{\Sigma}(d)$ le cardinal de l'ensemble
$$
\left\{
\bd\in \N^{\Sigma(1)}_{\bigast}, 
\quad 
\sumu{\alpha\in \Sigma(1)} \da=d
\right\}.
$$

Pour tout $d\geq 0$, on note $W_{d}$ le $k$-sch\'ema
$$
\coprod_{
\substack{
\bd\in \N^{\Sigma(1)}_{\bigast}, 
\\ \\ \sumu{\alpha\in \Sigma(1)} \da=d
}
} 
\widetilde{\left(\PP^1\right)^{B_{\Sigma}}_{\bd}}/T_{\Pic(\xs)}.
$$
\end{nota}

$W_{d}$ est donc une union 
disjointe de $n_{\Sigma}(d)$ vari\'et\'es d\'efinies sur $k$,
chacune de ces vari\'et\'es \'etant g\'eom\'etriquement 
irr\'eductible de dimension 
$d-\rg(\Pic(\xs)$.

Rappelons que $\HOM^{U,\ecL_0,d}_{k}(\PP^1,\xs)$ d\'esigne le foncteur qui 
\`a un $k$-sch\'ema $S$ associe
$$
\left\{
\varphi\in \Hom_k(\PP^1_S,X_{\Sigma}),
\quad 
\forall s\in S,\quad 
\deg\left(\varphi_s^{\ast}(\ecL_{0})\right)=d
\wedge
\varphi_s\in U(\kappa_s(t))
\right\},
$$ 
et que, d'apr\`es le lemme \ref{lm:repr},
ce foncteur est repr\'esentable par un sch\'ema quasi-projectif not\'e $U_{0,d}$.

\begin{lemme}\label{lm:param_mot}

Pour tout $d\geq 0$, $W_{d}$ 
est isomorphe \`a $U_{\ecL_0,d}$.

On a l'\'egalit\'e 
$$
\left[U_{0,d}\right]
=
(\LL-1)^{\,\dim(\xs)}
\,
\sum_{
\substack{
\bd\in \N^{\Sigma(1)}_{\bigast}, 
\\ \\ \sumu{\alpha\in \Sigma(1)} \da=d
}
} 
\left[\left(\PP^1\right)^{B_{\Sigma}}_{\bd}\right].$$
\end{lemme}
\begin{proof}[D\'emonstration]
La premi\`ere assertion d\'ecoule de la proposition \ref{prop:repda}
et du fait que $\sum_{\alpha}\mathfrak{D}_{\alpha}$ est un diviseur anticanonique sur $\xs$.

L'application 
$$
\widetilde{\left(\PP^1\right)^{B_{\Sigma}}_{\bd}}
/
T_{\Pic(\xs)}
\to 
\widetilde{\left(\PP^1\right)^{B_{\Sigma}}_{\bd}}
/
\G_m^{\Sigma(1)}
=\left(\PP^1\right)^{B_{\Sigma}}_{\bd}
$$
est un torseur sous $\G_m^{\Sigma(1)}/T_{\Pic(\xs)}=U$, localement trivial pour la topologie de Zariski car $U$ est d\'eploy\'e. 
Ainsi on a
$$
\left[
\widetilde{\left(\PP^1\right)^{B_{\Sigma}}_{\bd}}/T_{\Pic(\xs)}
\right]
=
\left[
\left(\PP^1\right)^{B_{\Sigma}}_{\bd}
\right]
\,
[U].
$$
On en déduit le r\'esultat. 
\end{proof}

\begin{cor}\label{cor:rep_pos_ques_1}
Soit $V$ une variété torique projective lisse et déployée définie sur
un corps $k$, d'orbite ouverte $U$. Pour tout entier $d\geq 1$, soit $U_{0,d}$ la variété
quasi-projective paramétrant les $k$-morphisme $\PP^1_k\to V$ dont
l'image rencontre $U$ et de degré anticanonique $d$. Soit
$\rho(U_{0,d})$ le nombre de composantes géométriquement irréductibles
de dimension maximale de $U_{0,d}$.
On a alors
$$
\overline{\lim_{d\to \infty}} \,\, \,\,\frac{\dim(U_{0,d})}{d}=1 
$$
et
$$
\overline{\lim_{d\to \infty}} \,\, \,\,\frac{\log(\rho(U_{0,d}))}{\log(d)}=\rg(\Pic(V)).
$$
En d'autres termes, dans le cas d'une courbe rationnelle, la r\'eponse \`a la question 
\ref{qu:conj_manin_mot_irr} est positive pour les variétés toriques
projectives, lisses et déployées en prenant pour ouvert
l'orbite ouverte.
\end{cor}
\begin{proof}
Ceci découle du lemme \ref{lm:param_mot} 
et de la th\'eorie
du polyn\^ome d'Ehrahrt qui permet 
de montrer qu'on a
$$
\lim_{d\to \infty} \,\, \,\,\frac{\log(n_{\Sigma}(d))}{\log(d)}=\card{\Sigma(1)}-r=\rg(\Pic(\xs)).
$$
\end{proof}
On d\'eduit \'egalement imm\'ediatement du lemme \ref{lm:param_mot} qu'on a l'expression 
\begin{equation}\label{eq:expr_zetamot}
\zetamthz(T)
=(\LL-1)^{\,\dim(\xs)}\,
\sum_{
\bd\in \N^{\Sigma(1)}_{\bigast}
}
\,
\left[\left(\PP^1\right)^{B_{\Sigma}}_{\bd}\right]\,
T^{\,\sumu{\alpha}\da}.
\end{equation}
La description de $\xs$ comme quotient g\'eom\'etrique permet par ailleurs de montrer
la formule suivante, qui est une version motivique d'une formule donnant le nombre
de points d'une vari\'et\'e torique sur un corps fini.
On renvoie \`a la sous-section \ref{subsec:mob_ev} pour la d\'efinition de $\bs$ et au d\'ebut de 
la section \ref{sec:fonc_mob_mot} pour la d\'efinition de $\mu^0_{\bs}$.
\begin{prop}\label{prop:form_phinxs}
Soit $k$ un corps de caractéristique zéro, $\Sigma$ un éventail
projectif et lisse et $\xs$ la $k$-variété torique associée.
On a l'\'egalit\'e
$$
\sum_{(\na)\in \{0,1\}^{\Sigma(1)}}
\mu^0_{\bs}((\na))\,\LL^{\,-n\,\sum \na}
=
\left(1-\LL^{-n}\right)^{\rg(\Pic(\xs))}
\,
\frac{\Phi_n^{\chi}(\xs)}{\LL^{\,-n\,\dim(\xs)}}.
$$
\end{prop}
\begin{proof}[D\'emonstration]
On sait que $\ts$ est un torseur sous le tore d\'eploy\'e 
$T_{\Pic(\xs)}\isom\G_m^{\rg(\Pic(\xs))}$
au-dessus de $\xs$. 
D'apr\`es la proposition \ref{prop:phixy}, on en
d\'eduit la relation
$$
\Phi_n^{\chi}(\ts)=\Phi_n^{\chi}(\xs)\,\Phi_n^{\chi}(\G_m)^{\rg(\Pic(\xs))}.
$$
D'apr\`es le corollaire \ref{cor:phichiad}, on a $\Phi_n^{\chi}(\G_m)=\LL^n-1$.
Calculons $\Phi_n^{\chi}(\ts)$. 
Pour $(\na)\in \{0,1\}^{\Sigma(1)}$
notons $\A_{(\na)}\eqdef\capu{\na=1}\{X_{\alpha}=0\}$.
On a ainsi
$$
\Phi_n^{\chi}\left(\A_{(\na)}\right)
=
\Phi_n^{\chi}\left(\A^{\sum \na}\right)
=
\LL^{\,n\,\sum \na},
$$
la derni\`ere \'egalit\'e provenant du corollaire \ref{cor:phichiad}.
On a alors
$$
\ts
=
\A^{\Sigma(1)}\setminus \bigcup_{(\na) \in \bs}\A_{(\na)}
=
\A^{\Sigma(1)}\setminus \bigcup_{(\na) \in \bs^{\min}}\A_{(\na)}
.
$$
D'apr\`es la proposition \ref{prop:phiuf}, 
on a
$$
\Phi_n^{\chi}
\left(\bigcup_{(\na) \in \bs}\A_{(\na)}\right)
=
\sum_{\vide\,\neq \,J \,\subset \,\bs^{\min}}
(-1)^{\,1+\card{J}}
\,
\Phi_n^{\chi}
\left(
\bigcap_{\alpha\in J}\A_{(\na)}
\right).
$$
Compte tenu du fait que pour $(\na)\in \bs$ on a 
$$
\A_{(\na)}
=
\bigcap_{(\na')\in \bs^{\min},\quad (\na')\leq (\na)}
\A_{(\na')},
$$
et rappelant que
$$
\ell_{\bs}((\na))
=
\card{\{(\na')\in \bs^{\min},\quad (\na')\leq (\na)\}},
$$
ceci se r\'e\'ecrit
\begin{align*}
\Phi_n^{\chi}
\left(\bigcup_{(\na) \in \bs}\A_{(\na)}\right)
&=
\sum_{(\na)\in \bs}
(-1)^{\,1+\ell_{\bs}((\na))}
\Phi_n^{\chi}
\left(
\A_{(\na)}
\right)
\\
&=
-\sum_{(\na)\in \{0,1\}^{\Sigma(1)}}
\mu^0_{\bs}((\na))
\,
\LL^{\,n\,\sum \na}.
\end{align*}
On a donc bien la formule annonc\'ee.
\end{proof}
\subsection{Le cas des surfaces de Hirzebruch}\label{subsec:hirz}

Nous traitons ce cas particulier s\'epar\'ement, 
car d'une part on peut travailler ici dans l'anneau $\mk$,
d'autre part on obtient le fait remarquable que la s\'erie $\zetamh(T)$ 
est une fonction rationnelle en $T$.  
Plus précisément, nous démontrons le point
\ref{item:1:theoprinc} du théorème \ref{thm:theoprinc},
dont nous rappelons l'énoncé.
\begin{thm}\label{thm:mot_hirz}
Soit $k$ un corps.
Soit $m\geq 0$ un entier.
Soit $\Sigma$ l'\'eventail de $\Z^2\otimes \R$ dont les rayons 
sont engendr\'es par $\rho_1=(1,0)$, $\rho_2=(-1,m)$, 
$\rho_3=(0,1)$, $\rho_4=-\rho_3$. 
La $k$-vari\'et\'e torique d\'eploy\'ee $\xs$
associ\'ee est la $m$-\`eme surface de Hirzebruch $\hirz{m}$.
On note $U$ son orbite ouverte.

Alors l'élément de $\mk[[T]]$  
\begin{equation}
(1+\LL\,T)\,(1+\LL\,T+\LL^2\,T^2+\dots+\LL^{\,m+1}\,T^{\,m+1})
\,(1-\LL\,T)^{2}\,Z^{\text{mot}}_{\PP^1,U,h_0}(T)
\end{equation}
est un polynôme dont la valeur en $\LL^{-1}$ est 
$
\LL^{\,2}\,(1-\LL^{-2})^2
$.
\end{thm}
\begin{rem}
On a $\rg(\Pic(\hirz{m}))=2$ et $\alpha^{\ast}(\hirz{m})=\frac{1}{m+1}$.
Par ailleurs on a 
\begin{align*}
\LL^{\,2}\,(1-\LL^{-2})^2
&=
\LL^{\,2}\,
\left(
\frac{1}{1-\LL^{-1}}
\right)^2
\,
\left[
(1-\LL^{-2})(1-\LL^{-1})
\right]^2
\\
&=
\LL^{\,2}\,
\left(
\frac{1}{1-\LL^{-1}}
\right)^2\,
\left(
\frac{1}{\zetapm(\LL^{-2})}
\right)^2
\\
&=
\LL^{\,\dim(\xs)}
\,
\left(
\left[(1-\LL\,T)\,\zetapm(T)\right](\LL^{-1})
\right)^{\,\rg(\Pic(\xs))}
\,\times\,
\sum_{\bd}
\mubsmp(\bd)\,\LL^{\,-\sumu{\alpha}\da},
\end{align*}
la derni\`ere \'egalit\'e provenant de l'expression de $\mubsmp$ obtenue au paragraphe \ref{par:musm_hirz}.
Nous obtenons donc bien une version motivique du 
r\'esultat principal de \cite{Bou:hirz} dans le cas où la courbe
est suppos\'ee rationnelle. Ce dernier résultat se déduit d'ailleurs du
théorème \ref{thm:mot_hirz} en sp\'ecialisant via le morphisme $\#_k$,
grâce à la rationnalité de la fonction z\^eta des hauteurs motivique.
\end{rem}
\begin{proof}[D\'emonstration]
On a d'apr\`es le lemme \ref{lm:param_mot}, 
la formule \eqref{eq:rel_mub_mot_2} 
et la description de $\Sigma$ 
\begin{align*}
\zetamthz(T)
&=
(\LL-1)^2
\sum_{
\substack{
(d_i)\in \N^{4}
\\ \\ 
d_1=d_2 
\\ \\ d_3+m\,d_2=d_4
}
}
\left[
U_{(d_i)}
\right]\,
T^{\,d_1+d_2+d_3+d_4}\\
&=
(\LL-1)^2
\negthickspace
\sum_{
\substack{
(d_i)\in \N^{4},
\,(e_i)\in \N^{\,4} 
\\ \\ 
d_1+e_1=d_2+e_2 
\\ \\ 
d_3+m\,d_2+e_3+m\,e_2=d_4+e_4
}
}
\negthickspace
\mubsmp(e_1,e_2,e_3,e_4)
\left[\PP^{\,d_1}\right]\,
\left[\PP^{\,d_2}\right]\,
\left[\PP^{\,d_3}\right]\,
\left[\PP^{\,d_4}\right]
\place{1cm}
{\times\,T^{\,d_1+d_2+d_3+d_4+e_1+e_2+e_3+e_4}}.
\end{align*}
En utilisant l'expression de $\mubsm$ obtenue \`a la section \ref{par:musm_hirz},
la relation \eqref{eq:muxi} et la relation 
$$(\LL-1)\,\left[\PP^d\right]=\LL^{\,d+1}-1,$$
on obtient la relation (\cf  \cite[section 4.3.2]{Bou:PhD} pour les d\'etails du calcul)
$$
\zetamthz(T)
=(\LL-1)\,\left(
\LL\,\frac{(1-T^2)\,(1+\LL^{\,m+1}\,T^{\,m+2})}{(1-\LL^2\,T^2)\,(1-\LL^{\,2+m}\,T^{\,2+m})}
-
\frac{1+\LL\,T^{\,2+m}}{1-\LL^{\,2+m}\,T^{\,2+m}}
\right)
$$
d'o\`u le r\'esultat annonc\'e.
\end{proof}

\subsection{Le cas g\'en\'eral}\label{subsec:cas_gen}

Soit $k$ un corps de caract\'eristique z\'ero et $\Sigma$ un éventail
projectif et lisse.
On conserve les notations des sous-sections \ref{subsec:mot_vtetor}
et \ref{subsec:param_mot}. 
Le but de cette partie est de 
 démontrer le point 
\ref{item:2:theoprinc} du théorème \ref{thm:theoprinc}. Pour cela, on va reprendre au niveau de l'anneau de
motifs virtuels $\mkc$ la strat\'egie employ\'ee dans \cite{Bou:vtetor}.
Soulignons que tous les calculs qui suivent sont 
valables sur l'anneau $\mk$ pour un corps $k$ quelconque
(en fait si $k$ est fini leur sp\'ecialisation via $\#_k$ \og
redonne\fg\  les calculs effectu\'es dans \cite{Bou:vtetor}), 
mais pas a priori les r\'esultats de convergence 
pour lesquels on aurait besoin d'une r\'eponse positive \`a la question
\ref{ques:muxb}. 

\begin{nota}
Pour all\'eger l'écriture, pour toute $k$-vari\'et\'e $V$, le motif
virtuel associ\'ee \`a $V$
sera noté $\symb{V}$ en lieu et
place de $\chi([V])$. 
Par ailleurs, on notera $\musc$ la fonction $\mubsp$.
\end{nota}

D'apr\`es les formules \eqref{eq:expr_zetamot} et \eqref{eq:rel_mub_mot},
on a  
\begin{align*}
\zetamthz(T)
&=(\LL-1)^{\,\dim(\xs)}
\negthickspace\negthickspace
\sum_{
\bd\in \N^{\Sigma(1)}_{\bigast}
}
\left[\left(\PP^1\right)^{B_{\Sigma}}_{\bd}\right]\,
T^{\,\sumu{\alpha}\da}
\\
&=
(\LL-1)^{\,\dim(\xs)}
\negthickspace\negthickspace
\sum_{
\substack{
\bd,\be\in \N^{\Sigma(1)}
\\ \\
(\bd+\be)\in \N^{\Sigma(1)}_{\bigast}
}}
\negthickspace\negthickspace
\musc(\be)
\,
\prod_{\alpha\in \Sigma(1)} 
\,
\left[\PP^{\da}\right]\,
T^{\,\sumu{\alpha}(\da+\ea)}
\\
&=
(\LL-1)^{\,\dim(\xs)}
\,
\sum_{
\be\in \N^{\Sigma(1)}
}
\musc(\be)
\,
\sum_{
\substack{
\bd\in \N^{\Sigma(1)}_{\bigast}
\\ \\ 
\bd\geq \be
}
}
\,
\prod_{\alpha\in \Sigma(1)} 
\,
\left[\PP^{\da-\ea}\right]\,
T^{\,\sumu{\alpha}\da}.
\end{align*}
Rappelons que l'injection naturelle
$
\N^{\Sigma(1)}
\longto 
\left(\oplus \Z\,{\mathfrak{D}_{\alpha}}
\right)^{\vee}
$
et le dual de la suite exacte \eqref{eq:suite_exacte_pic_tor} 
permettent d'identifier $\N^{\Sigma(1)}_{\bigast}$ 
au sous-ensemble de $\Pic(\xs)^{\vee}$ constitué des éléments
$y$ v\'erifiant $\acc{y}{\cDa}\geq 0$
pour tout $\alpha\in \Sigma(1)$.
Ainsi, pour $\be\in \N^{\Sigma(1)}$, on a 
$$
\sum_{
\substack{
\bd\in \N^{\Sigma(1)}_{\bigast}
\\ \\ 
\bd\geq \be
}
}
\,\,\,
\prod_{\alpha\in \Sigma(1)} 
\,
\left[\PP^{\da-\ea}\right]\,
T^{\,\sumu{\alpha}\da}
=
\sum_{
\substack{ 
y\in \Pic(\xs)^{\vee}
\\ \\
\acc{y}{\cDa}\geq \deg(\Ea)
}
}
\prod_{\alpha\in \Sigma(1)}
\left(
\frac{\LL^{1+\acc{y}{\cDa}-\ea}-1}{\LL-1}
\right)
\, 
T^{\,\,\acc{y}{\underset{\alpha}{\sum} \cDa}}.
$$
On a donc
\begin{align*}
\frac{\zetacthz(T)}{(\LL-1)^{\rg(\cara{U})-\card{\Sigma(1)}}}\hskip-0.2\textwidth&\\
&=
\sum_{ \be \in \N^{\Sigma(1)}}
\,
\musc( \be )
\,
\Big(
\sum_{
\substack{ 
y\in \Pic(\xs)^{\vee}
\\ \\
\acc{y}{\cDa}\geq \ea
}
}
\,
\prod_{\alpha\in \Sigma(1)}
\left(
\LL^{1+\acc{y}{\cDa}-\ea}-1
\right)
\, 
T^{\,\,\acc{y}{\underset{\alpha}{\sum} \cDa}}
\Big)
\\
&=
\sum_{ \be \in \N^{\Sigma(1)}}
\,
\musc( \be )
\,
\Big(
\sum_{
\substack{ 
y\in \ceff(\xs)^{\vee}\cap \Pic(\xs)^{\vee}
\\ \\
\acc{y}{\cDa}\geq \ea
}
}
\,
\prod_{\alpha\in \Sigma(1)}
\left(
\LL^{1+\acc{y}{\cDa}-\ea}-1
\right)
\, 
T^{\,\,\acc{y}{\ecL_0}}
\Big).
\end{align*}
On d\'ecompose \`a pr\'esent $\zetacthz(T)$ en une somme de plusieurs termes,
dont on estimera ensuite le comportement  s\'epar\'ement. 
On \'ecrit
$$
\zetacthz(T)
=
(\LL-1)^{\,-\rg(\Pic(\xs))}
\sum_{A\subset \Sigma(1)}(-1)^{\,\card{A}}
\,
Z_A(T)
$$
avec, pour $A\subset\Sigma(1)$,
$$
Z_A(T)=\sum_{ \be \in \N^{\Sigma(1)}} \musc( \be )\,Z_{A, \be }(T),
$$
$Z_{A, \be }(T)$ d\'esignant la s\'erie
$$
\sum_{
\substack{
y\in\ceff(\xs)^{\vee}\cap \Pic(\xs)^{\vee}
\\ \\
\forall\,\alpha\in\Sigma(1),
\quad
\acc{y}{\cDa} \,\geq \,\ea
}
}
\,
\LL^{\,\sumu{\alpha\notin A}1+\acc{y}{\cDa} -\ea}
\,
T^{\,\acc{y}{\ecL_0}}
=
\sum_{
\substack{
y\in\ceff(\xs)^{\vee}\cap \Pic(\xs)^{\vee}
\\ \\
\forall\,\alpha\in\Sigma(1),
\quad
\acc{y}{\cDa} \,\geq \,\ea
}
}
\,
\LL^{\,\sumu{\alpha\notin A}1+\acc{y}{\cDa} -\ea}
\,
T^{\,\acc{y}{\ecL_0}}.
$$

\begin{nota}\label{nota:dec_ceff}
On \'ecrit $\ceff(\xs)^{\vee}$ comme le support d'un 
\'eventail r\'egulier $\Delta$.
Pour $i\in\Delta(1)$ on note $m_i$ le g\'en\'erateur du rayon $i$. 
Pour toute partie $I$ de $\Delta(1)$ on note 
$
\cone(I)\eqdef\sum_{i\in I}\,\N_{>0}\,m_i
$
(avec la convention $\cone(\vide)=\{0\}$),
de sorte que $\cone(\delta(1))$ est l'ensemble des points du r\'eseau $\Pic(\xs)^{\vee}$ 
contenu dans l'int\'erieur relatif du c\^one $\delta$. 
\end{nota}
\begin{rem}
On a alors
$$
Z_{\ceff(V),\ecL_0}(T)
=
\sum_{\delta\in \Delta}
\,
\prod_{i\in \delta(1)}
\left(
\frac{1}{1-T^{\,\acc{m_i}{\ecL_0}}}-1
\right)
$$
d'o\`u
\begin{equation*}
\alpha^{\ast}(\xs)
=
\sum_{\substack{\delta\in \Delta \\ \dim(\delta)=\rg(\Pic(V))}}\,
\prod_{i\in \delta(1)}\frac{1}{\acc{m_i}{\ecL_0}}.
\end{equation*}
\end{rem}

Soit $A\subset \Sigma(1)$. 
On écrit
$$
Z_A(T)=\sum_{\delta\in\Delta}\,Z_{A,\delta}(T)
$$
avec
$$Z_{A,\delta}(T)=\sum_{ \be \in \N^{\Sigma(1)}} \musc( \be )\,Z_{A,\delta, \be }(T),$$
$Z_{A,\delta, \be }(T)$ d\'esignant la s\'erie
$$
\sum_{
\substack{
y\in \cone(\delta(1))
\\ \\
\forall\alpha\in\Sigma(1),
\quad
\acc{y}{\cDa} \geq \ea
}
}
\,
\LL^{\,\sumu{\alpha\notin A}1+\acc{y}{\cDa} -\ea}
\,
T^{\,\acc{y}{\ecL_0}}.
$$
De la m\^eme fa\c con que le th\'eor\`eme 1 (p.178) de \cite{Bou:vtetor} se d\'eduisait 
des propositions 1-(3) (p.181), 3 (p.190) et 4 (p.195) de (op.cit.),  
le point \ref{item:2:theoprinc} du th\'eor\`eme  \ref {thm:theoprinc} se d\'eduit alors de la proposition \ref{prop:form_phinxs} 
et des propositions \ref{prop:Avide_mot} et \ref{prop:Anonvide_mot} \'enonc\'ees
et d\'emontr\'ees ci-dessous.
Les démonstrations, tr\`es similaires 
\`a celles des propositions analogues de \cite{Bou:vtetor}, 
sont en outre simplifi\'ees par le fait qu'on \'etudie des convergences pour 
une norme non-archim\'edienne. 

\subsubsection{Le cas $A=\vide$.}
\begin{nota}
Soient $\be\in\N^{\,\Sigma(1)}$,
$\delta$  un cône de $\Delta$ et 
 $J$ une partie de $\Sigma(1)$.
On pose 
$$
\cone(\delta(1))_{J,\, \be }
\eqdef
\left\{
y\in \cone(\delta(1)),
\quad
\forall\alpha\in J,
\quad
\acc{y}{\cDa} \leq \ea\right
\}.
$$
\end{nota}
Soit $\delta$ un cône de $\Delta$. 
Nous écrivons
$$
Z_{\vide,\delta}(T)
=
\sum_{J\subset\Sigma(1)}
\,
(-1)^{\card{J}}
\,
Z_{\vide,\delta,J}(T)
$$
avec
\begin{equation}\label{eq:zvidedeltaj}
Z_{\vide,\delta,J}(T)
=
\sum_{ \be \in \N^{\Sigma(1)}} 
\musc( \be )
\,
Z_{\vide,\delta,J, \be }(T),
\end{equation}
$Z_{\vide,\delta,J, \be }(T)$ d\'esignant la s\'erie
$$
\sum_{
\substack{
y\in \cone(\delta(1))
\\ \\
\forall\,\alpha\in J,\,\acc{y}{\cDa} \,\,<\,\ea
}
}
\,
\LL^{\,\sumu{\alpha\in \Sigma(1)}\,1+\acc{y}{\cDa} -\ea\,}
\,
T^{\,\acc{y}{\ecL_0}}
=
\sum_{
\cone(\delta(1))_{J,\, \be }
}
\,
\LL^{\,\sumu{\alpha\in \Sigma(1)}\,1+\acc{y}{\cDa} -\ea\,}
\,
T^{\,\acc{y}{\ecL_0}}.
$$
\begin{prop}\label{prop:Avide_mot}
Soit $\delta$ un cône de $\Delta$.
La s\'erie
$$
(1-\LL\,T)^{\dim(\delta)}\,Z_{\vide,\delta,J}(T)
$$
converge dans $\mhatcq$ en $T=\LL^{-1}$. 
Si $J\neq\vide$ et $\dim(\delta)=\rg(\Pic(\xs))$,
sa valeur en $\LL^{-1}$ est nulle. 

La s\'erie 
$$
(1-\LL\,T)^{\,\rg(\Pic(\xs))}
\,
\sum_{
\substack
{
\delta\in\Delta
\\ \\ 
\dim(\delta)
=
\rg(\Pic(\xs))
}
}
\,
Z_{\vide,\delta,\vide}(T)
$$
converge dans $\mhatcq$ en $T=\LL^{-1}$ 
vers
$$
\alpha^{\ast}(\xs)
\,\
\LL^{\,\card{\Sigma(1)}}
\,
\sum_{ \be \in \N^{\Sigma(1)}}
\musc( \be )
\LL^{\,-\sum \ea}.
$$ 
\end{prop}
Cette proposition d\'ecoule de \eqref{eq:zvidedeltaj}, 
du corollaire \ref{cor:muxb} 
et du lemme \ref{lm:gen1_mot} ci-dessous.
\begin{lemme}\label{lm:gen1_mot}
Soient $ \be $ un \'el\'ement de $\N^{\,\Sigma(1)}$,
$\delta$ un c\^one de $\Delta$ et  $J$ une partie de
$\Sigma(1)$. 
Il existe alors un sous-ensemble $I'$ de $\delta(1)$ et 
un polyn\^ome $P_{I'}$ à coefficients dans $\Z[\LL]$ tels que~:
\begin{enumerate}
\item
pour tout entier $\kappa\geq 1$ on a
$$
P_{I'}
\left(\LL^{-\kappa}\right)
\in 
\fil^{\,-\card{\Sigma(1)}+\sumu{\alpha\in\Sigma(1)} \ea}\mlocc\quad;
$$
\item
on a la relation
\begin{equation}
\label{eq:sumycone}
\sum_{y\in \cone(\delta(1))_{J,\, \be }}
\LL^{\,\,\sumu{\alpha\in\Sigma(1)}1+\acc{y}{\cDa} -\ea}
\,
T^{\,\acc{y}{\ecL_0}}
=
\left(
\,
\prod_{i\in I'}
\,
\left[
\left(
\,1-(\LL\,T)^{\,-\acc{m_i}{\omega}}
\right)^{-1}-1
\right]
\right)
\,
\times
\,
P_{I'}(T).
\end{equation}
\end{enumerate}
Supposons en outre $\delta$ de dimension maximale. 
Alors si $J$ est non vide, 
$I'$ est un sous-ensemble strict de $\delta(1)$. 
Si $J$ est vide, 
on a $I'=\delta(1)$ et $P_{\delta(1)}$ est constant
\'egal \`a
$
\LL^{\,\card{\Sigma(1)}-\sum_{\alpha\in\Sigma(1)}\ea}.
$ 
\end{lemme}
\begin{proof}[D\'emonstration]  
On a $\cone(\delta(1))_{\vide,\, \be }=\cone(\delta(1))$.
Donc si $J$ est vide le membre de gauche \eqref{eq:sumycone} 
s'\'evalue imm\'ediatement et vaut 
$$
\LL^{
\,
\card{\Sigma(1)}
-
\sumu{\alpha\in\Sigma(1)}\ea
} 
\left(
\,
\prod_{i\in \delta(1)}
\,
\left[
\left(\,1-(\LL\,T)^{\,-\acc{m_i}{\omega}}
\right)^{-1}-1
\right]
\right).
$$ 
Si $J$ n'est pas vide, posons
$$
I_{J,1}
=
\{
\,
i\in \delta(1),
\quad
\forall\,\alpha\in J,
\quad
\acc{m_i}{\cDa} 
=
0
\}
$$
et $I_{J,2}=\delta(1)\setminus I_{J,1}$. 
En particulier on a $\cone(I_{J,1})_{J,\, \be }=\cone(I_{J,1})$.
Si on note
$$
P(T)=\sum_{y_2\in \cone(I_{J,2})_{J,\, \be }}
\LL^{\,\,\sumu{\alpha}\,1+\acc{y_2}{\cDa} -\ea\,}
\,\,
T^{\,\acc{y_2}{\ecL_0}},
$$
le membre de gauche de  \eqref{eq:sumycone} est égal au produit
\begin{equation}\label{eq:prod}
P(T)
\,
\times\,
\sum_{y_1\in \cone(I_{J,1})}
\,
\LL^{\,\acc{y_1}{\sumu{\alpha\notin J} \cDa}}
\,
T^{\,\acc{y_1}{\ecL_0}}.
\end{equation}
Si $y_1\in \cone(I_{J,1})$, 
on a 
$
\acc{y_1}{\sumu{\alpha\notin J} \cDa} 
=
\acc{y_1}{\ecL_0}. 
$ 
Donc le deuxi\`eme facteur de \eqref{eq:prod}
est \'egal \`a
$$
\prod_{i\in I_{J,1}}
\,
\left[
\left(
\,1-(\LL\,T)^{\,-\acc{m_i}{\omega}}
\right)^{-1}-1
\right].
$$
Passons au facteur $P(T)$. 
De la même manière que dans la preuve du lemme 3
de \cite{Bou:vtetor},
on voit facilement que $\cone(I_{J,2})_{J,\, \be }$ est fini. 
Ainsi $P$ est un polynôme à coefficients dans $\Z[L]$, et 
pour tout entier $\kappa$ on a 
$$
P\left(\LL^{-\kappa}\right)
=
\LL^{
\,\card{\Sigma(1)}-\sumu{\alpha}\ea
}
\,
\sum_{y_2 \in \cone(I_{J,2})_{J,\, \be }}
\LL^{\,(1-\kappa)\,\acc{y_2}{\ecL_0}}.
$$
Pour tout $y_2\in \cone(I_{J,2})_{J,\, \be }$ on a $\acc{y_2}{\ecL_0} \geq 0$.
Donc si $\kappa\geq 1$ on a 
$$
P\left(\LL^{-\kappa}\right)
\in 
\fil^{\,-\card{\Sigma(1)}+\sum \ea}
\,\,
\mlocc.
$$
Enfin, si $\delta$ est un c\^one de dimension maximale, 
les $(m_i)_{i\in \delta(1)}$ forment une $\Z$-base de $\Pic(\xs)^{\vee}$.
Ainsi, si $J$ n'est pas vide on ne peut avoir $I_{J,2}=\vide$. 
Ceci joint au calcul pour $J=\vide$ montre les deux derni\`eres assertions du lemme.
\end{proof}
\paragraph{Le cas $A\neq\vide$.}
Soit $\delta$ un cône de $\Delta$ et $A$ une partie non vide de $\Sigma(1)$.
On \'ecrit
$$
Z_{A,\delta}(T)
=
\sum_{J\,\subset\,\Sigma(1)\setminus A}
\,
(-1)^{\card{J}}
\,
Z_{A,\delta,J}(T)
$$
avec 
\begin{equation}\label{eq:zadeltaj}
Z_{A,\delta,J}(T)
=
\sum_{ \be \in\N^{\Sigma(1)}} 
\,
\musc( \be )
\,
Z_{A,\delta,J, \be }(T),
\end{equation}
l'expression $Z_{A,\delta,J, \be }(T)$ désignant la s\'erie
$$
\sum_{
\substack{
y\in \cone(\delta(1))
\\ \\
\forall\,\alpha\in A,\quad\acc{y}{\cDa} \,\geq \ea
\\ \\
\forall\,\alpha\in J,\quad\acc{y}{\cDa} \,<\,\ea
}
}
\,
\LL^{\,\sumu{\alpha\notin A}1+\acc{y}{\cDa} -\ea}
\,T^{\,\acc{y}{\ecL_0}}.
$$
\begin{prop}\label{prop:Anonvide_mot}
Soit $\delta$ un cône de $\Delta$, $A$ une partie non vide de
$\Sigma(1)$ et $J$ une partie de $\Sigma(1)\setminus A$.
La s\'erie
$$
(1-\LL\,T)^{\dim(\delta)}\,Z_{A,\delta,J}(T)
$$
converge dans $\mhatc$
en $T=\LL^{-1}$.
En outre, si $\dim(\delta)=\rg(\Pic(\xs))$, sa valeur en $\LL^{-1}$ est nulle.
\end{prop}
Cette proposition d\'ecoule de \eqref{eq:zadeltaj}, 
du corollaire \ref{cor:muxb} 
et du lemme \ref{lm:gen2_mot} ci-dessous.
\begin{lemme}\label{lm:gen2_mot}
Soit $\delta$ un cône de $\Delta$, $A$ une partie non vide de
$\Sigma(1)$ et $J$ une partie de $\Sigma(1)\setminus A$.
Soit $ \be \in \N^{\Sigma(1)}$.
Il existe alors  un sous-ensemble $I'$ de $\delta(1)$, 
et un \'el\'ement $R_{I'}$ de  $\Z[\LL][[T]]$    
tels que :
\begin{enumerate}
\item
pour tout entier $\kappa\geq 1$, 
$
R_{I'}\left(\LL^{-\kappa}\right)
$ 
converge dans $\mhatc$ vers un \'el\'ement de 
$\fil^{\,-\card{\Sigma(1)}+\sumu{\alpha\in\Sigma(1)} \ea}\,\,\mhatc$ ;
\item
on a la relation
\begin{equation}\label{eq:sumy:Ri}
\sum_{
\substack{
y\in \cone(\delta(1)) 
\\ \\ 
\forall\,\alpha\in A,
\quad
\acc{y}{\cDa} \,\geq \ea
\\ \\
\forall\alpha\in J,
\quad
\acc{y}{\cDa} \,< \,\ea
}
}
\LL^{
\,\,
\sumu{\alpha\notin A}\,1+\,\acc{y}{\cDa} -\ea\,
}
\,\,
T^{\,\acc{y}{\ecL_0}}
=
\left(
\,
\prod_{i\in I'}
\,
\left[
\left(
\,1-(\LL\,T)^{\,-\acc{m_i}{\omega}}
\right)^{-1}-1
\right]
\right)
\,
\times
\,
R_{I'}(T).
\end{equation}
\end{enumerate}
En outre, si $\delta$ est un c\^one de dimension maximale, 
$I'$ est un sous-ensemble strict de $\delta(1)$. 
\end{lemme}
\begin{proof}[D\'emonstration] On pose
$$
\cone(\delta(1))^{\,A}_{J,\, \be }
=
\left\{
y\in \cone(\delta(1))_{J,\, \be },
\quad
\forall\alpha\in A,
\quad
\acc{y}{\cDa} \geq \ea
\right\}.
$$
Le membre de gauche de \eqref{eq:sumy:Ri}
s'\'ecrit donc
\begin{equation}\label{eq:sumyconeAJ}
\sum_{y\in \cone(\delta(1))^A_{J,\, \be }}
\LL^{
\,\,\sumu{\alpha\notin A}\,(1+\acc{y}{\cDa} -\ea\,)
}
\,\,
T^{\,\acc{y}{\ecL_0}}.
\end{equation}
On pose
$$
I_{A,J,1}
=
\{\,i\in \delta(1),\quad
\forall\,\alpha\in A\cup J,
\quad
\acc{m_i}{\cDa} 
=0\}
$$
et $I_{A,J,2}=\delta(1)\setminus I_{A,J,1}$. 
Ainsi tout \'el\'ement de $\cone(\delta(1))^A_{J,\, \be }$ 
s'\'ecrit de mani\`ere unique $y_1+y_2$ avec 
$
y_1 \in \cone(I_{A,J,1})
$
et
$$
y_2 \in \coprod_{\substack{(\ha)\in \N^{A} \\ \\ \ha\geq \ea }} 
\left\{
y\in \cone(I_{A,J,2})_{J,\, \be },
\quad
\forall \,\alpha\in A,
\quad
\acc{y}{\cDa} \,
=\,\ha
\right\}
.
$$
Pour $(\ha)\in \N^A$, on voit facilement qu'il n'y a qu'un nombre fini 
d'éléments $y_2$ de $\cone(I_{A,J,2})_{J,\, \be }$ vérifiant 
$\acc{y_2}{\cDa}=\ha$ pour tout $\alpha\in A$. 
Pour un tel $y_2$, 
on a pour tout entier $\kappa$
\begin{align*}
\LL^{
\,\,
\sumu{\alpha\notin A}\,(1+\acc{y_2}{\cDa} -\ea\,)
}
\,
\LL^{
\,
-\kappa\,\acc{y_2}{\ecL_0}
}
&
=
\LL^{
\,\,\card{\Sigma(1)}-\card{A}-\sumu{\alpha\notin A}e_a-\sumu{\alpha\in A}\ha
}
\,
\LL^{\,(1-\kappa)\,\acc{y_2}{\ecL_0}
}
\\
&
=
\LL^{
\,\,\card{\Sigma(1)}
-
\card{A}
-
\sumu{\alpha\in \Sigma(1)}\ea-\sumu{\alpha\in A}(\ha-\ea)
}
\,
\LL^{
\,(1-\kappa)\,\acc{y_2}{\ecL_0}
}.
\end{align*}
Notons par ailleurs que $\acc{y_2}{\ecL_0}$ est positif.
Si on pose
$$
R_{(\ha)}(T)
=
\sum_{
\substack{
y_2\in \cone(I_{A,J,2})_{J,\, \be } 
\\ \\ 
\forall\,\alpha\in A,\quad\acc{y_2}{\cDa} \,=\,\ha
}
}
\LL^{\,\,\sumu{\alpha\notin A}\,(1+\acc{y_2}{\cDa} -\ea\,)}
\,
T^{\,\acc{y_2}{\ecL_0}},
$$
ce qui précède montre qu'on a pour tout entier $\kappa\geq 1$
\begin{equation}\label{eq:Rha:in:fil}
R_{(\ha)}
(\LL^{-\kappa})
\in
\fil^{
\,
-
\card{\Sigma(1)}
+
\sumu{\alpha\in\Sigma(1)}\ea
+
\sumu{\alpha\in A}(\ha-\ea)
}
\,\,
\mlocc.
\end{equation}
Par ailleurs la s\'erie \eqref{eq:sumyconeAJ} s'écrit   comme le produit
$$
\Big(
\sum_{
\substack{
(\ha)\in \N^{A} 
\\ \\ \ha\geq \ea 
}
} 
R_{(\ha)}(T)
\Big)
\times
\Big(
\sum_{y_1\in \cone(I_{A,J,1})}
\,
\LL^{\,\acc{y_1}{\sumu{\alpha\notin A\cup J} \cDa}}
\,
T^{\,\acc{y_1}{\ecL_0}}
\Big)
$$
dont on note $R(T)$ le premier facteur.
Le deuxi\`eme facteur est \'egal \`a
$$
\prod_{i\in I_{A,J,1}}
\,
\left[
\left(
\,1-(\LL\,T)^{\,-\acc{m_i}{\omega}}
\right)^{-1}-1
\right].
$$
Pour tout entier $\kappa\geq 1$,  \eqref{eq:Rha:in:fil} montre que 
$R(\LL^{-\kappa})$
converge dans $\mhatc$ vers un \'el\'ement de 
$$
\fil^{\,-\card{\Sigma(1)}+\sumu{\alpha\in \Sigma(1)}\ea}
\,\,
\mhatc.
$$
Notons enfin que comme $A$ est non vide, 
si $I'=I_\delta$ pour un c\^one $\delta$ maximal alors
$I_{A,J,2}$ ne peut \^etre vide. 
Ceci montre la derni\`ere assertion du lemme. 
\end{proof}
Comme d\'ej\`a annonc\'e, 
on d\'eduit le point \ref{item:1:theoprinc} du th\'eor\`eme \ref{thm:theoprinc}
des propositions \ref{prop:form_phinxs}, 
\ref{prop:Avide_mot} et \ref{prop:Anonvide_mot}.

\providecommand{\bysame}{\leavevmode\hbox to3em{\hrulefill}\thinspace}
\providecommand{\MR}{\relax\ifhmode\unskip\space\fi MR }
\providecommand{\MRhref}[2]{%
  \href{http://www.ams.org/mathscinet-getitem?mr=#1}{#2}
}
\providecommand{\href}[2]{#2}

\end{document}